\documentclass[11pt,reqno]{amsart}
\usepackage[T1]{fontenc}
\usepackage[utf8]{inputenc}
\usepackage{lmodern}
\usepackage[a4paper,textwidth=158mm,textheight=235mm,centering]{geometry}
\usepackage{amsmath,amssymb,mathtools,mathrsfs,booktabs}
\usepackage{longtable,array,enumitem,xcolor}
\usepackage[protrusion=true,expansion=false]{microtype}
\usepackage[hidelinks,pdfusetitle]{hyperref}
\usepackage{bookmark}
\usepackage{tikz}
\usepackage[section]{placeins}
\hypersetup{pdftitle={Noncritical multi-bubbles and hyperbolic double bubbles},pdfsubject={Isoperimetric clusters in constant-curvature spaces},pdfkeywords={multi-bubbles, isoperimetric clusters, profile inequalities, curvature lifting, hyperbolic geometry}}
\numberwithin{equation}{section}
\newtheorem{theorem}{Theorem}[section]
\newtheorem{proposition}[theorem]{Proposition}
\newtheorem{lemma}[theorem]{Lemma}
\newtheorem{corollary}[theorem]{Corollary}
\theoremstyle{definition}

\theoremstyle{remark}
\newtheorem{remark}[theorem]{Remark}
\newcommand{\R}{\mathbb R}
\newcommand{\Sn}{\mathbb S^n}
\newcommand{\Hn}{\mathbb H^n}
\newcommand{\HH}{\mathcal H}
\newcommand{\Id}{\operatorname{Id}}
\newcommand{\tr}{\operatorname{tr}}
\newcommand{\rank}{\operatorname{rank}}
\newcommand{\im}{\operatorname{Im}}

\newcommand{\Vol}{\operatorname{Vol}}
\newcommand{\Per}{\operatorname{Per}}
\newcommand{\supp}{\operatorname{supp}}
\newcommand{\dist}{\operatorname{dist}}
\newcommand{\interior}{\operatorname{int}}
\newcommand{\one}{\mathbf1}
\newcommand{\dd}{\,\mathrm d}
\newcommand{\Istd}{I_{\mathrm{std}}}
\newcommand{\Qzero}{Q^0}
\newcommand{\dcl}{d_{\mathrm{cl}}}
\newcommand{\symdiff}{\mathbin{\triangle}}
\newcommand{\cK}{\mathcal K}
\newcommand{\cM}{\mathcal M}
\newcommand{\cS}{\mathcal S}
\newcommand{\eps}{\varepsilon}

\allowdisplaybreaks[2]
\setlist{itemsep=3pt,topsep=5pt}
\title[Noncritical multi-bubbles and hyperbolic double bubbles]{Noncritical multi-bubbles and hyperbolic double bubbles}

\author{Zihao Wang}
\address{School of Mathematical Sciences, Fudan University, Shanghai 200433, China}
\email{wangzh25@m.fudan.edu.cn}
\date{31 October 2026}

\AtBeginDocument{\hypersetup{pdfauthor={\authors}}}
\begin{document}
\raggedbottom
\begin{abstract}
We prove that standard $k$-bubbles uniquely minimize perimeter in
$\R^n$ and $\Sn$ for $n\ge2$ and $1\le k\le n$, and that standard
double bubbles uniquely minimize perimeter in $\Hn$ for $n\ge3$.
For nonnegative curvature, the proof combines the structure and spectral
theory of Milman--Neeman with a uniform semialgebraic exceptional set,
an almost-everywhere profile comparison, and localized variational
arguments. A curvature-variable lift gives a uniformly elliptic
comparison at the Euclidean limit. The independent hyperbolic argument
uses the rotational structure of Cotton--Freeman, a scalar meridian
force potential, and a transverse conformal variation. The critical
range $k=n+1$ in $\R^n$ and $\Sn$, and hyperbolic clusters with more
than two finite cells, are not covered.
\end{abstract}
\maketitle
\tableofcontents

\section{Introduction and main statements}
\label{sec:intro}
The multi-bubble isoperimetric problem asks how to enclose several
prescribed volumes while minimizing the total area of the separating
interfaces. The single-volume problem selects a round ball in each
simply connected space form. With several volumes, the interfaces can
meet, the chambers can have several components, and the possible
singular junctions become part of the problem. A proof of optimality
must therefore control arbitrary finite-perimeter partitions before
it can identify the familiar standard configuration.

A $q$-cluster is a measurable partition $\Omega=(\Omega_1,\ldots,\Omega_q)$,
identified up to null sets, whose cells have locally finite perimeter and
whose total interface area is finite. We set $k=q-1$. In Euclidean and hyperbolic space,
$\Omega_q$ is the exterior cell and the other cells have prescribed positive
finite volumes. On the sphere all $q$ volumes are prescribed. We count an
interface only once:
\[
 P_g(\Omega)=\frac12\sum_{i=1}^q\Per_g(\Omega_i)
 =\sum_{i<j}\HH_g^{n-1}(\partial^*\Omega_i\cap\partial^*\Omega_j).
\]
No connectedness or regularity condition is imposed on competitors.

A standard spherical cluster is a conformal image of the spherical
regular-simplex partition. A standard Euclidean cluster is a stereographic
image of such a partition, with the projection pole in its exterior cell.
For prescribed positive volumes the standard family has a unique member
modulo isometry; this is uniqueness \emph{within the model family}, not a
classification of arbitrary minimizers; see \cite{Montesinos,MNS,MNP}.

\begin{theorem}[Noncritical spherical multi-bubbles]\label{thm:sphere-main}
Let $n\ge2$ and $2\le q\le n+1$. For every positive vector of $q$ spherical
cell volumes summing to $\Vol(\Sn)$, the standard cluster minimizes total
perimeter among all finite-perimeter clusters with those volumes. Every
minimizer is standard, up to an ambient isometry and null sets, with the
prescribed labels retained.
\end{theorem}

\begin{theorem}[Noncritical Euclidean multi-bubbles]\label{thm:euclidean-main}
Let $n\ge2$ and $1\le k\le n$. For every $v\in(0,\infty)^k$, the standard
Euclidean $k$-bubble minimizes total perimeter among all finite-perimeter
clusters with finite-cell volume vector $v$. Every minimizer is standard,
up to an ambient Euclidean isometry and null sets.
\end{theorem}

A standard hyperbolic double bubble consists of three compact totally umbilical
caps with a common $(n-2)$-dimensional spherical boundary, meeting there
at $120$ degrees. The two outer caps lie on geodesic spheres; the
separating cap may lie on any totally umbilical hypersurface
\cite[Definitions 2.2 and 2.5, Proposition 2.6]{CF}.

\begin{theorem}[Hyperbolic double bubbles]\label{thm:hyperbolic-main}
Let $n\ge3$ and let $v_1,v_2>0$. A standard double bubble in hyperbolic
space of sectional curvature $-1$ minimizes total perimeter among all
finite-perimeter clusters enclosing these two volumes. Every minimizer is
standard, up to an ambient hyperbolic isometry and null sets.
\end{theorem}

\subsection*{Earlier work}
Sullivan's standard-cluster conjecture appears as Problem~2 in the
collection edited by Sullivan and Morgan \cite{SullivanMorgan}.
Montesinos Amilibia \cite{Montesinos} established existence and uniqueness
of a standard Euclidean cluster with prescribed positive volumes.
This identifies the model family against which arbitrary competitors
must be compared. Milman's survey \cite{Survey} describes the subsequent
development of the multi-bubble problem in model spaces.

The first complete results for several prescribed volumes arose in low
dimensions. Foisy--Alfaro--Brock--Hodges--Zimba \cite{Foisy} proved the
planar double-bubble theorem, and Wichiramala \cite{Wichiramala} proved
the planar triple-bubble theorem. Hutchings \cite{Hutchings} developed
symmetry, profile-concavity and component-counting arguments for
minimizing double bubbles. Hass--Schlafly \cite{HassSchlafly} proved the
equal-volume theorem in $\R^3$ with rigorous numerical estimates.
Hutchings--Morgan--Ritor\'e--Ros \cite{HMRR} then proved the arbitrary-volume
theorem in $\R^3$, and Reichardt \cite{Reichardt} established the result
in every Euclidean dimension $n\ge2$.

In Gaussian space, Milman--Neeman \cite{MNGaussian} proved that simplicial
clusters are the unique minimizers for $q$ prescribed positive Gaussian
measures when $2\le q\le n+1$. Their subsequent work \cite{MNS}
established spherical interfaces, a common reflection symmetry and
Voronoi structure for minimizing clusters in $\R^n$ and $\Sn$ in the
range $q\le n+1$. It also settled the double-, triple- and quadruple-bubble
conjectures for $n\ge2$, $n\ge3$ and $n\ge4$, respectively.
In \cite{MNP}, Milman--Neeman proved the spherical quintuple-bubble
theorem for $n\ge5$ and, by approximation, Euclidean perimeter optimality
in the same range. The Euclidean uniqueness question requires a further
argument. That paper develops the Jacobi-operator theory and exact index
theorem used below, proves concavity of the spherical profile, and obtains
higher-bubble results under additional structural hypotheses. Their
boundary-connectedness theorem \cite{Connected} is also used in our
localization arguments. Milman--Xu \cite{MilmanXu} proved stability of
standard $k$-bubbles in all three model spaces for $n\ge3$ and
$2\le k\le n+1$; this concerns the second variation at the standard
configuration.

The hyperbolic problem has been studied through rotational structure
and comparisons between isoperimetric profiles. Cotton--Freeman
\cite{CF} proved the equal-volume double-bubble theorem in $\mathbb H^3$
and established the rotational and topological structure used in
Section~\ref{sec:hyperbolic}. Corneli et al. \cite{Corneli} treated
unequal volumes in $\mathbb H^3$ when the smaller volume is at least
$85\%$ of the larger; they also proved optimality in $\mathbb S^3$ when
each of the three regions occupies strictly more than $10\%$ of the
total volume. Theorem~\ref{thm:hyperbolic-main} addresses arbitrary
positive volumes in every dimension $n\ge3$.

The variational framework comes from geometric measure theory.
Almgren's work \cite{Almgren} established existence and almost-everywhere
regularity for constrained minimizing problems of this kind.
Taylor \cite{Taylor} classified the local singular models for
two-dimensional soap films and bubbles. The regularity theory near
polyhedral cones of Colombo--Edelen--Spolaor \cite{CES}, the local
compactness and convergence results of Cicalese--Leonardi--Maggi
\cite{CLM}, and the finite-perimeter framework in Maggi \cite{Maggi}
are used in Section~\ref{sec:local-tools}. At a density-one tangent
plane we use Allard regularity \cite{Allard}, in the precise formulation
given in Simon \cite{SimonGMT}. The hypotheses of these results are
recorded at their applications.

\subsection*{Contributions and method}
The new step in Sections~\ref{sec:sard}--\ref{sec:profile} is to
organize failure of the linearized structure equations into one closed
semialgebraic null set of pressures, uniformly over all adjacency graphs
and all frameworks. The implication from complete LSE solvability to
the identity $F=F_0$ is already discussed in \cite[Section 9]{MNP};
Proposition~\ref{prop:F} supplies an explicit bounded-field limit proof
adapted to the present argument. The exact index theorem then gives
strictly negative Alexandrov Hessians. A quantitative inverse-gradient
estimate transfers the exceptional pressure set to a null set of
volumes. A tilted lower-contact argument turns the resulting
almost-everywhere inequality into global spherical profile equality.

Generic rigidity by itself does not classify minimizers at exceptional
volumes. We address this by minimizing in an $L^1$ neighborhood of a
compact family of hypothetical nonstandard equality cases. Volume
repair, local regularity, reflection and first contact are established
before the stable-cluster structure theorem is applied. A strong
comparison principle then excludes the localized family. In Euclidean
space one must additionally retain curvature as a parameter: the flat
profile equation alone loses ellipticity. The substitution
$s=|z|^2$, with $z\in\R^{n+2}$, converts the curvature-dependent
inequality into an elliptic comparison that remains nondegenerate at
$s=0$.

For hyperbolic double bubbles, the new mechanism is a face potential
for the conserved meridian force. Its sign follows from a
minimum-radius argument that uses exactly two finite phases. A single
ambient conformal variation expresses the stability index as a weighted
sum of the face potentials. The opposite signs supplied by stability
and the potential argument force every meridian force to vanish;
the nonstandard alternatives in the imported structure theorem are
then excluded.

\subsection*{Proof architecture}
The proofs are logically separated as follows.
\begin{enumerate}[label=(\roman*),leftmargin=*]
\item Sections~\ref{sec:sard}--\ref{sec:profile} use global spherical
minimizers to prove profile equality and generic rigidity. The only
exceptional set is a fixed semialgebraic null set in pressure space.
\item Sections~\ref{sec:local-tools}--\ref{sec:sphere-unique} establish the
local geometric properties required to apply the same spectral theory to
restricted minimizers. A strong comparison then excludes exceptional
spherical equality cases.
\item Sections~\ref{sec:flat}--\ref{sec:euclidean-unique} prove the
Euclidean inequality and its separate equality statement. The curvature
argument uses spherical profile equality and generic rigidity, but not
all-volume spherical rigidity.
\item Section~\ref{sec:hyperbolic} is independent of those arguments.
It uses the rotational double-bubble structure of \cite{CF}, a planar
force potential, and a direct ambient second-variation calculation.
\end{enumerate}
The proofs impose no a priori connectedness, symmetry or smoothness
assumption on competitors. Geometric structure is imported for global
minimizers or established for the restricted minimizers before it is used.
The results cover the following ranges:
\[
\begin{array}{c|c|c}
\text{ambient space}&\text{dimension}&\text{bubble count }k\\\hline
\Sn,\ \R^n& n\ge2&1\le k\le n\\
\Hn& n\ge3&k=2.
\end{array}
\]
On the sphere ``$k$ bubbles'' means $k+1$ prescribed finite cells; there
is no infinite exterior there.

\section{Notation and preliminary results}
\label{sec:inputs}
Throughout Sections~\ref{sec:inputs}--\ref{sec:sphere-unique}, assume
$n\ge2$ and $2\le q\le n+1$, unless a smaller range is explicitly stated.
The induction starts at $q=1$, where $\Omega_1=M$ and $P(\Omega)=0$.
In particular, the imported Jacobi and canonical-operator statements below
are used within the range $q\le n+2$ stipulated in \cite[Section 7]{MNP}.
Put
\[
 m=n-1,\qquad d=q-1,\qquad E=\one^\perp\subset\R^q,
 \qquad \omega_n=\HH^n(\Sn).
\]
Until Section~\ref{sec:flat}, use $\mu=\omega_n^{-1}\HH^n$ and
$\sigma=\omega_n^{-1}\HH^{n-1}$ on the unit sphere. The same normalization
is applied to area and volume: $V_i(\Omega)=\mu(\Omega_i)$ and
$P(\Omega)=\omega_n^{-1}P_g(\Omega)$, where $g$ is the unit round metric
and $P_g$ is the unnormalized perimeter defined in the introduction. Set
\[
 \begin{gathered}
 \Delta_q=\{v\in[0,1]^q:\one\cdot v=1\},\qquad
 \Delta_q^\circ=\{v\in(0,\infty)^q:\one\cdot v=1\},\\
 I(v)=\inf\{P(\Omega):V(\Omega)=v\}.
 \end{gathered}
\]
Here $\Delta_q^\circ$ denotes relative interior in $\{\one\cdot v=1\}$;
for the initial case $q=1$, it consists of the single point $(1)$.
Gradients, Hessians, traces and Lebesgue measure on this simplex refer to its
Euclidean affine span, whose tangent space is $E$. Inverse matrices are
 inverses on $E$. Throughout this article, local semiconcavity means
 semiconcavity with linear modulus: on each sufficiently small convex
 neighborhood $U$ there is $C_U<\infty$ such that
 $u(x)-C_U|x|^2/2$ is concave. At an Alexandrov point $x\in U$ we mean
 the expansion
 \[
  u(x+h)=u(x)+Du(x)\cdot h+\tfrac12h^TD_A^2u(x)h+o(|h|^2).
 \]
 These are interior notions; no Alexandrov expansion at a simplex
 boundary point is assumed.

A normalized spherical Voronoi representation is
\begin{equation}\label{eq:spherical-scores}
 \Omega_i=\{x\in\Sn:c_i\cdot x+\kappa_i<c_j\cdot x+\kappa_j
 \text{ for all }j\ne i\},\qquad \sum_i c_i=0,\quad\sum_i\kappa_i=0.
\end{equation}
Write $c_{ij}=c_i-c_j$, $\kappa_{ij}=\kappa_i-\kappa_j$, and $e_{ij}=e_i-e_j$.
On an actual interface,
\[
 |c_{ij}|^2-\kappa_{ij}^2=1,\qquad
 \nu_{ij}=c_{ij}+\kappa_{ij}x.
\]
Here $\Sigma_{ij}$ denotes the regular interface between cells $i$ and $j$,
and $\nu_{ij}$ points from $\Omega_i$ into $\Omega_j$. We use
\[
 \mathrm{II}_{ij}(v,w)=\langle\nabla_v\nu_{ij},w\rangle,
 \qquad H_{ij}=\operatorname{tr}\mathrm{II}_{ij}=m\kappa_{ij}.
\]
Thus $H_{ij}$ is the sum of the principal curvatures. The pressure multiplier
is $\lambda=m\kappa$, with sign convention $\delta P=\lambda\cdot\delta V$.
Since $\delta V\in E$, adding a multiple of $\one$ to $\lambda$ does not
change this identity; the condition $\lambda\in E$ fixes that freedom. Perpendicularity means that $c_i\cdot N=0$
for a common unit vector $N$; its equator is $\Sn\cap N^\perp$.
A cell is equatorial if
\[
 \Omega_i\cap(\Sn\cap N^\perp)\ne\varnothing.
\]
We spell out the Plateau terminology of
\cite[Definitions 1.4, 3.18 and 3.26; Lemmas 3.15, 3.20 and 3.25]{MNP}.
For $p\in\Sn$, set
\[
 \begin{gathered}
 I_p=\{i:p\in\overline{\Omega_i}\},
 \qquad c_i^\top(p)=c_i-(c_i\cdot p)p,\\
 C_{p,i}=\{y\in T_p\Sn:(c_i-c_j)\cdot y<0
                  \text{ for all }j\in I_p\setminus\{i\}\},\quad i\in I_p,\\
 d_p=\dim\operatorname{span}\{c_i^\top(p)-c_j^\top(p):i,j\in I_p\}.
 \end{gathered}
\]
For the perpendicular clusters with nonempty cells considered here,
Lemma 3.15 of that source also gives
$I_p=\operatorname*{arg\,min}_i(c_i\cdot p+\kappa_i)$.
The cluster is nondegenerate when every incident label survives in the
tangent partition, namely $C_{p,i}\ne\varnothing$ for every $p$ and
$i\in I_p$. Its nonempty interfaces are then relatively open in their
supporting geodesic spheres. A nondegenerate cluster is $3$-Plateau if
\[
 d_p\le3\quad\Longrightarrow\quad |I_p|=d_p+1
 \qquad (p\in\partial\Omega).
\]
Equivalently, at such a point the tangent partition is the regular
simplex partition of $\R^{d_p}$ times $\R^{n-d_p}$.
For $d_p=1,2,3$, its total interface is respectively a hyperplane,
the product of a three-ray $120$-degree junction with $\R^{n-2}$,
or the tetrahedral cone times $\R^{n-3}$.

We specify the scalar-field spaces and signs following
\cite[Sections 5.1--5.5, Definitions 5.1 and 5.4, and Proposition 5.9]{MNP}.
These definitions apply to the regular stationary clusters with bounded
second fundamental forms and Lipschitz face and triple-stratum domains
in Proposition~\ref{ext:physical-approximation} below.
Write $\Sigma^1=\bigcup_{i<j}\Sigma_{ij}$ and
$\sigma_2=\omega_n^{-1}\HH^{n-2}$. We use real scalar fields
$u_{ji}=-u_{ij}$, the direct-sum spaces
\[
 H^1(\Sigma^1)=\bigoplus_{i<j}H^1(\Sigma_{ij}),\qquad
 L^2(\Sigma^1)=\bigoplus_{i<j}L^2(\Sigma_{ij},\sigma),
\]
and the Sobolev trace $\operatorname{Tr}u_{ij}$ on the face boundary.
The form domain is
\[
 \mathcal D_Q=\left\{u\in H^1(\Sigma^1):
 \operatorname{Tr}u_{ij}+\operatorname{Tr}u_{jk}
       +\operatorname{Tr}u_{ki}=0\text{ on every }\Sigma_{ijk}\right\}.
\]
The trace equality is understood in $H^{1/2}$. Here $\Sigma_{ijk}$ is
the triple stratum shared by the three indicated cells.
On it, let $\eta_{ij}$ be the outward unit conormal of $\Sigma_{ij}$,
and use the second fundamental form
$\mathrm{II}_{ij}(v,w)=\langle\nabla_v\nu_{ij},w\rangle$.
The junction identities and boundary coefficient are
\[
 \eta_{ij}=\frac{\nu_{ik}+\nu_{jk}}{\sqrt3},\qquad
 \beta_{ij}=
 \frac{\mathrm{II}_{ik}(\eta_{ik},\eta_{ik})+
       \mathrm{II}_{jk}(\eta_{jk},\eta_{jk})}{\sqrt3}.
\]
They are equations (2.14)--(2.15) of \cite{MNP} with the same outward
normal convention. The index form on the unit sphere is
\[
 \begin{split}
 \Qzero(u,v)=\sum_{i<j}\Bigg[&\int_{\Sigma_{ij}}
 \bigl(\langle\nabla_\Sigma u_{ij},\nabla_\Sigma v_{ij}\rangle
       -(m+|\mathrm{II}_{ij}|^2)u_{ij}v_{ij}\bigr)\dd\sigma\\
 &-\int_{\partial\Sigma_{ij}}\beta_{ij}
       \operatorname{Tr}u_{ij}\operatorname{Tr}v_{ij}\dd\sigma_2\Bigg],
 \qquad \Qzero(u)=\Qzero(u,u).
 \end{split}
\]
The term $m=n-1$ is $\operatorname{Ric}_{\Sn}(\nu,\nu)$.
For the spherical Voronoi faces used below,
$\mathrm{II}_{ij}=\kappa_{ij}\Id$ and
$\beta_{ij}=(\kappa_{ik}+\kappa_{jk})/\sqrt3$ on $\Sigma_{ijk}$.
The scalar first-volume map is
\[
 \mathcal Vu=\delta_uV=\sum_{i<j}\int_{\Sigma_{ij}}u_{ij}\dd\sigma\,e_{ij}.
\]
The operator domain and its sign are fixed by
\[
 \begin{split}
 \mathcal D_{\rm con}=\{u\in\mathcal D_Q:
   &\text{there is }f\in L^2(\Sigma^1)\text{ such that}\\
   &\Qzero(u,v)=-\langle f,v\rangle_{L^2(\Sigma^1)}
                    \text{ for every }v\in\mathcal D_Q\},
 \qquad L_{\rm Jac}u=f.
 \end{split}
\]
The function $f$ is unique because $\mathcal D_Q$ is dense in $L^2$.
By \cite[Definition 5.4, Remark 5.6, Lemma 5.8 and Proposition 5.9]{MNP},
this is equivalently the facewise operator
\[
 (L_{\rm Jac}u)_{ij}=\Delta_{\Sigma_{ij}}u_{ij}
                      +(m+|\mathrm{II}_{ij}|^2)u_{ij}
\]
with $\Delta_\Sigma u\in L^2$, the trace compatibility above, and
the conormal condition
\[
 \partial_{\eta_{ij}}u_{ij}-\beta_{ij}\operatorname{Tr}u_{ij}
 =\partial_{\eta_{jk}}u_{jk}-\beta_{jk}\operatorname{Tr}u_{jk}
 =\partial_{\eta_{ki}}u_{ki}-\beta_{ki}\operatorname{Tr}u_{ki}
 \quad\text{on }\Sigma_{ijk}.
\]
For Sobolev fields these conormal derivatives are the $H^{-1/2}$
functionals defined by Green's formula; for smooth fields they are the
usual outward conormal derivatives. In particular,
$\Qzero(u,v)=-\langle L_{\rm Jac}u,v\rangle$ for
$u\in\mathcal D_{\rm con}$ and $v\in\mathcal D_Q$.
 We use the following results.

\begin{proposition}[Milman--Neeman: global minimizing clusters]
\label{ext:global-structure}
Let $n\ge2$ and $2\le q\le n+1$. For every
$v\in\Delta_q^\circ$, a perimeter-minimizing $q$-cluster on $\Sn$
with $V(\Omega)=v$ exists. After modification on null sets, every such
minimizer is a regular, stable, perpendicular $3$-Plateau spherical
Voronoi cluster whose open cells are nonempty, connected and equatorial.
Here area and volume have the common normalization fixed above.

For each prescribed vector of $q-1$ positive finite cell volumes in
$\R^n$, a global perimeter minimizer also exists. Its total boundary is
bounded; after modification on null sets it is regular, stable and
perpendicularly spherical Voronoi, and its open cells, including the
exterior cell, are connected and equatorial.
The spherical statements are
\cite[Theorems 1.6 and 2.2]{MNP}; the Euclidean statements follow from
\cite[Theorems 1.9, 2.3, 2.4 and 2.15; Lemma 2.6]{MNS}.
\end{proposition}

\begin{proposition}[Milman--Neeman: the spherical profiles]
\label{ext:spherical-profiles}
Let $n\ge2$ and $2\le q\le n+1$. The true spherical profile $I$
is concave on $\Delta_q$. For every $v\in\Delta_q^\circ$, there is
a standard spherical $q$-cluster of volume $v$, unique up to ambient
isometry with its labels retained. Its normalized perimeter defines
$\Istd$. This function is smooth on $\Delta_q^\circ$, is continuous
on $\Delta_q$ with the recursive boundary values obtained by deleting
zero-volume labels, and satisfies
\begin{equation}\label{eq:model-pde}
 D^2\Istd\prec0,\qquad
 \tr\left[(-D^2\Istd)^{-1}
 \left(\Id_E+\frac{2}{m^2}\nabla\Istd\otimes\nabla\Istd\right)\right]
 =\frac2m\Istd.
\end{equation}
The concavity assertion is \cite[Theorem 1.9]{MNP}; the standard-family
statements and the equation are
\cite[Lemmas 3.39--3.40; Proposition 7.25]{MNP}. Smoothness and
boundary continuity are also stated in the proof of
\cite[Proposition 4.2]{MNP}.
\end{proposition}

\begin{proposition}[Milman--Neeman: approximation by ambient fields]
\label{ext:physical-approximation}
Let $n\ge2$ and $2\le q\le n+1$. Let $\Omega$ be a bounded,
regular, stationary $q$-cluster on $\Sn$ with
$V(\Omega)\in\Delta_q^\circ$. Assume that its second fundamental
forms are bounded on the regular interfaces and that both the face
domains $\Sigma_{ij}$ and the triple-stratum domains $\Sigma_{ijk}$
have Lipschitz topological boundaries, as in
\cite[Section 5.1, condition (Lip)]{MNP}. With the Sobolev form domain
$\mathcal D_Q$ and index forms defined above, for every
$u\in\mathcal D_Q$ and every $\eps>0$ there is a smooth compactly
supported ambient field $X_\eps$ such that
\begin{equation}\label{eq:physical-approx}
 \delta_{X_\eps}V=\mathcal Vu,
 \qquad |Q(X_\eps)-\Qzero(u)|\le\eps.
\end{equation}
This is the unit-sphere specialization of
\cite[Theorem 5.18]{MNP}, with the standing assumptions of Sections 5.1
and 5.8 of that source. In particular, the first-volume equality in
\eqref{eq:physical-approx} is exact.
\end{proposition}

\begin{proposition}[Milman--Neeman: the exact negative index]
\label{ext:exact-index}
Let $n\ge2$ and $2\le q\le n+1$, and let $\Omega$ be a regular stable,
perpendicular $3$-Plateau spherical Voronoi $q$-cluster on $\Sn$
with equatorial cells and $V(\Omega)\in\Delta_q^\circ$. Then $(L_{\rm Jac},\mathcal D_{\rm con})$
has exactly $d=q-1$ positive eigenvalues, counted with multiplicity.
Equivalently, the largest dimension of a subspace of $\mathcal D_Q$
on which $\Qzero$ is negative definite is $d$. The sign translation
uses $\Qzero(f,g)=-\langle L_{\rm Jac}f,g\rangle$ on the operator
domain. This is \cite[Theorem 7.22]{MNP}.
\end{proposition}

\begin{proposition}[Milman--Neeman: conformal Jacobi fields and their response]
\label{ext:canonical-response}
Let $n\ge2$ and $2\le q\le n+1$, and let $\Omega$ be a
$3$-Plateau spherical Voronoi $q$-cluster on $\Sn$ with positive
cell volumes. Suppose that for every $a\in E$ there is
$f^a\in\mathcal D_{\rm con}$ satisfying
$L_{\rm Jac}f^a_{ij}=ma_{ij}$ on every actual interface. Choose such
fields linearly in $a$. There is then a canonical symmetric operator
$F:E\to E$, independent of this choice, for which
\begin{equation}\label{eq:canonical}
 \mathcal Vf^a=Fa,\qquad \Qzero(f^a)=-m a^TFa.
\end{equation}
If $\Omega$ is in addition stable and perpendicular with equatorial
cells, the assumed family exists. These assertions are
\cite[Theorem 7.7; Propositions 7.9 and 7.14]{MNP}.
Canonicity in the first part does not require stability of $\Omega$.
\end{proposition}

\begin{proposition}[Milman--Neeman: the relaxed conformal response]
\label{ext:relaxed-response}
Let $n\ge2$ and $2\le q\le n+1$, and let $\Omega$ be a
perpendicular $3$-Plateau spherical Voronoi $q$-cluster on $\Sn$
with positive cell volumes and North pole $N$. Let $\Phi_\tau$ be
the flow of the spherical vector field
$W_N(x)=N-(x\cdot N)x$, and set $\Omega_\tau=\Phi_\tau(\Omega)$.
For $\tau\ne0$ this cluster has a compatibility vector
$\xi_\tau=\coth(\tau)N$ satisfying
$c_i(\tau)\cdot\xi_\tau+\kappa_i(\tau)=0$ for every $i$;
in particular, it is pseudo-conformally flat and its canonical response
$F_\tau$ exists. As $\tau\to0$, in any operator norm on $E$,
\begin{equation}\label{eq:relaxed}
 \begin{gathered}
 F_\tau\longrightarrow F_0
 =n\sum_{i<j}\int_{\Sigma_{ij}}(x\cdot N)^2\dd\sigma\,e_{ij}\otimes e_{ij},\\
 F_0\succ0,\qquad
 \tr\left[F_0\left(\tfrac12\Id_E+\kappa\otimes\kappa\right)\right]=P(\Omega).
 \end{gathered}
\end{equation}
These are \cite[Corollary 8.3; Propositions 7.23 and 8.4--8.5]{MNP}.
The notation $F_0$ denotes the relaxed limit; this statement alone does
not identify it with the canonical response of $\Omega$. Neither
stability nor a minimizing property of $\Omega_\tau$ is assumed.
\end{proposition}

\begin{proposition}[Milman--Neeman: Jacobi fields from the LSE]
\label{ext:lse-jacobi}
Let $n\ge2$ and $2\le q\le n+1$, and let $\Omega$ be an actual
$3$-Plateau spherical Voronoi $q$-cluster on $\Sn$ with
$V(\Omega)\in\Delta_q^\circ$. Let $a\in E$
and let $b_1,\ldots,b_q\in\R^{n+1}$ satisfy $\sum_i b_i=0$ and
$c_{ij}\cdot b_{ij}=\kappa_{ij}a_{ij}$ on every actual interface.
Then the scalar field $f^a_{ij}(x)=a_{ij}+b_{ij}\cdot x$ belongs
to $\mathcal D_{\rm con}$ and satisfies
$L_{\rm Jac}f^a_{ij}=ma_{ij}$ on every actual interface. This is
\cite[Proposition 9.6]{MNP}. The hypothesis concerns a geometric cluster,
not merely a collection of algebraic framework data.
\end{proposition}

\begin{proposition}[Milman--Neeman: rigidity from the canonical response]
\label{ext:response-rigidity}
Let $n\ge2$ and $2\le q\le n+1$. Let $\Omega$ be a global
perimeter-minimizing $q$-cluster on $\Sn$ with volume vector
$v\in\Delta_q^\circ$. Suppose that $I=\Istd$ on a relatively
open neighborhood of $v$ in $\Delta_q^\circ$ and that the canonical
response and relaxed response of $\Omega$ satisfy $F=F_0$.
Then $\Omega$ agrees up to null sets with a standard cluster.
For the prescribed labeled volumes, this identifies it up to ambient
isometry. This is \cite[Proposition 12.2]{MNP}, together with uniqueness
within the standard family in \cite[Lemma 3.40]{MNP}.
\end{proposition}

For the unit-sphere clusters in Proposition~\ref{ext:global-structure}, the curvature and domain
hypotheses in Proposition~\ref{ext:physical-approximation} are available before it is invoked.
Indeed $\mathrm{II}_{ij}=\kappa_{ij}\Id$ on each nonempty interface, so
\[
 \sup_{\Sigma^1}|\mathrm{II}|^2
 \le m\max_{\Sigma_{ij}\ne\varnothing}|\kappa_{ij}|^2<\infty.
\]
The maximum is over finitely many label pairs. The face and triple-stratum
Lipschitz hypotheses follow from the $3$-Plateau Voronoi structure as
explained in \cite[Section 5.1]{MNP}.
The following additional structure and topology results enter the uniqueness arguments.

\begin{proposition}[Milman--Neeman: structure under stability]
\label{ext:stable-structure}
Let $n\ge2$ and $2\le q\le n+1$. Let $\Omega$ be a bounded,
stable, regular $q$-cluster on $\Sn$ with $V(\Omega)\in\Delta_q^\circ$, where regularity has the meaning
of \cite[Definition 2.19]{MNS}. Assume that every cell is invariant
under reflection in a common great sphere, that its total boundary
$\bigcup_i\partial\Omega_i$ is connected, and that each open cell has
finitely many connected components. Then $\Omega$ is a perpendicular
spherical Voronoi cluster, and all its open cells are connected and
equatorial. This is the stable-cluster part of
\cite[Theorem 1.9]{MNS}. The boundedness assumption refers to the total
boundary. The $3$-Plateau property is not part of this conclusion and is
established separately when this proposition is applied below.
\end{proposition}

\begin{proposition}[Milman--Neeman: the planar cases]
\label{ext:planar-bubbles}
Let $n=2$ and $2\le q\le3$. On $\Sn$, for every
$v\in\Delta_q^\circ$, the standard cluster is the unique global
perimeter minimizer up to ambient isometry and null sets. On $\R^2$,
the same assertion holds for every prescribed vector of $q-1$ positive
finite cell areas. This is the planar specialization of
\cite[Theorem 1.13]{MNS}; the conjectures resolved there include
the uniqueness assertion.
\end{proposition}

\begin{proposition}[Milman--Neeman: full interfaces]
\label{ext:full-interfaces}
Let $n\ge2$ and $2\le q\le n+1$. A spherical Voronoi $q$-cluster
on $\Sn$ is standard if and only if every pair interface is nonempty.
Moreover, a global perimeter-minimizing $q$-cluster in $\R^n$ is
standard whenever every pair interface is nonempty. The spherical
statement is \cite[Proposition 3.38]{MNP}; the Euclidean minimizing
statement is \cite[Corollary 1.12]{MNS}.
\end{proposition}

\begin{proposition}[Milman--Neeman: boundary connectedness of a global minimizer]
\label{ext:connected-global}
Let $(M,g)$ be a smooth, connected, simply connected, homogeneous
Riemannian manifold. Suppose that $M$ is compact or has one end.
Let $\Omega=(\Omega_1,\ldots,\Omega_q)$ be a finite-perimeter cluster
that globally minimizes total perimeter at its prescribed cell volumes,
with at most the exterior cell having infinite volume. Choose its open
representatives so that
\[
 \Omega_i=\interior\overline{\Omega_i},\qquad
 \partial\Omega_i=\overline{\partial^*\Omega_i}
 \quad(1\le i\le q).
\]
Then the total boundary
\[
 \Sigma=\bigcup_{i=1}^q\partial\Omega_i
\]
is connected and bounded \cite[Theorem 1.1]{Connected}.
\end{proposition}

\begin{proposition}[Milman--Neeman: separation of a disconnected partition boundary]
\label{ext:connected-partition}
Let $X$ be a connected, simply connected, locally path-connected metric
space, and let $q\ge1$ be an integer. Suppose that
$U_1,\ldots,U_q\subset X$ are nonempty connected open sets satisfying
\[
 U_i\cap U_j=\varnothing\quad(i\ne j),\qquad
 U_i=\interior\overline{U_i},\qquad
 \bigcup_{i=1}^q\overline{U_i}=X.
\]
For $K\subset\{1,\ldots,q\}$ write
$A_K=\bigcup_{i\in K}\overline{U_i}$. If
$\Sigma=\bigcup_{i=1}^q\partial U_i$ is disconnected, then there are
nonempty sets $I,J$ and an index $\ell$ with
\[
 \{1,\ldots,q\}=I\mathbin{\dot\cup}\{\ell\}\mathbin{\dot\cup}J
\]
such that
\[
 X=A_I\mathbin{\dot\cup}U_\ell\mathbin{\dot\cup}A_J,
 \qquad
 \partial U_\ell=\partial A_I\mathbin{\dot\cup}\partial A_J,
 \qquad
 \partial A_I\ne\varnothing,\quad\partial A_J\ne\varnothing.
\]
Here every dotted union is disjoint. This is
\cite[Theorem 1.5]{Connected}.
\end{proposition}

The ambient spaces used here are connected, simply connected, locally
path-connected homogeneous manifolds; the round sphere is compact and
Euclidean space has one end.

For Euclidean minimizing clusters the basic compactness, infiltration
and regularity results are developed in \cite[Chapters 29--30]{Maggi}.
For local perimeter quasiminimizers, the infiltration, almost-everywhere
regularity and local boundary-convergence results are
\cite[Lemma 4.5, p.~2021; Corollary 4.6, p.~2022;
Theorem 4.9, p.~2024]{CLM}. These belong to the arbitrary-dimensional
Euclidean theory of \cite[Section 4, p.~2013]{CLM}.
The minimizing and conical properties of tangent partitions are proved
by rescaling, gluing and monotonicity in Lemma~\ref{lem:regularity}.
The low-codimension regularity and tangent-cone classification are stated
in Propositions~\ref{ext:ces-regularity} and \ref{ext:taylor-cones}.
The smooth-metric extensions and the stationarity bootstrap
are explained in \cite[Remarks 2.16--2.17]{MNS}.

\begin{lemma}[Physical volume surjectivity]\label{lem:volume-onto}
Let $\Omega$ be a finite-perimeter partition of a connected smooth ambient
manifold, with every label of positive measure. Compactly supported smooth
ambient fields have first-volume range $E$; in the presence of an infinite
exterior the equivalent finite-cell range is $\R^k$.
\end{lemma}
\begin{proof}
Suppose $a\in E$ annihilates all first-volume vectors. Then
\[
 0=\sum_i a_i\int_{\Omega_i}\operatorname{div}_gX\,\dd\Vol_g
 \qquad (X\in C_c^\infty(TM)).
\]
The distributional gradient of $\sum_i a_i\mathbf1_{\Omega_i}$ is zero.
A locally integrable function with zero distributional gradient is constant
on a connected manifold. Positivity of all labels implies that all $a_i$
are equal, and $a\in E$ makes them zero. The range is a linear subspace of
a finite-dimensional space, so the vanishing annihilator proves surjectivity.
For an infinite exterior, define its variation as minus the sum of the
finite-cell variations; the same compactly supported identity applies.

At a regular cluster a useful concrete choice is obtained from a spanning
tree of the interface graph \cite[Lemma 2.29]{MNP}. Choose disjoint smooth
interior interface patches and fields with unit transverse flux. Their
volume vectors are the tree incidence vectors $e_i-e_j$, which form a basis
of $E$. This construction will be used for local repairs, but connectedness
of individual cells is not needed for surjectivity.

Every patch center $z\in\Sigma_{ij}$ is chosen inside a genuine two-phase
neighborhood: for some $\varepsilon_z>0$,
\[
 B_{\varepsilon_z}(z)\cap\Omega_h=\varnothing\quad(h\ne i,j),
 \qquad
 \lim_{t\downarrow0}t^{1-n}
   \HH_g^{n-1}(\partial^*\Omega_h\cap B_t(z))=0\quad(h\ne i,j).
\]
The second condition follows from the first. These are the interface
points used in \cite[Lemma 29.13 and Theorem 29.14]{Maggi}, with the
point condition in \href{https://web.ma.utexas.edu/users/maggi/errata.pdf}{Erratum 32}
(21 October 2024). For nearby $L^1$ competitors we retain the same
ambient fields. Continuity of the volume Jacobian then preserves its
invertibility; the quantitative estimates are proved in Lemma~\ref{lem:repair}.
\end{proof}

\begin{lemma}[Boundary values and compactness]\label{lem:profile-boundary}
The true spherical profile is lower semicontinuous on $\Delta_q$. On a face
of the simplex, both profiles reduce to the problem obtained by deleting
zero-volume labels. At a vertex both profiles vanish.
\end{lemma}
\begin{proof}
Let $v_j\to v$ in the closed simplex. If $\liminf_j I(v_j)=+\infty$,
lower semicontinuity is automatic. Otherwise pass to a subsequence along
which $I(v_j)$ tends to its finite lower limit and choose minimizers
$\Omega^j$. Their cell perimeters have uniformly bounded sum and the
ambient sphere has finite volume. BV compactness in a finite family of
smooth charts (the Euclidean statement is \cite[Theorem 12.26]{Maggi}),
applied to the finitely many characteristic functions, gives a further subsequence converging
in $L^1$ cell by cell. A further almost-everywhere convergent subsequence
shows that the limits are characteristic functions of a partition
$\Omega$. Moreover,
\[
 |\mu(\Omega_i^j)-\mu(\Omega_i)|
 \le \|\mathbf1_{\Omega_i^j}-\mathbf1_{\Omega_i}\|_{L^1}\longrightarrow0.
\]
Thus $V(\Omega)=v$, and lower semicontinuity of every cell perimeter yields
$I(v)\le P(\Omega)\le\liminf_jI(v_j)$.

If $v_i=0$, then $\mathbf1_{\Omega_i}=0$ almost everywhere, so its
distributional derivative and perimeter vanish. Deleting such labels or
inserting null labels preserves both admissibility and total perimeter.
This proves the face identity for the true profile. The corresponding
standard-family boundary statement is included in Proposition~\ref{ext:spherical-profiles}.
At a vertex there is only one nonnull label and the perimeter is zero.
In particular, once the two profiles agree on a boundary face, the
inequality $I\le\Istd$, lower semicontinuity, and continuity of $\Istd$
give continuity of $I$ at every point of that face. This is the boundary
continuity needed in the induction below.
\end{proof}

\section{A uniform exceptional set in pressure space}
\label{sec:sard}
We seek a closed null set of pressures $\mathcal B_{n,q}$ with the
following uniform property:
\[
 \kappa\notin\mathcal B_{n,q}
 \quad\Longrightarrow\quad
 \text{the LSE is solvable for every normalized cluster of pressure }\kappa.
\]
To obtain a set independent of the cluster, we work with all normalized
interface frameworks on the fixed label set. Their defining equations
are polynomial, so semialgebraic decomposition and Sard's theorem can
be applied to their pressure projections.

For a graph $J$ on $\{1,\ldots,q\}$ consider the real algebraic set
\[
 \mathcal X_J=\{(C,\kappa):\ \sum_i c_i=0,\ \kappa\in E,
 \ |c_{ij}|^2-\kappa_{ij}^2=1\ (ij\in J)\}.
\]
Here $c_i\in\R^{n+1}$, $C=(c_1,\ldots,c_q)$, and $J$ ranges over
all simple graphs on the fixed label set. Thus $\mathcal X_J$ is a subset
of a fixed finite-dimensional Euclidean vector space. It includes
frameworks which are not clusters; allowing these extra points makes the
exceptional set independent of any choice of a minimizing cluster. Let
\[
 \mathscr C=\{(b_1,\ldots,b_q)\in(\R^{n+1})^q:\sum_i b_i=0\}.
\]
For $B=(b_i)\in\mathscr C$ and $a\in E$, set
$b_{ij}=b_i-b_j$ and $a_{ij}=a_i-a_j$. The two linear maps are
\[
 \begin{aligned}
 A_z:\mathscr C\to\R^J,&\qquad A_zB=(c_{ij}\cdot b_{ij})_{ij\in J},\\
 T_z:E\to\R^J,&\qquad T_za=(\kappa_{ij}a_{ij})_{ij\in J},
 \end{aligned}
 \qquad \pi(C,\kappa)=\kappa.
\]
Complete solvability of the linearized structure equations (LSE) means
$\im T_z\subset\im A_z$.

\begin{proposition}[Coste: semialgebraic decomposition, projection and dimension]
\label{ext:semialgebraic}
Let $A\subset\R^N$ be semialgebraic. Its image under a coordinate
projection is semialgebraic, and it admits a finite partition into smooth
embedded semialgebraic submanifolds, each semialgebraically diffeomorphic
to an open cube. If $A$ is nonempty, its dimension is
the largest dimension of the pieces in such a partition, and
$\dim\overline A=\dim A$.
These are \cite[Theorem 2.3, Corollaries 2.4 and 3.8, and
Propositions 3.15--3.16]{Coste}.
\end{proposition}

\begin{proposition}[Sard: smooth critical values]
\label{ext:sard}
Let $M$ be a second-countable smooth manifold without boundary of dimension $r\ge d\ge1$,
and let $f:M\to\R^d$ be smooth. Then
\[
 \mathcal L^d\bigl(f(\{x\in M:\rank\mathrm df_x<d\})\bigr)=0.
\]
The Euclidean statements are \cite[Theorems 4.1 and 7.2]{Sard};
the formulation on $M$ follows by applying them in a countable smooth atlas.
\end{proposition}

\begin{lemma}[Exceptional pressures]\label{lem:sard}
There is a closed semialgebraic null set $\mathcal B_{n,q}\subset E$ such that
complete LSE solvability holds for every graph and every framework with
pressure $\kappa\notin\mathcal B_{n,q}$.
\end{lemma}
\begin{proof}
Fix $J$ and use fixed bases for the centered $B$-variables and for $E$.
The entries of $A_z$ and $T_z$ are polynomial functions of $z=(C,\kappa)$.
Consequently
\[
 \mathcal Z_J=\{z\in\mathcal X_J:
       \rank[A_z\ T_z]>\rank A_z\}
\]
is semialgebraic: each rank condition is a finite combination of
vanishing and nonvanishing conditions on minors. Take a finite decomposition
$\mathcal Z_J=\bigcup_\alpha S_\alpha$ into smooth semialgebraic submanifolds
\cite[Corollary 3.8]{Coste}. For $z\in S_\alpha$ and a tangent vector $(B,a)$, choose a
smooth curve in $S_\alpha$ with that initial velocity. Differentiating
$|c_{ij}|^2-\kappa_{ij}^2=1$ along the curve gives
\[
       2c_{ij}\cdot b_{ij}-2\kappa_{ij}a_{ij}=0
       \quad(ij\in J),\qquad A_zB=T_za.
\]
If $\mathrm d(\pi|_{S_\alpha})_z$ were onto $E$, there would be such a
vector $(B,a)$ for every $a\in E$, contradicting the definition of
$\mathcal Z_J$. Thus the derivative has rank strictly smaller than $d$
everywhere. Put $s_\alpha=\dim S_\alpha$. If $s_\alpha<d$, a smooth
coordinate parametrization is Lipschitz on each relatively compact
coordinate cube; its pressure image has zero $d$-dimensional measure.
A countable exhaustion handles the whole stratum. If $s_\alpha\ge d$,
every point is critical in the maximal-rank sense of Sard's original
formulation. The required differentiability holds because the map is
smooth. Sard's theorem (Proposition~\ref{ext:sard}), applied in a
countable atlas, again gives $\mathcal L^d(\pi(S_\alpha))=0$. This also
covers noncompact strata. The finite union $\pi(\mathcal Z_J)$
is null and is semialgebraic by the projection theorem
\cite[Theorem 2.3 and Corollary 2.4]{Coste}.

There are $2^{q(q-1)/2}$ possible graphs. Set
\begin{equation}\label{eq:exceptional-set}
 \mathcal B_{n,q}=\overline{\bigcup_J\pi(\mathcal Z_J)}.
\end{equation}
A null semialgebraic subset of $E$ has dimension at most $d-1$.
Indeed, a $d$-dimensional stratum would be open in $E$, contradicting
nullity; dimension is the largest stratum dimension
\cite[Corollary 3.8 and Proposition 3.15]{Coste}.
Semialgebraic closure preserves dimension \cite[Proposition 3.16]{Coste}. Hence this finite union and
its closure remain null. If $\kappa\notin\mathcal B_{n,q}$, no framework
over that pressure belongs to any $\mathcal Z_J$, which proves the
assertion simultaneously for every graph and every framework.
\end{proof}

The following implication is due to Milman--Neeman \cite[Section 9.1]{MNP}.

\begin{proposition}[Bounded representatives identify the limit]\label{prop:F}
For a perpendicular 3-Plateau normalized spherical Voronoi cluster with
positive-volume cells, complete LSE solvability implies $F=F_0\succ0$.
\end{proposition}
\begin{proof}
Choose solutions $b_i$ of the LSE for a basis of $E$ and extend them
linearly in $a$. Subtracting the common mean does not change differences.
Orthogonal projection of each $b_i$ onto $N^\perp$ also preserves the
equations because $c_{ij}\cdot N=0$. We may therefore arrange
$\sum_i b_i=0$ and $b_i\cdot N=0$.

Write $\Phi_\tau$ for the actual perpendicular M\"obius flow. Its
parameter transformation is
\[
 c_i(\tau)=c_i-\kappa_i\sinh\tau\,N,
 \qquad \kappa_i(\tau)=\kappa_i\cosh\tau.
\]
For fixed $a\in E$ put
\[
 b_i(\tau)=\frac{b_i}{\cosh\tau}-a_i\tanh\tau\,N.
\]
Using both orthogonality relations, direct multiplication gives
\[
 c_{ij}(\tau)\cdot b_{ij}(\tau)-\kappa_{ij}(\tau)a_{ij}
 =\frac{c_{ij}\cdot b_{ij}-\kappa_{ij}a_{ij}}{\cosh\tau}=0.
\]
The diffeomorphism $\Phi_\tau$ takes each actual interface onto its
corresponding interface and preserves the $3$-Plateau property. By
Proposition~\ref{ext:lse-jacobi}, $f^a_{ij,\tau}(x)=a_{ij}+b_{ij}(\tau)\cdot x$ is a
conformal Jacobi field. This gives a complete family at every $\tau$,
including $\tau=0$, so its volume response is the canonical operator
of Proposition~\ref{ext:canonical-response}. Stability of the moved cluster is not needed.

To justify the limit on the moving interfaces, pull every integral back
to the original $\Sigma_{ij}$. Let $j_{ij,\tau}$ be the tangential
Jacobian of $\Phi_\tau$ there. Smooth dependence of the ambient
diffeomorphism on the compact sphere implies
\[
 \Phi_\tau\longrightarrow\Id\ \text{in }C^1,\qquad
 \sup_{\Sigma_{ij}}|j_{ij,\tau}-1|\longrightarrow0.
\]
Since $b_{ij}(\tau)\to b_{ij}$ and all these vectors are bounded for
$|\tau|\le\tau_0$, the integrand
\[
 \big(a_{ij}+b_{ij}(\tau)\cdot\Phi_\tau(x)\big)j_{ij,\tau}(x)
\]
converges uniformly to $a_{ij}+b_{ij}\cdot x$ and is bounded by a
constant depending only on the fixed field and $\tau_0$. Each original
interface has finite area. Dominated convergence and the finite number
of interfaces therefore give $F_\tau a\to Fa$. The independent
relaxed-limit formula in Proposition~\ref{ext:relaxed-response} gives $F_\tau a\to F_0a$.
Equality on a basis proves $F=F_0$; in finite dimension the convergence
also holds in operator norm. Positive definiteness follows from
Proposition~\ref{ext:relaxed-response}. No continuity assertion for a degenerating inverse Jacobi
operator is being used.
\end{proof}

\section{From pressure genericity to all-volume optimality}
\label{sec:profile}
Two distinct measure-theoretic steps are needed. First, a null set of
pressures must pull back to a null set of volumes. Concavity alone does
not give this: a gradient can map an open set to a point. The index
theorem supplies the strictly negative Hessian needed for a quantitative
inverse estimate. Second, an inequality valid almost everywhere must
control a contact point. Rather than assume that an arbitrary minimum
is a good point, we construct a positive-measure family of tilted
contacts and choose a good point from that family.

\begin{proposition}[Semiconcavity and second derivatives;
Cannarsa--Sinestrari]\label{ext:semiconcavity}
Let $N\ge1$, let $G\subset\R^N$ be open, and let $u:G\to\R$ be locally semiconcave
with linear modulus. On every convex neighborhood $B\Subset G$ on
which the semiconcavity constant is $C$, the function
$u-C|x|^2/2$ is concave and
\[
 \int_B u\,\partial_{\xi\xi}\varphi\,\dd x
 \le C|\xi|^2\int_B\varphi\,\dd x
 \qquad(\xi\in\R^N,\ 0\le\varphi\in C_c^\infty(B)).
\]
Moreover $u$ is locally Lipschitz, its almost-everywhere gradient belongs
to $BV_{\rm loc}(G;\R^N)$, and for almost every $x\in G$ there are
$Du(x)\in\R^N$ and $D_A^2u(x)\in\operatorname{Sym}_N$ such that
\[
 u(x+h)=u(x)+Du(x)\cdot h+\tfrac12h^TD_A^2u(x)h+o(|h|^2).
\]
These local statements follow from \cite[Proposition 1.1.3,
Theorems 2.1.7 and 2.3.1]{Cannarsa}, applied on convex neighborhoods.
\end{proposition}

\begin{lemma}[Strict Alexandrov Hessians]\label{lem:hessian}
Assume $n\ge2$ and $2\le q\le n+1$. At every Alexandrov point
$v\in\Delta_q^\circ$ of the true spherical profile,
$D_A^2I(v)\prec0$, and every minimizer at $v$ has pressure
$m\kappa=\nabla I(v)$.
\end{lemma}
\begin{proof}
Fix an Alexandrov point $v$ and any minimizer $\Omega$ there. Put
$\lambda=\nabla I(v)$, $A=D_A^2I(v)$ and let $\lambda_\Omega=m\kappa$
be its stationary pressure. For a smooth compactly supported field $X$,
write
\[
 V(\Phi_t\Omega)=v+tw+\tfrac12t^2z+o(t^2),\qquad
 P(\Phi_t\Omega)=I(v)+t\lambda_\Omega\cdot w
                         +\tfrac12t^2P''(0)+o(t^2).
\]
The Alexandrov expansion, evaluated along this curve, gives
\[
 I(V(\Phi_t\Omega))=I(v)+t\lambda\cdot w
       +\tfrac12t^2\big(\lambda\cdot z+w^TAw\big)+o(t^2).
\]
The remainder is $o(t^2)$ because $V(\Phi_t\Omega)-v=O(t)$.
Minimality compares the last two expressions for both signs of $t$.
Their linear terms agree, and Lemma~\ref{lem:volume-onto} gives
$\lambda=\lambda_\Omega$. Comparing quadratic terms now yields
\[
 w^TAw\le P''(0)-\lambda_\Omega\cdot V''(0)=Q(X).
\]
For a fixed $u\in\mathcal D_Q$, use \eqref{eq:physical-approx}:
the approximating fields have exactly the same first-volume vector
$\mathcal Vu$. Letting $\eps\downarrow0$ proves
\begin{equation}\label{eq:hessian-form}
 (\mathcal Vu)^TA(\mathcal Vu)\le\Qzero(u)
           \qquad(u\in\mathcal D_Q).
\end{equation}
Proposition~\ref{ext:exact-index} supplies a $d$-dimensional subspace
$G\subset\mathcal D_Q$ on which $\Qzero$ is negative definite. If
$u\in G$ and $\mathcal Vu=0$, stability and the same physical
approximation imply $\Qzero(u)\ge0$, so $u=0$. Thus
$\mathcal V|_G:G\to E$ is an isomorphism. For every nonzero $w\in E$
choose $u_w\in G$ with $\mathcal Vu_w=w$; then
$w^TAw\le\Qzero(u_w)<0$. The argument holds for the initially arbitrary
minimizer, proving both assertions.
\end{proof}

\begin{lemma}[Pullback of null sets]\label{lem:null-pullback}
Let $u$ be a finite convex function on an open convex subset of $\R^d$.
Among its Alexandrov points with positive-definite Hessian, the inverse
image under $\nabla u$ of any Lebesgue-null set is null.
\end{lemma}
\begin{proof}
At an Alexandrov point $x$ with positive-definite Hessian, its
second-order expansion supplies $c_x,r_x>0$ such that $B_{r_x}(x)$
lies in the domain and
\[
 u(y)\ge u(x)+\nabla u(x)\cdot(y-x)+c_x|y-x|^2
             \quad (|y-x|<r_x).
\]
At any differentiability point $y$ in this ball, convexity gives
$u(x)\ge u(y)+\nabla u(y)\cdot(x-y)$. Adding yields
\[
 (\nabla u(y)-\nabla u(x))\cdot(y-x)\ge c_x|y-x|^2,
 \qquad |\nabla u(y)-\nabla u(x)|\ge c_x|y-x|.
\]
Cover the points in question by countably many pieces on which
$c_x\ge1/j$, $r_x\ge1/\ell$, and which lie in a fixed rational cube
of diameter less than $1/\ell$. The displayed inequality holds for any
two points of each piece. The gradient is injective on that piece, and
its inverse on its image is $j$-Lipschitz.

Let $N$ be a null subset of the target, and fix one piece $S$ with
inverse map $T:\nabla u(S)\to S$ of Lipschitz constant at most $j$.
For $\eps>0$, cover $N$ by cubes $Q_\ell$ with
$\sum_\ell\ell(Q_\ell)^d<\eps$, where $\ell(Q_\ell)$ is the side length.
For each nonempty intersection $Q_\ell\cap N\cap\nabla u(S)$, its
image under $T$ has diameter at most $j\sqrt d\,\ell(Q_\ell)$ and is
contained in a cube of volume at most
$C_dj^d\ell(Q_\ell)^d$. Thus, writing $\mathcal L^{d,*}$ for outer measure,
\[
 \mathcal L^{d,*}\bigl(T(N\cap\nabla u(S))\bigr)
 \le C_dj^d\sum_\ell\ell(Q_\ell)^d<C_dj^d\eps.
\]
Let $\eps\downarrow0$ and take the countable union over the pieces.
The resulting preimage has zero outer measure, proving the assertion. Neither a global inverse for
$\nabla u$ nor Lipschitz continuity of $\nabla u$ is assumed.
\end{proof}

The open-domain Alexandrov theorem is given in
\cite[Theorem 2.3.1(i)]{Cannarsa}; for a finite convex function $u$,
apply its semiconcave formulation to the concave function $-u$.
For completeness, its local use also
reduces to the whole-space formulation in
\cite[Theorem 1.1]{Alexandrov}. If $u$ is finite convex
on an open convex set $U$ and $K\Subset U$ is a closed ball, local
boundedness of the subgradients gives
$\sup\{|p|:x\in K,\ p\in\partial u(x)\}<\infty$. Thus
\[
 \widetilde u(y)=\sup_{x\in K,\ p\in\partial u(x)}
       \{u(x)+p\cdot(y-x)\},\qquad y\in\R^d,
\]
is finite convex on $\R^d$ and equals $u$ on $K$.
Apply Alexandrov's theorem to $\widetilde u$ and exhaust $U$ by
countably many such ball interiors. This proves the local statement
needed for $u=-I$ on $\Delta_q^\circ$. Apply
Lemma~\ref{lem:null-pullback} to $u=-I$ and the target null set
$-m\mathcal B_{n,q}$. Since $Du=-\nabla I$, Lemmas
\ref{lem:hessian} and \ref{lem:null-pullback} imply that almost every
interior volume has $D_A^2I\prec0$ and
$m^{-1}\nabla I\notin\mathcal B_{n,q}$. All minimizers there share that
pressure. Proposition~\ref{prop:F} applies to each of them.

\begin{proposition}[Almost-everywhere sharp inequality]\label{prop:pdi}
At almost every $v\in\Delta_q^\circ$,
\begin{equation}\label{eq:true-pdi}
 D_A^2I\prec0,\qquad
 \tr\left[(-D_A^2I)^{-1}
 \left(\Id_E+\frac2{m^2}\nabla I\otimes\nabla I\right)\right]
 \le\frac2m I.
\end{equation}
\end{proposition}
\begin{proof}
Fix a volume where the Alexandrov expansion and $F=F_0\succ0$ hold,
and write $A=D_A^2I$ and $p=\nabla I=m\kappa$.
For each $a\in E$, Proposition~\ref{ext:canonical-response} gives a field
$f^a$ with $\mathcal Vf^a=Fa$ and $\Qzero(f^a)=-ma^TFa$.
Thus \eqref{eq:hessian-form} gives, for every $a\in E$,
\[
 a^TFAFa\le-ma^TFa.
\]
Congruence by $F^{-1}$, followed by inversion of positive-definite
operators on $E$, yields
\[
 -A\succeq mF^{-1},\qquad (-A)^{-1}\preceq F/m.
\]
Set $H=\Id_E+2p\otimes p/m^2=\Id_E+2\kappa\otimes\kappa\succ0$.
Since $H^{1/2}(F/m-(-A)^{-1})H^{1/2}\succeq0$, its trace is nonnegative.
Using \eqref{eq:relaxed}, we obtain
\[
 \tr\bigl((-A)^{-1}H\bigr)
 \le\frac1m\tr(FH)
 =\frac2m\tr\bigl(F_0(\tfrac12\Id_E+\kappa\otimes\kappa)\bigr)
 =\frac2mP(\Omega)=\frac2m I(v).
\]
\end{proof}

\begin{lemma}[Positive measure of tilted contacts]\label{lem:contacts}
Let $D\subset\R^d$ be compact convex with nonempty interior, and let $w$
be continuous on a convex neighborhood $U$ of $D$, with
$w-C|x|^2/2$ concave on $U$ for some finite $C\ge0$. Suppose every
minimizer on $D$ of $w(x)-a\cdot x$, for $|a|<\rho$, lies in a fixed
$K\Subset D^\circ$. Then the union of these minimizers has positive measure.
\end{lemma}
\begin{proof}
Use the closed ball of slopes $|a|\le\rho/2$. The set of pairs
$(x,a)$ for which $x$ minimizes $w-a\cdot x$ on $D$ is closed in the
compact product $D\times\overline B_{\rho/2}$, hence compact. Its
projection $\Gamma$ lies in $K$ and is compact. At each such $x$ there
is a lower affine support with slope $a$. Semiconcavity gives an upper
quadratic support at $x$; evaluating both supports in directions $h$
and $-h$, then letting $|h|\downarrow0$, shows that their linear parts
agree. In fact
\[
 0\le w(x+h)-w(x)-a\cdot h\le C|h|^2/2,
\]
so $w$ is differentiable at $x$ and $\nabla w(x)=a$.

Take $C\ge0$ to be a semiconcavity constant on a convex neighborhood
of $D$. For $x,y\in\Gamma$, write $h=y-x$ and $p_x=\nabla w(x)$.
If $p_y\ne p_x$, put $e=(p_y-p_x)/|p_y-p_x|$ and
$z=y+|h|e$. When $|h|<\dist(K,\partial D)$, this point belongs to
$D$. The lower support at $y$, upper support at $x$, and lower support
at $x$ give
\[
 w(y)+p_y\cdot(z-y)\le w(z)
 \le w(x)+p_x\cdot(z-x)+\tfrac C2|z-x|^2,
 \qquad w(y)\ge w(x)+p_x\cdot h.
\]
Since $|z-x|\le2|h|$, subtraction yields
$|p_y-p_x|\le2C|h|$. This estimate is automatic if $p_y=p_x$.
Thus $\nabla w$ restricted to $\Gamma$ is Lipschitz on each of
finitely many sufficiently small pieces. Its image contains the entire
ball $\overline B_{\rho/2}$, because every slope has a minimizer.
If $\Gamma$ were null, its image would be null by the Lipschitz image
estimate, a contradiction. The union for $|a|<\rho$ contains
$\Gamma$ and therefore has positive measure.
\end{proof}

\begin{proposition}[All-volume profile equality]\label{prop:profile-equality}
For $q\le n+1$, $I=\Istd$ on $\Delta_q$.
\end{proposition}
\begin{proof}
Proceed by induction on $q$, with the one-cell identity
$I_1(1)=0=I_{{\rm std},1}(1)$ as the base case. In the induction step
$2\le q\le n+1$. Deleting zero-volume labels and applying the inductive
hypothesis gives equality on every boundary face. Lemma
\ref{lem:profile-boundary} and $I\le\Istd$ give continuity at the
boundary; interior continuity follows from the finiteness and concavity
of $I$. Thus $w=I-\Istd$ is continuous on the closed simplex, vanishes
on its boundary, and is nonpositive.

Suppose $\min w=-\eta<0$, attained at $v_*$. The set
$\{w\le-\eta/2\}$ is compact in the interior. Choose a compact
simplex $D\Subset\Delta_q^\circ$ with this set in $D^\circ$, and put
$K=\{w\le-3\eta/4\}\cap D$. The function $w$ is semiconcave near
$D$, since $I$ is concave and $\Istd$ is smooth. Write
$B_*=-D^2\Istd\succ0$ and $p_* =\nabla\Istd$. Choose $\rho>0$ so that
\[
 \rho\operatorname{diam}D<\eta/4,\qquad
 \frac2m\max_K|B_*^{-1}p_*|\rho<\eta/2.
\]
If $x$ minimizes $w-a\cdot v$ on $D$, $|a|<\rho$, comparison with
$v_*$ gives
$w(x)\le-\eta+a\cdot(x-v_*)<-3\eta/4$.
All tilted minimizers therefore lie in $K\Subset D^\circ$.

Lemma~\ref{lem:contacts} gives a positive-measure set of such contacts.
The points where the Alexandrov expansion or \eqref{eq:true-pdi}
fails form a null set. Choose a contact outside that set, with its
supporting slope $a$. At this interior minimum the first and second
terms in the Alexandrov expansion of $w-a\cdot v$ give
\[
 \nabla I=p_*+a,\qquad D_A^2I\succeq D^2\Istd.
\]
Both negative Hessians are positive definite, so inversion reverses
their order and yields $(-D_A^2I)^{-1}\succeq B_*^{-1}$.
Taking traces against the positive matrix in \eqref{eq:true-pdi},
then using \eqref{eq:model-pde}, gives
\begin{align*}
 \frac2mI
 &\ge\tr\left[B_*^{-1}\left(\Id_E+
              \frac2{m^2}(p_*+a)(p_*+a)^T\right)\right]\\
 &=\frac2m\Istd+\frac4{m^2}p_*^TB_*^{-1}a
                         +\frac2{m^2}a^TB_*^{-1}a.
\end{align*}
The last term is nonnegative. Consequently
$w\ge-(2/m)|B_*^{-1}p_*|\rho>-\eta/2$, in contradiction with
$w<-3\eta/4$ at this contact. Hence $w=0$ everywhere.
\end{proof}

\begin{corollary}[Generic spherical uniqueness]\label{cor:generic}
The set
\[
 U_{n,q}=\{v\in\Delta_q^\circ:m^{-1}\nabla\Istd(v)\notin\mathcal B_{n,q}\}
\]
is open, dense and of full measure. Every global minimizer at such a volume
is standard.
\end{corollary}
\begin{proof}
Profile equality identifies every minimizer's pressure with
$m^{-1}\nabla\Istd$. Proposition~\ref{prop:F} and Proposition~\ref{ext:response-rigidity} apply. Since
$D^2\Istd\prec0$, the pressure map is a local diffeomorphism. A countable
cover by inverse-function neighborhoods pulls the closed null exceptional
set back to a relatively closed null set with empty interior.
\end{proof}

\section{Local variational tools}
\label{sec:local-tools}
For the local geometric constructions in Sections~\ref{sec:local-tools}
and \ref{sec:sphere-unique}, assume $n\ge3$.
When $n=2$, the range of Theorems~\ref{thm:sphere-main} and
\ref{thm:euclidean-main} is $2\le q\le3$; their optimality and uniqueness
assertions already follow from Proposition~\ref{ext:planar-bubbles}.
The new bounded-chart geometric construction in Section~\ref{sec:euclidean-local}
will likewise be used only for $n\ge3$. Purely analytic statements retain
their explicitly stated dimension ranges.
We now study fixed-volume minima in an $L^1$ neighborhood. The aim is
to prove volume correction, local perimeter comparison, regularity and
boundary connectedness for these restricted minimizers.
For a compact sphere we use
\[
 \dcl(\Omega,\Lambda)=\sum_{i=1}^q\mu(\Omega_i\symdiff\Lambda_i).
\]
In a fixed Euclidean chart the same notation uses Lebesgue measure.
For a compact family of smooth metrics on a common compact set $K$, there
are constants $0<c\le C<\infty$ such that
\[
 c\le\sqrt{\det g(x)}\le C\quad(x\in K),\qquad
 c\,\dcl^{E}(\Omega,\Lambda)\le\dcl^{g}(\Omega,\Lambda)
       \le C\,\dcl^{E}(\Omega,\Lambda).
\]
Here the superscripts denote the metric used to measure volume.
Thus these distances define the same convergence and their estimates
can be transferred with uniform constants.

\begin{lemma}[Uniform finite-dimensional volume repair]\label{lem:repair}
Let $\cK$ be an $L^1$-compact family of regular clusters with fixed positive
volumes and bounded perimeters. Suppose their films lie in a compact
ambient manifold or in a common compact chart. There are finitely many
smooth flow families and an $L^1$ neighborhood of $\cK$ on which one of
these families repairs any sufficiently small volume error. The parameter
is bounded by a constant times the error, and the volume maps have uniformly
bounded first and second parameter derivatives and inverse derivatives.
The perimeter and $L^1$ costs are bounded by $C P(\Lambda)$ times the
parameter size.
\end{lemma}
\begin{proof}
At each $\Omega\in\cK$, choose $d=q-1$ compactly supported smooth
fields $X_1,\ldots,X_d$ as in Lemma~\ref{lem:volume-onto}, whose
first-volume vectors form a basis. Denote their flows by $\Phi_t^{X_j}$
and define
\[
 F_a=\Phi_{a_d}^{X_d}\circ\cdots\circ\Phi_{a_1}^{X_1},\qquad
 A_\Omega=\bigl(\delta_{X_1}V(\Omega),\ldots,
                    \delta_{X_d}V(\Omega)\bigr).
\]
We use independent volume coordinates, so $A_\Omega$ is invertible.
Each flow fixes every point outside its field's compact support;
therefore $F_a=\Id$ outside the union of these finitely many supports.
For any finite-perimeter cluster,
\[
 V_i(F_a\Lambda)=\int_{\Lambda_i}\operatorname{Jac}_gF_a\dd\mu_g.
\]
Here and below, an exterior cell of infinite volume is handled by taking
the negative sum of the finite-cell increments; no infinite integral is
differentiated. The flows are the identity off one fixed compact set.
For every parameter multi-index $\gamma$ of bounded order,
$\partial_a^\gamma\operatorname{Jac}_g F_a$ is uniformly bounded there.
Differentiation under the integral follows from dominated convergence, and
the resulting integral changes by at most $C\dcl(\Lambda,\Omega)$ when
the cluster is changed. Thus the derivative matrix of the volume map is
continuous in $\Lambda$ in the $L^1$ topology, without any convergence
assumption on its normals or boundaries.

Write $A_\Omega$ for this matrix at $a=0$. Choose a parameter ball so small
that $\|A_\Omega^{-1}(D_aV(F_a\Lambda)-A_\Omega)\|<1/2$ for all
clusters in a neighborhood of $\Omega$. The quantitative inverse function
theorem then solves the prescribed volume increment $e$ with
$|a|\le 2\|A_\Omega^{-1}\||e|$. Differentiating the inverse equation
gives uniform first and second derivative bounds on this inverse.
Explicitly, set $f_\Lambda(a)=V(F_a\Lambda)-V(\Lambda)$, using the
independent finite-volume coordinates. The map
\[
 T_{\Lambda,e}(a)=a+A_\Omega^{-1}\bigl(e-f_\Lambda(a)\bigr)
\]
has derivative norm at most $1/2$ on the chosen parameter ball.
For $|A_\Omega^{-1}e|$ smaller than half its radius it maps that closed
ball into itself. Iteration from $a=0$ is Cauchy, since successive
increments decrease by a factor at most $1/2$; its unique fixed point
satisfies $f_\Lambda(a)=e$ and the stated bound on $|a|$.
The identities $D_ea=(D_af_\Lambda)^{-1}$ and
$D_e^2a[h,k]=-(D_af_\Lambda)^{-1}
D_a^2f_\Lambda[D_ea h,D_ea k]$ give the claimed derivative bounds.
Compactness of $\cK$ selects finitely many such neighborhoods and flow
families. Taking the minimum of their admissible radii and the maximum of
their constants makes every bound uniform. The area formula gives
\[
 P_g(F_a\Lambda)=\sum_{i<j}\int_{\Sigma_{ij}(\Lambda)}
 \operatorname{Jac}_{n-1,g}(DF_a|T\Sigma_{ij})\dd\HH_g^{n-1}.
\]
For a tangent $m$-plane $T$, put
$J_m(a,x,T)=\operatorname{Jac}_{m,g}(DF_a|T)$, where $m=n-1$.
The flow parameter lies in a fixed compact neighborhood, $x$ lies in
the common compact support, and the Grassmannian of $m$-planes is
compact. Smoothness and invertibility of $DF_a$ therefore give
\[
 \sup_{a,x,T}|\partial_a^\gamma J_m(a,x,T)|\le C
 \quad(1\le|\gamma|\le2),\qquad
 |\partial_a^\gamma P_g(F_a\Lambda)|\le C P_g(\Lambda).
\]
For one constituent flow $\Phi_t$ with generator $X$, the $BV$
transport formula and its area formula imply
\[
 \begin{split}
 \Vol_g(\Phi_t E\symdiff E)
 &\le\int_0^{|t|}\int_{\partial^*(\Phi_{\sigma\tau}E)}
             |X\cdot\nu|\,\dd\HH_g^{n-1}\dd\tau\\
 &\le C|t|\Per_g(E),\qquad \sigma=\operatorname{sign}(t).
 \end{split}
\]
This formula follows first for sets with smooth boundary by differentiation along
the flow, and then for $BV$ indicators by strict smooth approximation
in the flow's compact support and lower semicontinuity. Applying it
successively to the finitely many constituent flows gives
\[
 \dcl(F_a\Lambda,\Lambda)\le C|a|P_g(\Lambda).
\]
The same compactness argument makes these constants uniform for a
smooth compact family of metrics in the fixed chart.
\end{proof}

\begin{lemma}[Restricted minima are ordinary local minima]\label{lem:quasimin}
Suppose a restricted volume minimizer has a fixed positive $L^1$ distance
from the boundary of its admissible tube, uniform perimeter bound $M$, and
belongs to the repair neighborhood of Lemma~\ref{lem:repair}. Then it is
stationary and stable for the ordinary volume constraints. Uniformly over
such minimizers there are $\Lambda_0,r_0>0$ such that
\begin{equation}\label{eq:quasimin}
 P_g(\Lambda;B_r(x))\le P_g(\Lambda';B_r(x))
 +\Lambda_0\dcl(\Lambda',\Lambda)
\end{equation}
whenever $r<r_0$ and the two partitions agree off a compact subset of $B_r(x)$.
\end{lemma}
\begin{proof}
Write $d=\dcl(\Lambda',\Lambda)$, and let $\delta>0$ be a
common lower bound for the $L^1$ distance from the admissibility boundary.
If $P_g(\Lambda')\ge P_g(\Lambda)$, cancellation outside $B_r(x)$
proves \eqref{eq:quasimin}. Otherwise $P_g(\Lambda')\le M$.
The volume error $e=V(\Lambda)-V(\Lambda')$ satisfies
\[
 |e|\le C d,\qquad d\le C r^n.
\]
Lemma~\ref{lem:repair} gives a parameter $a$ with
\[
 V(F_a\Lambda')=V(\Lambda),\quad |a|\le C d,\quad
 \begin{cases}
 P_g(F_a\Lambda')\le P_g(\Lambda')+CMd,\\
 \dcl(F_a\Lambda',\Lambda')\le CMd.
 \end{cases}
\]
Choose $r_0$ so small that $(1+CM)Cr_0^n<\delta/2$ and that
both competitors lie in the repair neighborhood. Then
\[
 \dcl(F_a\Lambda',\Lambda)\le(1+CM)d<\delta/2,
 \qquad
 P_g(\Lambda)\le P_g(F_a\Lambda')\le P_g(\Lambda')+CMd.
\]
The two original partitions agree near the complement of $B_r(x)$.
Subtracting their equal perimeter contributions there gives
\eqref{eq:quasimin} with $\Lambda_0=CM$.

First-volume surjectivity and constrained first variation give
$\delta P=\lambda\cdot\delta V$. To verify stability, let $\Phi_t$
have first volume derivative zero, and put
\[
 v(a,t)=V(F_a\Phi_t\Lambda),\qquad
 p(a,t)=P_g(F_a\Phi_t\Lambda),\qquad A=D_av(0,0).
\]
The inverse function theorem gives $a(t)$ satisfying
$v(a(t),t)=v(0,0)$. All following derivatives are evaluated at $(0,0)$:
\[
 a(0)=0,\quad a'(0)=0,\quad Aa''(0)+v_{tt}=0,
 \qquad p_a=\lambda^TA.
\]
Consequently
\[
 \left.\frac{\dd^2}{\dd t^2}p(a(t),t)\right|_{t=0}
 =p_{tt}+p_a a''(0)
 =p_{tt}-\lambda\cdot v_{tt}=Q(X).
\]
For sufficiently small $|t|$, the corrected path stays in the admissible
tube and preserves volume exactly. Its perimeter has a minimum at zero,
so $Q(X)\ge0$.
\end{proof}

\begin{lemma}[Infiltration and regular-open representatives]\label{lem:infiltration}
For a family satisfying \eqref{eq:quasimin}, with a uniform finite number of
labels and uniformly controlled metrics, there are $r_1,\eta>0$ such that
\[
 \Vol_g(\Lambda_i\cap B_r(x))\le\eta r^n
 \quad\Longrightarrow\quad
 \Vol_g(\Lambda_i\cap B_{r/2}(x))=0
 \qquad (r\le r_1).
\]
The cells have regular-open representatives. At their boundary points both
the cell and its complement have volume at least $\eta r^n$ in $B_r(x)$,
after decreasing $\eta$. Component infiltration is established separately
in Lemma~\ref{lem:component-infiltration}, after interface regularity.
\end{lemma}
\begin{proof}
Put $u(t)=\Vol_g(\Lambda_i\cap B_t(x))$. The coarea formula makes
$u$ absolutely continuous. The $BV$ cutting and trace formulas hold for
almost every radius; the finite perimeter measure charges $\partial B_t$
for at most countably many $t$. Choose $t<r_0/2$ outside these exceptional
sets. The two traces then agree almost everywhere on $\partial B_t$,
and the coarea identity and gluing formulas hold at this radius.
For each $\ell\ne i$, let
\[
 A_{i\ell}(t)=\HH_g^{n-1}
   (\partial^*\Lambda_i\cap\partial^*\Lambda_\ell\cap B_t),
 \qquad a_i(t)=\HH_g^{n-1}(\Lambda_i\cap\partial B_t).
\]
Relabel the $i$-part inside $B_t$ by $\ell$, leaving every other label
unchanged, and call the resulting partition $\Gamma$. The trace and
gluing formulas give
\[
 P_g(\Gamma)-P_g(\Lambda)\le -A_{i\ell}(t)+a_i(t),
 \qquad \dcl(\Gamma,\Lambda)=2u(t).
\]
For $t<R<r_0$, the change is compactly contained in $B_R$.
Apply \eqref{eq:quasimin} in $B_R$. The relative perimeters on
$B_R\setminus\overline B_t$ agree and cancel; the possible new jump on
$\partial B_t$ is already included in $a_i(t)$.
It follows that $A_{i\ell}(t)\le a_i(t)+2\Lambda_0u(t)$.
Since
\[
 \Per_g(\Lambda_i;B_t)=\sum_{\ell\ne i}A_{i\ell}(t)
 \le(q-1)\max_{\ell\ne i}A_{i\ell}(t),
\]
and $a_i(t)=u'(t)$ for geodesic balls, we obtain
\[
 \Per_g(\Lambda_i;B_t)
 \le(q-1)\bigl(u'(t)+2\Lambda_0u(t)\bigr)
 \quad\text{for almost every }t.
\]
For coordinate balls write $d(y)=|y-x|$. Uniform metric equivalence
gives $0<c\le|\nabla_g d|_g\le C$ off the center. Coarea reads
\[
 u'(t)=\int_{\Lambda_i\cap\{d=t\}}
                 \frac{1}{|\nabla_g d|_g}\,\dd\HH_g^{n-1},
 \qquad c\,u'(t)\le a_i(t)\le C\,u'(t).
\]
Thus the same differential inequality holds with uniform constants.
Apply the local single-set isoperimetric inequality to $\Lambda_i\cap B_t$:
\[
 c u(t)^{(n-1)/n}\le C_1u'(t)+C_2u(t).
\]
Since $u(t)^{1/n}\le Ct$, the last term is absorbed for $t\le r_1$. If
$u(r/2)>0$, integration of $(u^{1/n})'\ge c'>0$ from $r/2$ to $r$ gives
$u(r)\ge(c'r/2)^n$, proving the contrapositive.

Define
\[
 U_i=\{x:\text{there is }r>0\text{ with }\Vol_g(B_r(x)\setminus\Lambda_i)=0\}.
\]
It is open. At every density-one point of $\Lambda_i$, all the other
labels have volume $o(r^n)$ in $B_r(x)$. Applying infiltration to the
finitely many other labels gives $B_{r/2}(x)\subset U_i$ for small $r$.
The Lebesgue density theorem thus gives $\Vol_g(U_i\symdiff\Lambda_i)=0$.
If $B_r(x)\subset\overline U_i$, the ball cannot meet another open
representative $U_\ell$, since the representatives are pairwise disjoint.
They cover the ambient space up to a null set, so $\Lambda_i$ has full
measure in $B_r(x)$. Hence
\[
 U_i=\interior\overline U_i.
\]
If $x\in\partial U_i$, every $B_{r/2}(x)$ meets $U_i$ in positive
volume. The contrapositive of infiltration gives
$\Vol_g(U_i\cap B_r(x))>\eta r^n$.
For the complement, if
$\Vol_g(B_r(x)\setminus U_i)\le\eta r^n$, each of the other labels
satisfies the infiltration hypothesis and is absent from $B_{r/2}(x)$.
Then $x\in U_i$, a contradiction. Therefore
\[
 \min\{\Vol_g(U_i\cap B_r(x)),\Vol_g(B_r(x)\setminus U_i)\}
 \ge\eta r^n\qquad(x\in\partial U_i).
\]

\end{proof}

\begin{proposition}[Euclidean interior regularity;
Colombo--Edelen--Spolaor]\label{ext:ces-regularity}
Let $m\ge2$, let $U\subset\R^{m+1}$ be open, and let $V$ be a
rectifiable varifold with multiplicity one, and put
$S=\supp\|V\|\cap U$. Assume $\|V\|=\HH^m\llcorner S$ in $U$.
Assume the following in $U$.
\begin{enumerate}[label=(\roman*),leftmargin=*]
\item There is a locally bounded generalized mean curvature $H$ and
no singular first-variation boundary, in the precise sense that
\[
 \delta V(X)=-\int H\cdot X\,\dd\|V\|
 \qquad(X\in C_c^1(U;\R^{m+1})).
\]
\item $V$ has an associated cycle structure: there are countably many
integral $m$-currents $T_j$ without boundary in $U$ whose rectifiable
carriers cover $S$ up to an $\HH^m$-null set and have no carrier outside
$S$ up to that same measure.
\item There are $C,\alpha,r_0>0$ such that for each ball
$B_r(x)\Subset U$, $r<r_0$, and each Lipschitz map
$\phi:B_r(x)\to B_r(x)$ equal to the identity near the boundary,
\[
 \HH^m(S\cap B_r(x))
 \le(1+Cr^\alpha)\HH^m\bigl(\phi(S)\cap B_r(x)\bigr).
\]
Here $\phi$ is extended by the identity outside the ball. The support
is reduced, $S=\supp(\HH^m\llcorner S)\cap U$.
\end{enumerate}
For $0\le j\le m$, let $S^j(V)$ be the points at which no tangent
cone splits off a $(j+1)$-dimensional Euclidean factor, and set
$S^{-1}(V)=\varnothing$. Then, locally in $U$, for some $\beta>0$,
\begin{enumerate}[label=(\alph*),leftmargin=*]
\item $S^m(V)\setminus S^{m-1}(V)$ is a locally finite union of
embedded $C^{1,\beta}$ $m$-manifolds;
\item $S^{m-1}(V)\setminus S^{m-2}(V)$ is a locally finite union of
embedded $C^{1,\beta}$ $(m-1)$-manifolds, near which the support is
locally diffeomorphic to $Y\times\R^{m-1}$, where $Y$ consists of three
half-lines in $\R^2$ meeting at $120$ degrees;
\item $S^{m-2}(V)\setminus S^{m-3}(V)$ is a locally finite union of
embedded $C^{1,\beta}$ $(m-2)$-manifolds, near which the support is
locally diffeomorphic to $T\times\R^{m-2}$, where $T$ is the
two-dimensional tetrahedral cone;
\item for $m\ge3$, $S^{m-3}(V)$ is relatively closed and
$(m-3)$-rectifiable, with locally finite $\HH^{m-3}$ measure; for
$m=2$ the residual set is empty.
\end{enumerate}
This is the local form of \cite[Theorem 3.15, p.~427]{CES}, using
the definitions of cycle structure and almost minimization
\cite[Definitions 3.6 and 3.11, pp.~424--425]{CES}.
\end{proposition}

\begin{lemma}[Regularity needed by the stable structure theorem]
\label{lem:regularity}
An interior restricted minimizer as in Lemma~\ref{lem:quasimin} has the
regularity required in \cite[Definition 2.19]{MNS}: smooth relatively open
interfaces, the regular triple stratum, the $C^{1,\alpha}$ tetrahedral
stratum, a closed residual set of dimension at most $n-4$, local area upper
bounds. Together with Lemma~\ref{lem:component-infiltration}, these give
all the component conditions in that definition. Its tangent partitions are
perimeter minimizing with no volume constraint.
\end{lemma}
\begin{proof}
In the Euclidean case $g=g_E$, comparison \eqref{eq:quasimin} is the
$(2\Lambda_0,r_0)$-minimizing condition of
\cite[equation (1.10), p.~1986]{CLM}, because their distance satisfies
\[
 d_{\rm CLM}(\Lambda,\Gamma)
 =\frac12\sum_i|\Lambda_i\symdiff\Gamma_i|
 =\frac12\dcl(\Lambda,\Gamma).
\]
The source uses $N$ chamber labels and a complementary label $0$;
our label convention does not change the half-sum over all chambers.
The cited cluster statements are formulated for $N\ge2$, corresponding
to $q\ge3$. For $q=2$, the direct infiltration, compactness and tangent
arguments below, together with two-phase Allard regularity, supply the
single-set case.
The source's infiltration and reduced-interface regularity are
\cite[Lemma 4.5, p.~2021; Corollary 4.6, p.~2022]{CLM}.
Its local boundary-convergence theorem assumes a fixed label set,
common quasiminimality constants and $L^1$ convergence on an open set
$A$, and concludes convergence on each $A'\Subset A$
\cite[Theorem 4.9, p.~2024]{CLM}.
All these results hold in Euclidean dimensions $n\ge2$.

For the smooth metric families here, the common chart gives
\[
 c\,\dcl^E\le\dcl^g\le C\,\dcl^E,
 \qquad \Lambda_0\dcl^g\le C\Lambda_0\dcl^E,
\]
with positive constants independent of the cluster and metric parameter.
These estimates control the volume-error term; the perimeter in
\eqref{eq:quasimin} retains its metric $g$.
The smooth weighted Riemannian regularity framework is explained in
\cite[Remark 2.16]{MNS}. We establish its required hypotheses below
by rescaling, gluing, an explicit Euclidean first-variation identity
and support comparisons. Infiltration was already proved in
Lemma~\ref{lem:infiltration}, and the total-boundary convergence needed
later is proved in Lemma~\ref{lem:hausdorff}. The tangent partitions'
minimizing and conical properties are established directly below.
Here are the compactness and comparison details for tangent partitions.
Filling a comparison ball with one label and using
\eqref{eq:quasimin} gives, for almost every small $t$,
\[
 P_g(\Lambda;B_t(x))\le\HH_g^{n-1}(\partial B_t(x))+Ct^n
                    \le C't^{n-1}.
\]
For an arbitrary sufficiently small $r$, choose a good radius
$t\in(r,2r)$. Inclusion of balls and monotonicity of the perimeter
measure give
\[
 P_g(\Lambda;B_r(x))\le P_g(\Lambda;B_t(x))
 \le C't^{n-1}\le 2^{n-1}C'r^{n-1}.
\] In normal coordinates at $x$, rescale by $r_j^{-1}$,
where $r_j\downarrow0$, and denote the partitions and metrics by
$\Lambda^{(j)}$ and $g_j(y)=g(r_jy)$. Then $g_j\to g_E$ smoothly on
compact sets, and the rescaled comparison reads
\[
 P_{g_j}(\Lambda^{(j)};B)
 \le P_{g_j}(\Gamma;B)+\Lambda_0r_j
       \sum_i\Vol_{g_j}(\Gamma_i\symdiff\Lambda_i^{(j)}).
\]
The area bound and $BV$ compactness give a subsequence converging
locally in $L^1$ to a partition $\mathcal C$.
Fix a competitor $\Gamma$ for $\mathcal C$ with compact difference
inside $B_R$. Thus $\Gamma=\mathcal C$ on an outer annulus. For
$\varepsilon>0$, coarea gives a radius $t_j\in(R,R+\varepsilon)$
at which the jump cost for gluing $\Gamma$ inside $B_{t_j}$ to
$\Lambda^{(j)}$ outside is bounded by
\[
 \frac{C}{\varepsilon}\sum_i
   |(\Lambda_i^{(j)}\symdiff\mathcal C_i)
                    \cap(B_{R+\varepsilon}\setminus B_R)|=o(1).
\]
Apply the rescaled comparison to this glued partition and cancel the
unchanged outer annulus. Uniform convergence of $g_j$ and the local
area bound allow replacement of $P_{g_j}$ by $P_E$ with an $o(1)$
error. Lower semicontinuity then gives
\[
 P_E(\mathcal C;B_R)
    \le P_E(\Gamma;B_{R+\varepsilon}).
\]
Choose $R$ with zero interface measure on $\partial B_R$ and let
$\varepsilon\downarrow0$. It follows that $\mathcal C$ is locally
perimeter minimizing. Taking $\Gamma=\mathcal C$ in the same gluing
argument also yields local convergence of perimeter measures at such
radii. No volume correction survives this argument: its coefficient
is $\Lambda_0r_j\to0$.
We next verify the Euclidean varifold and support hypotheses in
conformally flat coordinates.

Let $S=\bigcup_i\partial\Lambda_i$ for the regular-open
representatives. The two-sided density bounds of
Lemma~\ref{lem:infiltration} imply
$\partial\Lambda_i\subset\partial^e\Lambda_i$, where
$\partial^e$ denotes essential boundary. The finite-perimeter
structure theorem gives
\[
 \HH^{n-1}\!\left(S\setminus
        \bigcup_{i<j}(\partial^*\Lambda_i\cap\partial^*\Lambda_j)
            \right)=0.
\]
Thus the multiplicity-one rectifiable varifold on $S$ is precisely the
cluster interface measure. This identification uses density and finite
perimeter, and does not presuppose the higher-stratum regularity to be
proved.

Let $\lambda_i$ be the volume multipliers. Stationarity
under every ambient field implies
\[
 \delta_gP(Y)=\sum_{i<j}(\lambda_i-\lambda_j)
       \int_{\Sigma_{ij}}\langle Y,\nu_{ij}^g\rangle_g\,\dd A_g.
\]
Thus the multiplicity-one boundary varifold has bounded generalized mean
curvature and no singular first-variation boundary. To obtain a uniform bound, use independent volume coordinates and a
correction family with volume derivative $A$ and perimeter derivative
$b=D_aP_g(F_a\Lambda)|_{a=0}$. Stationarity gives
\[
 b=A^T\widehat\lambda,\qquad
 |\widehat\lambda|\le\|A^{-T}\|\,|b|\le CM,
\]
where $\widehat\lambda=(\lambda_i-\lambda_q)_{i<q}$ and
$P_g(\Lambda)\le M$. Thus all pressure differences are uniformly bounded.

All charts needed in this paper are conformally flat: write
$g=\rho^2g_E$ with $\rho>0$ smooth. Testing the last identity with
$Y=\rho^{-m}X$ gives the Euclidean first variation explicitly:
\begin{equation}\label{eq:euclidean-varifold-force}
 \delta_EP(X)=\sum_{i<j}\int_{\Sigma_{ij}}
 \left\langle \rho(\lambda_i-\lambda_j)\nu_{ij}^E
       -m(\nabla_E\log\rho)^\perp,X\right\rangle\dd\HH^{n-1}.
\end{equation}
This verifies bounded Euclidean mean curvature and absence of a boundary
measure in the actual coordinates where the regularity theorem is used.

For completeness, this also makes the conical nature of tangents explicit.
Use a translated and constantly rescaled conformal chart with $x=0$ and
$g(0)=g_E$. The preceding blow-up, gluing and perimeter-convergence
argument applies in this chart as well.
The bounded-mean-curvature monotonicity formula gives a finite limit of
$r^{-m}P_E(\Lambda;B_r(x))$ as $r\downarrow0$
\cite[Chapter 4, Theorem 3.20 and Remark 3.21]{SimonGMT}.
The local perimeter convergence proved above then gives the same constant
density ratio at every radius for the limiting minimizing partition.
The stationary monotonicity identity
\cite[Chapter 4, equation (3.8)]{SimonGMT} forces the radial direction
to be tangent to its film almost everywhere. Its cell-boundary currents
are cycles, so the conical-current criterion
\cite[Chapter 6, Lemma 2.40]{SimonGMT} makes them invariant under dilations.
For every cell $\mathcal C_i$ the zero radial normal trace gives
\[
 \int_{\mathcal C_i}\bigl(n\varphi+y\cdot D\varphi\bigr)\dd y=0
 \qquad(\varphi\in C_c^\infty(\R^n)).
\]
For any such test function, differentiation of
\[
 e^{-nt}\int_{\mathcal C_i}\varphi(e^{-t}y)\dd y
 =\int_{\R^n}\mathbf1_{\mathcal C_i}(e^t y)\varphi(y)\dd y
\]
therefore gives zero. Hence
$\mathbf1_{\mathcal C_i}(e^t\,\cdot)=\mathbf1_{\mathcal C_i}$ almost
everywhere for every $t\in\R$, proving dilation invariance of the labels.
Thus the tangent partitions used below really are minimizing cones.

Here is the support-deformation comparison. Let $S$ be the total support,
and let $\phi$ be Lipschitz, equal to the identity near the boundary of a
coordinate ball $B_r$, and map $B_r$ into itself. There is a partition
$\Lambda'$ equal to $\Lambda$ near the complement of $B_r$, whose boundary
inside the ball lies in $\phi(S)$ and whose relative perimeter is at most
the weighted area of $\phi(S)$. To give the construction explicitly, extend $\phi$ by the identity and
work with locally integral currents in a slightly larger connected
coordinate ball. The extension is the identity outside a compact set,
so it is proper on the current support: the inverse image of a compact
set, intersected with that support, is compact. The facts used here are
Lipschitz push-forward,
commutation with the boundary, the constancy theorem, and the $BV$
representation of top-dimensional normal currents; see
\cite[Chapter 6, equation (2.31), Lemmas 2.36 and 3.12,
Theorem 2.41, and Remark 2.42]{SimonGMT}.
For the integer-multiplicity push-forward statement, one may embed the
chart as $\R^n\times\{0\}\subset\R^{n+1}$ and extend the map by
$(x,t)\mapsto(\phi(x),t)$; its properness on the current support is
unchanged. This puts that use in the positive-codimension convention of
the cited statement.
The push-forward of $[[\Lambda_i]]$ is a top-dimensional integral current,
so it has an integer-valued $BV$ multiplicity $m_i$. The identity
$\partial\phi_\#=\phi_\#\partial$ implies that $Dm_i$ is supported on
$\phi(S)$. The current $\sum_i\phi_\#[[\Lambda_i]]$ has zero boundary
in the larger ball and agrees with its multiplicity-one integration current
off $B_r$. The constancy theorem therefore gives $\sum_i m_i=1$ almost
everywhere, and $m_i=\mathbf1_{\Lambda_i}$
off $B_r$. Define
\[
 E'_i=\{m_i\ge1\}\setminus\bigcup_{h<i}\{m_h\ge1\}.
\]
There is always a positive $m_i$, so these sets partition the chart up to
a null set. The truncation $t\mapsto\min(1,\max(0,t))$ converts an
integer-valued $BV$ function into the indicator of its positive set.
Finite intersections and differences preserve finite perimeter, proving
that $(E'_i)$ is a Caccioppoli partition. It agrees with the original one
near the complement of $B_r$ and has no reduced interface off $\phi(S)$.
At almost every reduced boundary point of a partition exactly two labels
meet. Thus its total perimeter counts that point once, regardless of the
possibly large multiplicities of the currents used in the construction:
\[
 P_g(\Lambda';B_r)\le
 \int_{\phi(S)\cap B_r}\rho^m\,d\HH^{n-1}.
\]
The construction applies directly to the given Lipschitz map $\phi$.

The BV Poincar\'e inequality for the label differences, which are supported
in $B_r$, gives
\[
 \dcl(\Lambda',\Lambda)
 \le Cr\bigl(P_g(\Lambda;B_r)+P_g(\Lambda';B_r)\bigr).
\]
Put $A=P_g(\Lambda;B_r)$ and $A'=P_g(\Lambda';B_r)$. These two
inequalities imply $(1-C\Lambda_0r)A\le(1+C\Lambda_0r)A'$.
Choose $r$ so that $C\Lambda_0r\le1/2$. Moreover,
$\sup_{B_r}\rho^m/\inf_{B_r}\rho^m=1+O(r)$ uniformly on the compact
chart. Combining these facts with the preceding weighted-area estimate
gives
\begin{equation}\label{eq:support-almost-minimal}
 \HH^{n-1}(S\cap B_r)
 \le (1+C'r)\HH^{n-1}(\phi(S)\cap B_r).
\end{equation}
Thus $S$ satisfies the local Almgren almost-minimizing comparison. Finally the currents
$\partial[[\Lambda_i]]$ have zero boundary, and their rectifiable supports
cover $S$ up to $\HH^{n-1}$-null sets. They provide the cycle structure.

The comparison \eqref{eq:support-almost-minimal} is precisely the
area inequality in \cite[Definition 3.11, equation (3.7), p.~425]{CES},
with surface dimension $n-1$ and $\varepsilon(r)=C'r$.
The support condition in part A of that definition follows from the
regular-open representatives and the density bounds above.
All hypotheses of Proposition~\ref{ext:ces-regularity} have now been
established: $\|V\|=\HH^{n-1}\llcorner S$, the cycle structure, bounded
mean curvature without a singular boundary, and
\eqref{eq:support-almost-minimal}. Apply it with $m=n-1\ge2$.
The resulting strata have dimensions $n-1$, $n-2$, $n-3$, and at most
$n-4$. The regular and triple strata are smooth by stationarity, the
junction equations and \cite[Remark 2.17]{MNS}; the tetrahedral stratum
has the stated $C^{1,\alpha}$ regularity.
The local area bound was proved at the start of the argument.
Lemma~\ref{lem:component-infiltration} below establishes the remaining
component conditions in the definition of a regular cluster.
\end{proof}

\begin{lemma}[Infiltration for connected components]\label{lem:component-infiltration}
Each open connected cell component satisfies the infiltration and lower
density conclusions of Lemma~\ref{lem:infiltration}. The constants are
uniform over families with fixed quasiminimality bounds and a bounded
number of original labels. This does not assume component finiteness.
\end{lemma}
\begin{proof}
By the interface regularity just established, the singular boundary has
zero $\HH^{n-1}$ measure, and a regular point locally separates two
different original labels. If $D$ is one open component of $\Lambda_i$, then
$\partial D\subset\partial\Lambda_i$: a small connected ball contained
in $\Lambda_i$ cannot meet both $D$ and a different component. The
$(n-1)$-measure of $\partial D$ is locally finite.
Near a regular point of $\partial\Lambda_i$, the $i$-side is connected;
if that point belongs to $\partial D$, the whole local $i$-side belongs
to $D$. Thus $D$ has almost $C^1$ boundary: its exceptional boundary
is contained in the original singular set and has zero
$\HH^{n-1}$ measure. The Gauss--Green theorem for almost smooth open sets
\cite[Theorem 9.6]{Maggi} gives, in every coordinate neighborhood,
\[
 \int_D\operatorname{div}_gX\,\dd\mu_g
 =\int_{\partial D\setminus S_{\rm sing}}
       \langle X,\nu_D\rangle_g\,\dd\HH_g^{n-1},
 \qquad X\in C_c^1.
\]
Here $S_{\rm sing}$ is the original singular set and
$\HH_g^{n-1}(S_{\rm sing})=0$. The right-hand side defines a locally
finite vector-valued Radon measure, so $D$ has locally finite perimeter;
on the common compact support its perimeter is finite.
At every regular boundary point the local $i$-side is contained in
exactly one component. Hence
\[
 \HH_g^{n-1}\bigl(\partial^*D\cap
                 \partial^*(\Lambda_i\setminus D)\bigr)=0,
 \quad
 \Per_g(D;O)+\Per_g(\Lambda_i\setminus D;O)=\Per_g(\Lambda_i;O)
\]
for every open $O$. Splitting $\Lambda_i$ into these two labels therefore
preserves the relative cluster perimeter on every open set.

Given a competitor $\Gamma$ for the refinement, merge its two new labels
back into label $i$, obtaining $\bar\Gamma$. Coarsening satisfies
$P_g(\bar\Gamma;B)\le P_g(\Gamma;B)$ and
$\dcl(\bar\Gamma,\Lambda)\le\dcl(\Gamma,\widetilde\Lambda)$.
Apply \eqref{eq:quasimin} to $\bar\Gamma$ and use
$P_g(\widetilde\Lambda;B)=P_g(\Lambda;B)$. This proves the same
quasiminimality inequality for the refined partition, with the same
constants; no volume constraint is imposed on either new label. Apply
Lemma~\ref{lem:infiltration} with at most $q+1$ labels. The component is
regular-open: an interior point of its closure belongs to the original
regular-open cell, and a small connected ball in that cell which meets
$D$ lies in $D$. Thus the lower density bound applies at every component
boundary point.

We make explicit the dimension-only lower density required by
\cite[Lemma 2.6 and Definition 2.19]{MNS}. In all applications here
$q\le n+1$, so the refined partition has at most $n+2$ labels.
At a fixed point $p$, choose normal coordinates and decrease their
radius until both the metric and its volume density are uniformly
comparable to the Euclidean ones, with comparison factors depending
only on $n$. The local isoperimetric constant in the preceding
infiltration argument then depends only on $n$. Decreasing the radius
again absorbs the term containing $\Lambda_0u(t)$. The resulting
constant $c'$ in $(u^{1/n})'\ge c'$ depends only on $n$; the permitted
radius may still depend on the point and cluster. Consequently there
is $\eta_n>0$ such that, for $p\in\partial D$ and all sufficiently
small $r>0$,
\[
 \Vol_g(D\cap B_r(p))\ge\eta_n r^n.
\]
Writing $\alpha_n=\mathcal L^n(B_1(0))$, normal coordinates also give
$\Vol_g(B_r(p))/r^n\to\alpha_n$. Thus, with
$\varepsilon_n=\min\{1/2,\eta_n/(2\alpha_n)\}>0$,
\[
 \Theta(D,p):=\liminf_{r\downarrow0}
 \frac{\Vol_g(D\cap B_r(p))}{\Vol_g(B_r(p))}
 \ge\varepsilon_n\qquad(p\in\overline D).
\]
Interior points have density $1$, so they satisfy the same estimate.
This is the component-density condition in that definition.
\end{proof}

\begin{lemma}[Uniformly finite components]\label{lem:finite-components}
On a compact round sphere the local families above have a uniformly
bounded number of open connected cell components. The same conclusion
holds in a common bounded chart, allowing one exterior component outside
the common support ball.
\end{lemma}
\begin{proof}
Every nonempty component $D$ on the sphere has a boundary point, since all
original cell volumes are positive. For a fixed small radius $r_1$,
Lemma~\ref{lem:component-infiltration} gives
\[
 \Vol_g(D)\ge\Vol_g(D\cap B_{r_1}(x))\ge\eta r_1^n
 \qquad(x\in\partial D).
\]
Thus every finite collection $D_1,\ldots,D_N$ of distinct components
on the sphere satisfies
\[
 N\eta r_1^n\le\sum_{a=1}^N\Vol_g(D_a)\le\Vol_g(M).
\]
The bound holds for every finite subcollection and therefore rules out
infinitely many components. In a bounded chart, all finite components
and all bounded exterior components lie in a fixed support ball $B_R$.
The complement of that ball is connected, so it belongs to at most one
additional exterior component. Consequently
\[
 N_{\rm comp}\le1+\frac{\sup_g\Vol_g(B_R)}{\eta r_1^n},
\]
and the numerator is finite by uniform metric equivalence.
\end{proof}

\begin{proposition}[Taylor's classification of minimizing cones]
\label{ext:taylor-cones}
Let $C\subset\R^3$ be a nonempty reduced two-dimensional cone with
locally finite $\HH^2$ measure. For each $R>0$ put
$S_R=C\cap B_R(0)$ and $B_0=C\cap\partial B_R(0)$.
Assume that $S_R$ is area minimizing with respect to $B_0$ in the
sense of \cite[Section I(8), p.~496]{Taylor}: it is relatively closed
outside $B_0$, has finite area, agrees with the support of its area
measure outside $B_0$, and for every Lipschitz map
$\phi:\R^3\to\R^3$, with $W=\{\phi\ne\Id\}$, such that
$W\cup\phi(W)$ is bounded and has positive distance from $B_0$, one has
\[
 \HH^2(S_R\cap W)\le\HH^2\bigl(\phi(S_R\cap W)\bigr).
\]
Then, up to an orthogonal transformation, $C$ is a plane, the union of
three half-planes meeting at $120$ degrees along their common boundary
line, or the cone over the one-skeleton of a regular tetrahedron centered
at the origin. This is \cite[Proposition II.3, pp.~500--501]{Taylor},
applied to the finite-area restricted cones $S_R$.
\end{proposition}

\begin{lemma}[The Plateau property from minimizing tangents]
\label{lem:local-plateau}
Let $\Omega$ be a regular spherical Voronoi cluster whose tangent
partitions are locally perimeter minimizing without volume constraints.
Then $\Omega$ is $3$-Plateau.
\end{lemma}
\begin{proof}
Regularity includes interface regularity and the lower cell-density
bound. Thus \cite[Lemma 3.19]{MNP} makes $\Omega$ non-degenerate.
For $p\in\partial\Omega$ write
\[
 I_p=\{i:p\in\overline\Omega_i\},\qquad
 n_i=c_i-(c_i\cdot p)p,\qquad
 W_p=\operatorname{span}\{n_i-n_j:i,j\in I_p\},\quad d_p=\dim W_p.
\]
The tangent description in \cite[Lemma 3.24]{MNP} gives
\[
 \mathcal C_p\simeq\mathcal C'_p\times\R^{n-d_p},
\]
where $\mathcal C'_p$ is a conical Voronoi partition in $W_p$ and
every label in $I_p$ has nonempty interior. The factor
$\mathcal C'_p$ is itself perimeter minimizing. To verify this, let $h=n-d_p>0$ and suppose that a competing
partition $\mathcal F'$ agrees with $\mathcal C'_p$ off a compact subset
of $U\Subset W_p$, with
\[
 P(\mathcal F';U)-P(\mathcal C'_p;U)=-\delta<0.
\]
Use $\mathcal F'$ on $U\times(-L,L)^h$ and the original product
partition outside this cylinder. The side boundary of $U$ produces no
jump. On the $2h$ end faces the total jump cost is at most
$2h(2L)^{h-1}|U|$. Therefore the perimeter difference is at most
\[
 -\delta(2L)^h+2h(2L)^{h-1}|U|<0
 \qquad\text{if }L>h|U|/\delta,
\]
contradicting minimality of $\mathcal C_p$. If $d_p=n$, the factor is
the entire tangent partition and there is nothing to prove.

Non-degeneracy gives $d_p\ge1$. We spell out the reduction to the
precise area-minimizing class in \cite[Section I(8), p.~496]{Taylor}.
At a point with $d_p\le3$, put $d=d_p$ and
$\widehat{\mathcal C}=\mathcal C'_p\times\R^{3-d}$.
For a competing partition $\mathcal F$ with compact difference, choose
a product cylinder $U\times I$ containing that difference, where
$U\subset\R^d$ and $I\subset\R^{3-d}$. For each cell $F_i$, the
$BV$ slicing identity and $|D_x\mathbf1_{F_i}|\le|D\mathbf1_{F_i}|$
give
\[
 \int_I\Per((F_i)_y;U)\,\dd y
 =|D_x\mathbf1_{F_i}|(U\times I)
 \le\Per(F_i;U\times I).
\]
Summing over $i$ and dividing by two yields
\[
 P(\mathcal F;U\times I)
 \ge\int_I P(\mathcal F_y;U)\dd y
 \ge |I|P(\mathcal C'_p;U)
 =P(\widehat{\mathcal C};U\times I).
\]
Here $\mathcal F_y$ agrees with $\mathcal C'_p$ outside a compact subset
of $U$ for almost every $y$; for $d=3$ no slicing is needed.
Thus $\widehat{\mathcal C}$ is a minimizing partition in $\R^3$.
Write $\widehat S$ for the support of its total perimeter measure.
Repeating the Lipschitz
push-forward partition construction in Lemma~\ref{lem:regularity}, now
with zero comparison error, gives for every compactly supported
Lipschitz deformation $\phi$, with $W=\{\phi\ne\Id\}$,
\[
 \HH^2(\widehat S\cap B)
 \le\HH^2(\phi(\widehat S)\cap B)
 \le\HH^2(\widehat S\cap(B\setminus W))
      +\HH^2(\phi(\widehat S\cap W)),
\]
where the comparison ball $B$ contains $W\cup\phi(W)$ compactly.
The last inequality is subadditivity of $\HH^2$. Since $W\subset B$,
the disjoint decomposition of the domain gives
\[
 \HH^2(\widehat S\cap B)
 =\HH^2(\widehat S\cap W)
  +\HH^2(\widehat S\cap(B\setminus W)).
\]
The second summand is finite. Subtracting it gives
\[
 \HH^2(\widehat S\cap W)
 \le\HH^2\bigl(\phi(\widehat S\cap W)\bigr),
\]
which is Taylor's moving-part comparison with zero error.

To meet the finite-total-area convention of that definition, fix $R>0$
and use $S_R=\widehat S\cap\{|x|<R\}$ with boundary set
$B_0=\widehat S\cap\partial B_R$. For any Taylor-admissible map $\phi$
let $W=\{\phi\ne\Id\}$ and, if $W\ne\varnothing$, choose
$0<\epsilon<\min\{\dist(W,B_0),R/4\}$. Choose a smooth radial cutoff
$\chi$ equal to $1$ on $[0,R-\epsilon]$ and to $0$ on
$[R-\epsilon/2,\infty)$, and set
\[
 F(x)=x+\chi(|x|)(\phi(x)-x).
\]
For $x\in\widehat S$ with $R-\epsilon<|x|<R$, its radial projection
onto $\partial B_R$ lies in $B_0$, so $\phi(x)=x$.
Consequently $F=\phi$ on $S_R$, $F=\Id$ on $\widehat S\setminus B_R$,
and the moving part of $\widehat S$ under $F$ is exactly $S_R\cap W$.
The preceding zero-error comparison for the compactly supported map $F$
therefore proves
\[
 \HH^2(S_R\cap W)\le\HH^2(\phi(S_R\cap W)).
\]
The case $W=\varnothing$ is immediate. The reduced conical support
$S_R$ is relatively closed outside $B_0$ and has finite area.
It is thus area minimizing in Taylor's sense, for every deformation
scale, rather than merely stationary or stable.

Fix $i_0\in I_p$. By the definition of $W_p$,
\[
 d_p=\operatorname{rank}\{n_i-n_{i_0}:i\in I_p\setminus\{i_0\}\}
 \le |I_p|-1.
\]
For $d_p\le3$, Taylor's cone classification (Proposition~\ref{ext:taylor-cones}), as used in
\cite[Lemma 3.32]{MNP}, applies to the finite-area restricted cone $S_R$.
Conical invariance then identifies $\widehat S$ as a plane,
a $Y$ cone times a line, or a tetrahedral $T$ cone.
Write $j=1,2,3$, respectively, for the dimension of the corresponding
facet-normal span. Every facet normal lies in $W_p\times\{0\}$, so
$j\le d_p$. Each complementary chamber is
connected and has zero partition perimeter in its interior, so it carries
one label. The chamber adjacency graphs of the plane, $Y$ and $T$ films
are $K_2$, $K_3$ and $K_4$, respectively. Two adjacent chambers cannot
carry the same label: their common relatively open face would then have
zero partition perimeter, contrary to its presence in the tangent film.
Consequently the map from the $j+1$ complementary chambers to their
labels is injective. It is also onto $I_p$, because every tangent cell
with a label in $I_p$ has nonempty interior. Hence
\[
 j\le d_p\le |I_p|-1=(j+1)-1=j.
\]
All inequalities are equalities. In particular,
\[
 |I_p|=d_p+1.
\]
The vectors $n_i$, after subtracting their common mean, are consequently
affinely independent. By \cite[Lemma 3.25]{MNP}, this is exactly the
Plateau condition at $p$. Since it holds at every point with
$d_p\le3$, Definition 3.26 of that source gives the $3$-Plateau
property. 
\end{proof}

\begin{lemma}[Strict convergence preserves pair areas]\label{lem:pair-areas}
If $\Omega_j\to\Omega$ in $L^1$ and $P(\Omega_j)\to P(\Omega)$, then every
pair-interface area converges. The assertion holds also for smoothly varying
metrics on a common compact chart.
\end{lemma}
\begin{proof}
Lower semicontinuity of every cell perimeter and convergence of their sum
imply convergence of each cell perimeter. For $i\ne l$,
\[
 2A_{il}=\Per(\Omega_i)+\Per(\Omega_l)-\Per(\Omega_i\cup\Omega_l).
\]
Lower semicontinuity of the last term implies $\limsup_j A_{il}^{(j)}\le
A_{il}$. Summing over the finitely many pairs gives convergence of each.
For clarity, if $a_{il}^{(j)}$ denotes a pair area, the argument has
shown $\limsup_j a_{il}^{(j)}\le a_{il}$ for every pair. Since
$\sum_{i<l}a_{il}^{(j)}\to\sum_{i<l}a_{il}$ and the number of pairs is
finite, a strict deficit along any subsequence would force a strict excess
for another pair. Hence every $a_{il}^{(j)}$ converges. On a common compact
chart, uniform convergence of the metrics gives
$(1-\varepsilon_j)P_g(E)\le P_{g_j}(E)\le(1+\varepsilon_j)P_g(E)$
for every finite-perimeter set, with $\varepsilon_j\to0$. The uniform
perimeter bound therefore reduces varying metrics to the fixed-metric
argument.
\end{proof}

\begin{lemma}[Boundary convergence and a gap estimate]\label{lem:hausdorff}
Let $\Omega_j$ and $\Omega$ be regular-open representatives of clusters
with the same finite label set, whose films lie in a common compact
ambient region. If the $\Omega_j$ satisfy uniform quasiminimality bounds
and converge cell by cell to $\Omega$ in $L^1$, their total boundaries
converge in Hausdorff distance.
If $S_j=A_j\mathbin{\dot\cup}B_j$ is disconnected, with $A_j,B_j$ nonempty
compact, and $S_j\to S$ with $S$ connected compact, then
\[
 \dist(A_j,B_j)\le2d_H(S_j,S)\longrightarrow0.
\]
\end{lemma}
\begin{proof}
For the first Hausdorff inclusion, suppose $x_j\in S_j$ and
$\operatorname{dist}(x_j,S)\ge\varepsilon>0$. Pass to a convergent
subsequence in the common compact region. A fixed ball about its limit
lies in one chamber of the limiting cluster. Uniform boundary density at
$x_j$ requires at least a fixed positive volume of another label in that
ball, contradicting $L^1$ convergence. For the other inclusion, if
$x_j\in S$ and $\operatorname{dist}(x_j,S_j)\ge\varepsilon$, a fixed
ball about a subsequential limit contains no boundary of the $j$th
cluster. Connectedness of the ball and the regular-open representatives
make it a single label. Select a subsequence with that label fixed.
The limiting partition then has the same label almost everywhere in the
ball, contradicting that its center lies on $S$. Compactness converts the
two sequential statements into $d_H(S_j,S)\to0$. For the second
claim, the closed $\eps$ neighborhoods of $A_j,B_j$, with
$\eps=d_H(S_j,S)$, would be disjoint if their distance were greater than
$2\eps$. They cover $S$ and both meet it, contradicting connectedness.
\end{proof}

\begin{proposition}[Allard regularity at a density-one point]
\label{ext:allard}
Let $m\ge1$. Let $V=v(M,\theta)$ be a
rectifiable $m$-varifold in an open set $U\subset\R^{m+1}$, write
$\mu_V=\|V\|$, and assume that $\theta\ge1$ almost everywhere for
$\mu_V$. Suppose the first variation in $U$ is represented by a generalized mean
curvature $H\in L^{2m}_{\rm loc}(\mu_V)$ without a singular boundary:
\[
 \delta V(X)=-\int H\cdot X\,\dd\mu_V
 \qquad(X\in C_c^1(U;\R^{m+1})).
\]
Let $x_0\in\supp\mu_V\cap U$ have density one:
\[
 \lim_{r\downarrow0}\frac{\mu_V(B_r(x_0))}{\omega_m^{\rm ball}r^m}=1,
\]
where $\omega_m^{\rm ball}$ is the volume of the unit $m$-ball.
Then there are $\rho>0$, an orthogonal map $Q$ of $\R^{m+1}$,
an open neighborhood $W\subset\R^m$ of $0$, and
$u\in C^{1,1/2}(W;\R)$ such that $u(0)=0$, $Du(0)=0$, and
\[
 Q(\supp\mu_V-x_0)\cap B_\rho(0)
 =\operatorname{graph}u\cap B_\rho(0).
\]
This is \cite[Chapter 5, Theorem 6.1, pp.~150--151]{SimonGMT}
with codimension one and $p=2m$, a consequence of the small-mass-ratio
and small-mean-curvature theorem
\cite[Chapter 5, hypotheses 5.1 and Theorem 5.2]{SimonGMT}.
\end{proposition}

\begin{lemma}[A supporting ball gives a regular contact]\label{lem:support-plane}
Let a cluster in $\R^n$ or a round $n$-sphere, $n\ge3$, be locally perimeter
quasiminimizing and stationary for perimeter with constant volume
multipliers. Suppose its boundary near $x$ lies in a closed region
admitting a smooth supporting hypersurface at $x$, and that this
hypersurface bounds an exterior ball. Then $x$ is a smooth two-phase point.
\end{lemma}
\begin{proof}
Every tangent boundary cone lies in the corresponding tangent half-space.
It is stationary and is the multiplicity-one film of a minimizing tangent
partition. Let $e$ be the inward normal to this half-space, and take a smooth
nonincreasing radial cutoff $\chi$. Stationarity of $X(y)=\chi(|y|)e$ gives
\[
 0=\int\frac{\chi'(|y|)}{|y|}\,e\cdot y\dd\|T\|(y),
\]
Indeed, if $T_y$ denotes the approximate tangent plane, conical
invariance gives $P_{T_y}y=y$ almost everywhere, and
\[
 \operatorname{div}_{T_y}X
 =\frac{\chi'(|y|)}{|y|}\,e\cdot P_{T_y}y
 =\frac{\chi'(|y|)}{|y|}\,e\cdot y\le0.
\]
Choose $\chi'<0$ on an arbitrary annulus $a<|y|<b$.
The zero integral forces $e\cdot y=0$ for $\|T\|$-almost every point
of that annulus. A countable collection of such annuli covers
$\R^n\setminus\{0\}$; since the film is the support of $\|T\|$,
it is contained in $e^\perp$. The two open half-spaces contain no interface, so each has one label.
They have different labels: uniform boundary density at $x$ excludes an
empty tangent film or a one-label tangent partition. Compactness for
quasiminimizers includes local convergence of perimeter measures, and
therefore the tangent film is this plane with multiplicity one.
Use conformal coordinates centered at $x$ and normalized by $g(0)=g_E$.
The Euclidean multiplicity-one film varifold $V$ has no first-variation
boundary and has locally bounded Euclidean mean curvature $H_E$, by
\eqref{eq:euclidean-varifold-force}. The planar tangent, perimeter
convergence and the density-ratio limit proved above give
\[
 \frac{\|V\|(B_r(0))}{\omega_m^{\rm ball}r^m}\longrightarrow1,
 \qquad
 \left(r^{p-m}\int_{B_r(0)}|H_E|^p\dd\|V\|\right)^{1/p}
 \le Cr\longrightarrow0,
 \quad m=n-1,\quad p=2m>m,
\]
where $\omega_m^{\rm ball}$ denotes the volume of the Euclidean unit
$m$-ball. The density is $1$ almost everywhere. Proposition~\ref{ext:allard} applies at the density-one point. Thus all the smallness,
density and mean-curvature hypotheses of Allard's theorem in
\cite[Chapter 5, hypotheses 5.1 and Theorem 5.2]{SimonGMT} hold after
choosing $r$ small. The entire support in a smaller ball is a single
$C^{1,1/2}$ graph. In a still smaller cylinder its complement has two
connected sides, each carrying one label; the density bounds force
these labels to be different. Stationarity with constant volume
multipliers gives the prescribed-mean-curvature equation, which is
uniformly elliptic on the graph because its gradient is locally bounded.
Its smooth coefficients give the smooth bootstrap
\cite[Remark 2.17]{MNS}.

\end{proof}

\begin{figure}[htbp]
\centering
\begin{tikzpicture}[x=1.65cm,y=1.3cm,>=stealth,font=\small]
  \fill[blue!6] (-1.6,1.45) -- (1.6,1.45) -- plot[domain=1.6:-1.6,samples=65] (\x,{.45*\x*\x}) -- cycle;
  \fill[blue!6] (-1.6,-1.45) -- (1.6,-1.45) -- plot[domain=1.6:-1.6,samples=65] (\x,{-.45*\x*\x}) -- cycle;
  \draw[thick,blue!55!black,domain=-1.6:1.6,samples=65] plot (\x,{.45*\x*\x});
  \draw[thick,blue!55!black,domain=-1.6:1.6,samples=65] plot (\x,{-.45*\x*\x});
  \draw[dashed,gray] (-1.8,0)--(1.85,0) node[right,black] {$T_x\Sigma$};
  \fill (0,0) circle (1.3pt) node[above left] {$x$};
  \node at (0,1.12) {$B$};
  \node at (0,-1.12) {$A_{\delta_j}$};
  \node[anchor=west] at (1.62,1.10) {$u_+$};
  \node[anchor=west] at (1.62,-1.10) {$u_-$};
  \draw[<->] (1.12,-.5645)--(1.12,.5645);
  \node[fill=white,inner sep=1pt] at (.82,.29) {$U_\ell$};
  \node[anchor=east,align=right] at (-1.85,.55) {$u_+(0)=u_-(0)=0$\\$Du_+(0)=Du_-(0)=0$};
\end{tikzpicture}
\caption{A local section at first contact in Lemma~\ref{lem:collision}.
The curves represent two smooth graphs over a common tangent plane;
their quadratic shape is schematic. The gap belongs to the label of
$U_\ell$. Tangency gives $\sup_{|x'|\le r}(u_+-u_-)=o(r)$, so its
volume is $o(r^n)$, contrary to the positive boundary-density bound.}
\label{fig:collision-gap}
\end{figure}
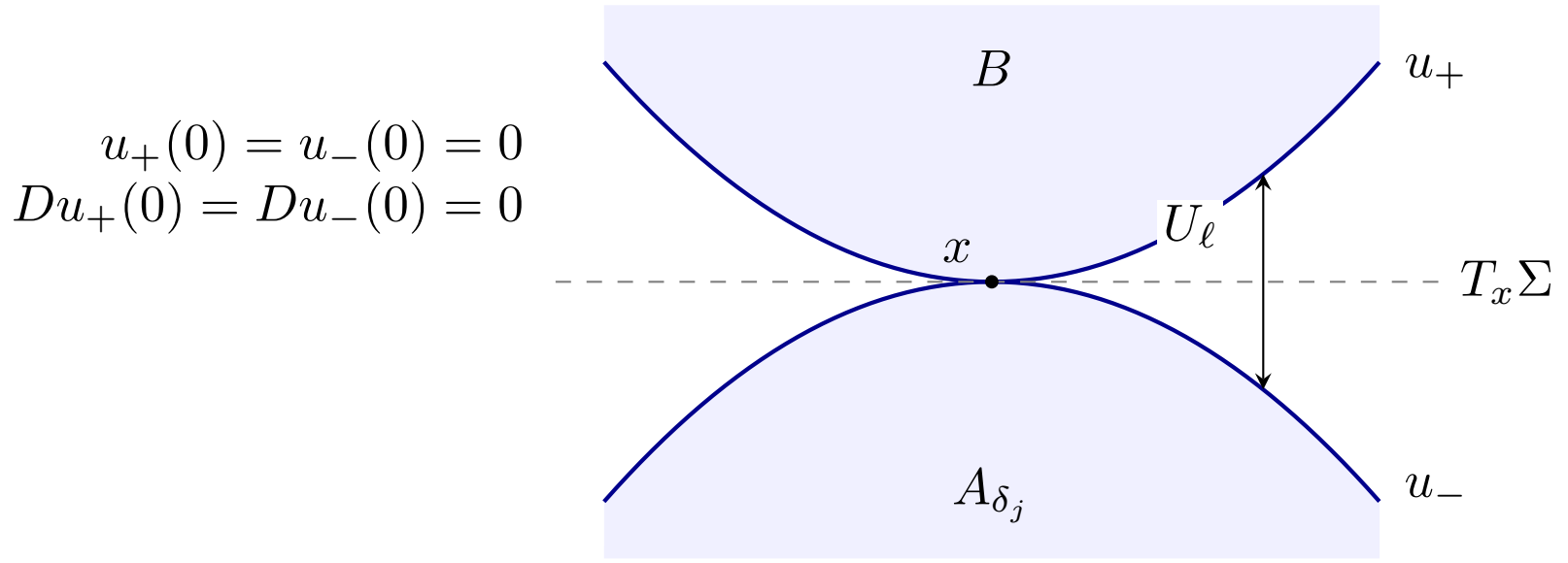

\begin{lemma}[Local collision lemma on round spheres]\label{lem:collision}
Let $\Lambda_j$ attain the perimeter minimum at their prescribed volumes
inside their admissible $L^1$ tubes on round spheres, and suppose they lie
in the interiors of these tubes. Assume uniform local correction, density
and perimeter bounds. The word ``restricted'' here means a minimum over
the specified tube, not merely stationarity or a local minimum with an
unspecified radius. Suppose their
films converge in Hausdorff distance to a connected compact film, and their distance from the
boundary of the admissible $L^1$ tube is bounded below. Then their total
boundaries are connected for all sufficiently large $j$. The same assertion
holds in bounded charts of spheres whose curvatures tend to zero, provided
exact isometric recentering is admissible and has displacement tending to
zero.
\end{lemma}
\begin{proof}
Suppose $S_j=\partial\Lambda_j$ is disconnected. Refine into the finitely
many connected open components provided by Lemma~\ref{lem:finite-components}.
The components are regular-open: this follows from the same property of
their original cells and connectedness of a sufficiently small ball inside
the interior of the component's closure. The refinement creates no
positive-area interface, so it has the same perimeter and boundary.
The topological partition theorem (Proposition~\ref{ext:connected-partition}) gives a
partition of the refined labels into nonempty sets $I,\{\ell\},J$ with
\[
 M=A\mathbin{\dot\cup}U_\ell\mathbin{\dot\cup}B,
 \quad A=\bigcup_{\alpha\in I}\overline{U_\alpha},
 \quad B=\bigcup_{\beta\in J}\overline{U_\beta},
 \quad\partial U_\ell=\partial A\mathbin{\dot\cup}\partial B.
\]
The boundary separates as $(S_j\cap A)\mathbin{\dot\cup}(S_j\cap B)$.
A closest pair of points of $A,B$ lies on their boundaries, so their distance
$\delta_j$ equals the distance between these two boundary pieces. Lemma
\ref{lem:hausdorff} gives $\delta_j\to0$.

Let $x\in\partial A,y\in\partial B$ be a closest pair, and write
$R=s^{-1/2}$ for the sphere radius (with $R=1$ on the unit sphere).
For large $j$, $0<\delta_j<\pi R$. Let $R_t$ be the rotation in the
ambient plane spanned by $x,y$, through angle $t/R$, in the direction
from $x$ to $y$, acting identically on its orthogonal complement.
Its generator has norm at most one on the sphere. Therefore
\[
 R_{\delta_j}x=y,\quad d_g(R_tz,z)\le t,\quad
 d_g(R_ta,b)\ge d_g(a,b)-t\ge\delta_j-t
 \quad(a\in A,\ b\in B).
\]
It follows that $R_tA\cap B=\varnothing$ for $0\le t<\delta_j$
and that first contact occurs at $t=\delta_j$.
Move all refined
cells in $I$, leave those in $J$ fixed, and define $U_\ell(t)$ as the
complement of the two closed groups. All refined volumes, and hence all
original volumes, are exactly preserved. Before collision every interface is either carried by an isometry or
left fixed, and the two moved groups are disjoint; hence the perimeter is
unchanged. The transport estimate below gives $L^1$ continuity up to the
first-contact time. Lower semicontinuity consequently gives perimeter at
collision at most the pre-collision value. Boundary intersections have
zero $n$-volume, so the refined volumes remain unchanged at that time too.

For an isometric motion whose generator has speed at most one,
\[
 \Vol(R_tD\symdiff D)\le t\Per(D).
\]
It follows that the distance of the moved original cluster from the initial
one is at most $4tP(\Lambda_j)$. This tends to zero uniformly along the
motion, so the whole path stays in the tube. In a bounded chart, choose the
group not containing the omitted point; uniform metric equivalence gives
the same estimate there. An isometric recentering with displacement tending
to zero does not affect area or volumes and preserves the tube margin.

Restricted minimality forces equality of perimeter at collision, making
the collision cluster another interior local minimizer. To check regularity before collision, the open ball of radius
$\delta_j$ centered at $y$ misses $A$. In a sufficiently small
neighborhood of $x$, the other group $B$ is absent and $U_\ell$ is the
only region outside $A$. Thus the \emph{entire} original film in that
neighborhood is on the $A$ side of the supporting sphere, including any
interfaces internal to $A$. Lemma~\ref{lem:support-plane} applies to the
original film and makes $x$ a regular two-phase point. The argument at
$y$ is identical. Their tangent planes are perpendicular to the minimizing
geodesic joining them. The chosen rotation carries the tangent plane at
$x$ to that at $y$, so the two regular sheets meet tangentially. The original label containing $U_\ell$
is locally the gap between them: at either regular point the other side
has a different original label. If the sheets overlap in positive
$(n-1)$-measure, perimeter has dropped. Otherwise write the sheets locally as graphs $u_-(x')\le u_+(x')$
over their common tangent plane, as in Figure~\ref{fig:collision-gap}. Both have value and derivative zero at
contact, so
$\sup_{|x'|\le r}(u_+(x')-u_-(x'))=o(r)$.
If the gap is absent on a neighborhood, the sheets coincide on a
positive-area set, the preceding alternative. Hence the gap label is
present arbitrarily close to contact, placing that point on its boundary
in the regular-open representative of the collision minimizer. Its volume
in a radius-$r$ ball is bounded by
$C r^{n-1}\sup_{|x'|\le Cr}(u_+-u_-)=o(r^n)$.
There are no other pieces of that label in a sufficiently small
neighborhood, because each pre-collision contact is a two-phase point.
This contradicts the uniform positive lower boundary density of the
collision minimizer.
Both possibilities are impossible.
\end{proof}

\section{Promotion to all-volume spherical uniqueness}
\label{sec:sphere-unique}
Fix $v_0\in\Delta_q^\circ$. Let $\cM$ be the set of global minimizers at
$v_0$ and $\cS_0$ the standard isometry orbit. Suppose
$\cK=\cM\setminus\cS_0$ is nonempty.

\begin{lemma}[Isolation and the restricted profile]\label{lem:sphere-tube}
The set $\cK$ is compact and positively separated from $\cS_0$. There is
$r>0$ and a neighborhood of $v_0$ such that
\[
 \dist_{\dcl}(\cK,\cS(v))>3r.
\]
For the tube $\mathcal T_r=\{\Lambda:\dist_{\dcl}(\Lambda,\cK)\le r\}$,
set
\[
 J(v)=\inf\{P(\Lambda):V(\Lambda)=v,\ \Lambda\in\mathcal T_r\}.
\]
For nearby $v$, this infimum is attained; all minimizers approach $\cK$
uniformly as $v\to v_0$, have a uniform perimeter bound, and lie, for
example, within distance $r/8$ of $\cK$.
\end{lemma}
\begin{proof}
BV compactness and profile equality make $\cM$ compact, with strict
perimeter convergence along its convergent sequences. If nonstandard
minimizers approached a standard one, Lemma~\ref{lem:pair-areas} would give
positive area for every pair interface for large indices. The global
structure theorem and the full-interface criterion (Proposition~\ref{ext:full-interfaces}) would make them standard. Thus $\cS_0$ is
relatively open in $\cM$, and $\cK$ is compact. The standard family depends
continuously on positive volumes, giving the asserted separation after
shrinking the volume neighborhood.

Use the uniform repair families on $\cK$ to construct competitors at $v$
with area at most $I(v_0)+C|v-v_0|$. At fixed $v$, a minimizing
sequence therefore has uniformly bounded perimeter. BV compactness on
the compact sphere gives a subsequence converging in cluster $L^1$.
The volume constraints pass to this limit. The tube constraint also
passes to the limit, since
\[
 |\dist_{\dcl}(\Lambda,\cK)-\dist_{\dcl}(\Gamma,\cK)|
 \le\dcl(\Lambda,\Gamma).
\]
Lower semicontinuity then proves attainment.

If $v_j\to v_0$ and $\Lambda_j$ are restricted minimizers, the same
perimeter bound gives a subsequential limit $\Omega$ with volume $v_0$,
and
\[
 I(v_0)\le P(\Omega)\le\liminf_jP(\Lambda_j)
 \le\limsup_jP(\Lambda_j)\le I(v_0).
\]
Thus $\Omega\in\cM$ and the perimeters converge. The closed tube
excludes $\cS_0$, so $\Omega\in\cK$. If the asserted uniform approach
failed, some $\varepsilon>0$ and such a sequence would satisfy
$\dist_{\dcl}(\Lambda_j,\cK)\ge\varepsilon$ for every $j$,
contradicting its subsequential convergence to $\Omega\in\cK$.
Taking $\varepsilon=r/8$ gives the stated interior margin.
\end{proof}

Lemmas~\ref{lem:quasimin}--\ref{lem:regularity} now give uniform ordinary
quasiminimality, stationarity, stability and regularity for these restricted
minimizers.

\begin{lemma}[Semiconcavity of $J$]\label{lem:local-semiconcave}
On a sufficiently small convex volume ball $U$ about $v_0$, the function
$J$ is finite, locally Lipschitz and semiconcave.
\end{lemma}
\begin{proof}
For each restricted minimizer $\Lambda$ at $v$, the uniform correction map
produces a $C^2$ competitor value $H_\Lambda(w)$ on a volume ball of a
uniform radius, with
\[
 J(w)\le H_\Lambda(w),\quad J(v)=H_\Lambda(v),\quad
 |DH_\Lambda|+|D^2H_\Lambda|\le C.
\]
All these competitors remain in the tube by the uniform margin. Shrink $U$
so that its diameter is smaller than the common radius. Then
$J=\inf_\Lambda H_\Lambda$ on $U$. The functions
$H_\Lambda-C|w|^2/2$ are concave, and so is their infimum. Two-point
comparison gives the Lipschitz bound.
\end{proof}

\begin{lemma}[Common bisectors of Voronoi clusters]\label{lem:bisector}
If a great hyperplane $N^\perp$ bisects every cell of a positive-cell
spherical Voronoi cluster, its reflection preserves each label.
\end{lemma}
\begin{proof}
Set $a_i=c_i\cdot N$ and $s=x\cdot N>0$. If the winning labels at $x$ and
its reflection $Rx$ are $i,j$, then
\[
 0\le f_j(x)-f_i(x)\le2s(a_j-a_i).
\]
Thus $a_j-a_i\ge0$. Its integral over the hemisphere is zero, because every
label occupies equal volume on both sides. Hence it vanishes almost
everywhere. The two scores at $x$ are equal, so uniqueness of the winning
label away from the measure-zero ties gives $i=j$. This is the labeled
reflection symmetry. The reflection preserves each oriented interface
sphere: it preserves the two adjacent labels and hence its chosen
normal. Uniqueness of the unit-normalized oriented sphere equation
implies $Rc_{ij}=c_{ij}$ on every actual interface. Thus
$c_{ij}\cdot N=0$ there. The interface graph is connected, so all
$c_i\cdot N$ are equal; $\sum_i c_i=0$ makes this common value zero.
Consequently $CN=0$ as well.
\end{proof}

\begin{lemma}[A symmetric structured restricted minimizer]\label{lem:structured}
After shrinking $U$, at every $v\in U$ there is a restricted minimizer
which is a regular stable perpendicular 3-Plateau spherical Voronoi cluster
with equatorial cells.
\end{lemma}
\begin{proof}
A common great-sphere volume bisector exists for $q\le n+1$: only $q-1$
of the bisection equations are independent, and the ham-sandwich theorem
applies; this is the bisection input recalled in \cite[Section 1.2, p.~78]{MNS}.
To see that reflection is admissible, let $\Lambda_j\to\cK$ and
choose any such bisecting normals $N_j$. Pass to $\Lambda_j\to\Omega\in\cK$
and $N_j\to N$. Let $H_j^\pm=\{\pm x\cdot N_j>0\}$ and let $R_j$ be the reflection.
The indicators of $H_j^\pm$ converge almost everywhere to those of
$H^\pm$, because the limiting equator has zero volume. Consequently
\[
 \mu(\Lambda_{j,i}\cap H_j^+)
 \longrightarrow\mu(\Omega_i\cap H^+).
\]
Lemma~\ref{lem:bisector} makes $R\Omega=\Omega$ label by label.
For the moving-reflection limit, write $\chi_i=\mathbf1_{\Omega_i}$
and $\chi_{j,i}=\mathbf1_{\Lambda_{j,i}}$; all norms use the normalization
$\mu(\Sn)=1$. Measure preservation gives
\[
 \|\chi_{j,i}\circ R_j-\chi_i\circ R\|_1
 \le\|\chi_{j,i}-\chi_i\|_1
       +\|\chi_i\circ R_j-\chi_i\circ R\|_1.
\]
For any $\eps>0$, choose a continuous $f_i$ on the sphere with
$\|f_i-\chi_i\|_1<\eps$. Then
\[
 \|\chi_i\circ R_j-\chi_i\circ R\|_1
 \le2\eps+\|f_i\circ R_j-f_i\circ R\|_\infty,
 \qquad
 \limsup_{j\to\infty}\|\chi_i\circ R_j-\chi_i\circ R\|_1\le2\eps.
\]
Uniform continuity of $f_i$ and uniform convergence $R_j\to R$ justify
this bound. Letting $\eps\downarrow0$ proves $L^1$ convergence of the
reflected cells. The identity
\[
 \mathbf1_{\Lambda_{j,i}^+}
 =\mathbf1_{H_j^+}\chi_{j,i}
    +\mathbf1_{H_j^-}(\chi_{j,i}\circ R_j)
 \quad\text{almost everywhere}
\]
then proves $\Lambda_j^+\to\Omega$, and the same calculation gives
$\Lambda_j^-\to\Omega$. If the tube margin failed uniformly near
$v_0$, a sequence witnessing that failure would have a subsequence as
above, contradicting these limits.

For the perimeter comparison, duplicate the restriction of $\Lambda$ to
$H^+$ by reflection, and denote the result by $\Lambda^+$. The $BV$
traces of the two copies agree on the equator, so the gluing formula gives
$P(\Lambda^+)=2P(\Lambda;H^+)$; similarly,
$P(\Lambda^-)=2P(\Lambda;H^-)$.
In particular their sum is at most $2P(\Lambda)$, even if the original
film has positive area on the equator. Bisection gives
$V(\Lambda^\pm)=V(\Lambda)$. Choose the copy with smaller perimeter.
Restricted minimality forces equality. We have obtained a symmetric
restricted minimizer.

Every restricted minimizer has finitely many components by Lemma
\ref{lem:finite-components}. Their films converge, along subsequences as
$v\to v_0$, to the connected film of a global minimizer in $\cK$; use
Lemma~\ref{lem:hausdorff} and Proposition~\ref{ext:connected-global}. Lemma
\ref{lem:collision} then gives connected total boundary for all nearby
restricted minimizers. Proposition~\ref{ext:stable-structure} applies to the symmetric one: boundedness,
regularity, stability, finite components and connected boundary have all
been checked. It gives the perpendicular spherical Voronoi representation
and equatorial cells.

Finally the tangent partitions are minimizing by
Lemma~\ref{lem:regularity}. Lemma~\ref{lem:local-plateau} therefore
supplies the $3$-Plateau property. At this stage the faces lie on
finitely many fixed geodesic spheres, so their curvatures are bounded.
The face and triple-stratum Lipschitz assumptions needed in
Proposition~\ref{ext:physical-approximation} follow from \cite[Section 5.1]{MNP}. In particular the
spectral theory is invoked only after these properties have been obtained.
\end{proof}

\begin{proposition}[Sharp inequality for the restricted profile]
\label{prop:local-pdi}
On a smaller $U$, the function $J$ is concave. Almost everywhere there,
\[
 D_A^2J\prec0,\qquad
 \tr\left[(-D_A^2J)^{-1}
 \left(\Id_E+\frac2{m^2}\nabla J\otimes\nabla J\right)\right]
 \le\frac2mJ.
\]
\end{proposition}
\begin{proof}
At an Alexandrov point of the semiconcave $J$, choose a structured restricted
minimizer. All sufficiently small physical variations remain admissible.
The proof of Lemma~\ref{lem:hessian}, including exact-volume physical
approximation and the $d$-dimensional negative space, gives
$D_A^2J\prec0$. A semiconcave function has a matrix-valued distributional
Hessian whose singular part is negative semidefinite
\cite[Proposition 1.1.3 and Theorem 2.3.1]{Cannarsa}.
We record why its absolutely continuous density equals the Alexandrov
Hessian. At a point $x$ with an Alexandrov expansion, which is also a
Lebesgue point of that density and a zero-density point of the total
variation of the singular part, set
\[
 J_r(y)=\frac{J(x+ry)-J(x)-rDJ(x)\cdot y}{r^2}.
\]
Then $J_r$ converges locally uniformly to
$y^TD_A^2J(x)y/2$. Consequently its distributional Hessians converge to
$D_A^2J(x)\,\mathcal L^d$. On the other hand, differentiation of the
Hessian measure shows that these same rescaled measures converge to
its absolutely continuous density at $x$, times $\mathcal L^d$.
These limits agree. The stated three pointwise properties hold almost
everywhere, which proves the identification.
Thus the entire distributional Hessian of $J$ is negative semidefinite,
and $J$ is concave. Lemma~\ref{lem:null-pullback} applies to $-J$ with
target null set $-m\mathcal B_{n,q}$, and therefore
$m^{-1}\nabla J(v)\notin\mathcal B_{n,q}$ for almost every $v\in U$.
The proof of Proposition
\ref{prop:pdi} now applies.
\end{proof}

\begin{proposition}[Weak Harnack inequality with zero data;
Gilbarg--Trudinger]\label{ext:weak-harnack}
Let $N\ge3$ and $R>0$, and let $B_{4R}(y)\subset\R^N$. Let $A$ be a measurable
symmetric matrix and let $\beta,c$ be measurable coefficients on that
ball. Assume that for some $0<\lambda\le\Lambda<\infty$,
\[
 \begin{gathered}
 \lambda|\xi|^2\le\xi^TA(x)\xi\le\Lambda|\xi|^2
 \quad\text{for every }\xi\in\R^N\text{ and almost every }x,\\
 \beta,c\in L^\infty(B_{4R}(y)),\qquad c\le0.
 \end{gathered}
\]
Suppose $u\in W^{1,2}(B_{4R}(y))\cap C(B_{4R}(y))$ is nonnegative
and is a weak supersolution of
$\operatorname{div}(A Du)+\beta\cdot Du+cu\le0$, meaning that
\[
 \int A Du\cdot D\varphi\,\dd x
 -\int(\beta\cdot Du+cu)\varphi\,\dd x\ge0
 \qquad(0\le\varphi\in C_c^1(B_{4R}(y))).
\]
Then
\[
 \frac1{|B_{2R}|}\int_{B_{2R}(y)}u\,\dd x
 \le C_0\operatorname*{ess\,inf}_{B_R(y)}u,
\]
where $C_0$ depends only on $N$, $\Lambda/\lambda$,
$R\|\beta\|_\infty/\lambda$ and $R^2\|c\|_\infty/\lambda$.
This is the $p=1$, zero-right-hand-data case of
\cite[Theorem 8.18, p.~194]{GT}. In that theorem's notation one may
take the data exponent $Q=N+1>N$; zero data belong to both required
spaces $L^Q$ and $L^{Q/2}$, and $1<N/(N-2)$.
\end{proposition}

\begin{lemma}[Strong minimum principle for semiconcave functions]
\label{lem:weak-minimum}
Let $G\subset\R^N$ be connected and open, where $N\ge1$. Let
$A\in C^1(G;\operatorname{Sym}_N)$ be positive definite and uniformly
elliptic on compact subsets, and let $b,c\in L^\infty_{\rm loc}(G)$,
with $c\le0$ almost everywhere. Suppose that $w\ge0$ is locally
semiconcave with linear modulus and
\[
 A:D_A^2w+b\cdot Dw+cw\le0
       \quad\text{almost everywhere in }G.
\]
Then either $w>0$ everywhere in $G$ or $w\equiv0$ in $G$.
\end{lemma}
\begin{proof}
Work first on a relatively compact ball. Semiconcavity with constant
$C$ means that $w-C|x|^2/2$ is concave. Hence $Dw\in BV_{\rm loc}$,
and its distributional derivative is a symmetric matrix-valued measure
\cite[Proposition 1.1.3 and Theorem 2.3.1]{Cannarsa}. Write
\[
 D^2w=D_A^2w\,\mathcal L^N+D_s^2w.
\]
Here the absolutely continuous density is the Alexandrov Hessian.
For every $\xi\in\R^N$, the measure inequality
$\xi^TD^2w\xi\le C|\xi|^2\mathcal L^N$ implies
$\xi^TD_s^2w\xi\le0$. A common dominating positive measure and a
countable dense set of $\xi$ show that the density matrix of $D_s^2w$
is negative semidefinite almost everywhere with respect to that
measure. Since $A$ is continuous and positive definite,
$A:D_s^2w\le0$. The assumed almost-everywhere inequality therefore
implies the distributional inequality
\[
 \operatorname{div}(ADw)
 +(b-\operatorname{div}A)\cdot Dw+cw\le0.
\]
In the notation of \cite[Chapter 8, equations (8.1)--(8.8)]{GT},
the coefficient identification is
\[
 a^{ij}_{\rm GT}=A_{ij},\qquad b^i_{\rm GT}=0,\qquad
 c^i_{\rm GT}=b_i-\sum_j\partial_jA_{ji},\qquad d_{\rm GT}=c.
\]
Thus the coefficient inside the divergence multiplying $w$ is zero;
the drift is $c_{\rm GT}$. Since $b_{\rm GT}=0$, the weak
zero-order sign condition reduces to
\[
 \int d_{\rm GT}\varphi=\int c\varphi\le0
 \qquad(0\le\varphi\in C_c^1(G)).
\]
The function $w$ is locally Lipschitz, so it belongs to
$W^{1,2}$ on every smaller ball. All coefficients in this
divergence-form inequality are bounded there; the coefficient of $w$
inside the divergence is zero, and the zero-order sign is $c\le0$.

If $N=1$ or $N=2$, let $y\in\R^{3-N}$ range over a small ball,
and set $\widetilde w(x,y)=w(x)$, $\widetilde c(x,y)=c(x)$,
$\widetilde A(x,y)=\operatorname{diag}(A(x),I_{3-N})$ and
$\widetilde\beta(x,y)=(b(x)-\operatorname{div}A(x),0)$.
For any nonnegative smooth test function $\varphi$ compactly supported
in this product neighborhood, Fubini's theorem gives
\begin{align*}
 &\int\bigl(\widetilde A D\widetilde w\cdot D\varphi
       -(\widetilde\beta\cdot D\widetilde w+
                         \widetilde c\widetilde w)\varphi\bigr)\dd x\dd y\\
 &=\int\left[\int\bigl(A Dw\cdot D_x\varphi(x,y)
       -((b-\operatorname{div}A)\cdot Dw+cw)\varphi(x,y)\bigr)\dd x\right]\dd y
 \ge0.
\end{align*}
The last inequality is the original weak inequality tested against
$x\mapsto\varphi(x,y)$ for each fixed $y$.
The extended function remains nonnegative and locally semiconcave.
Thus it suffices to apply the weak Harnack inequality in dimension at
least three. If $x_0$ is a zero of $w$, choose $R>0$ so that the
corresponding $4R$ ball is compactly contained in the domain. Proposition~\ref{ext:weak-harnack}, applied in the possibly enlarged dimension, gives
\[
 \frac{1}{|B_{2R}|}\int_{B_{2R}(x_0)}w
       \le C_0\operatorname*{ess\,inf}_{B_R(x_0)}w=0.
\]
Continuity and $w(x_0)=0$ justify the last equality. Consequently $w$
vanishes on a neighborhood of $x_0$; for a dummy-variable extension,
restrict this conclusion to the original slice. The zero set is
therefore open as well as closed in $G$, and connectedness concludes
the proof.
\end{proof}

\begin{lemma}[Strong local profile comparison]\label{lem:strong-sphere}
Let $n\ge2$, $m=n-1$, $2\le q\le n+1$, and let
$U\subset\Delta_q^\circ$ be connected and relatively open.
Suppose $J$ is locally semiconcave on $U$ and, almost everywhere on $U$,
\[
 B_J=-D_A^2J\succ0,\qquad
 \tr\!\left[B_J^{-1}
  \left(\Id_E+\frac2{m^2}DJ\otimes DJ\right)\right]
 \le\frac2mJ.
\]
If $J\ge\Istd$ on $U$ and $J(v_*)=\Istd(v_*)$ for some $v_*\in U$,
then $J=\Istd$ throughout $U$.
\end{lemma}
\begin{proof}
For $B\succ0$ and $p\in E$, define
\[
 \Psi(B,p)=\tr B^{-1}+\frac2{m^2}p^TB^{-1}p.
\]
It is jointly convex. Indeed each term is the supremum of affine functions,
using $\tr B^{-1}=\sup_X(2\tr X-\tr X^TBX)$ and
$p^TB^{-1}p=\sup_y(2p\cdot y-y^TBy)$.
Linearize the supporting plane at
$B_*=-D^2\Istd$, $p_* =\nabla\Istd$, and set $w=J-\Istd$. The inequality
for $J$ and equality for the model give
\begin{equation}\label{eq:sphere-strong-linear}
 a:D^2w+b\cdot Dw-\frac2m w\le0\quad\text{almost everywhere},
\end{equation}
where
\[
 a=B_*^{-1}\left(\Id_E+\frac2{m^2}p_*p_*^T\right)B_*^{-1}\succ0,
 \qquad b=\frac4{m^2}B_*^{-1}p_*.
\]
On compact subsets these coefficients are smooth and bounded, and $a$
is uniformly elliptic. They depend only on the model. Subtracting the
smooth model preserves local semiconcavity with linear modulus, so
Lemma~\ref{lem:weak-minimum} applies with $A=a$ and $c=-2/m<0$.
It gives $w\equiv0$ because $w$ has an interior zero. Its proof includes
the case $\dim E=q-1=1$.
\end{proof}

\begin{proof}[Proof of Theorem~\ref{thm:sphere-main}]
For $n=2$ the assertion follows from Proposition~\ref{ext:planar-bubbles}, since
$2\le q\le3$. Henceforth assume $n\ge3$.
Optimality was proved in Proposition~\ref{prop:profile-equality}. If a
nonstandard minimizer existed at $v_0$, the preceding tube construction
would give $J\ge I=\Istd$ with equality at $v_0$. Proposition
\ref{prop:local-pdi} and Lemma~\ref{lem:strong-sphere} force $J=\Istd$
near $v_0$. Choose a nearby volume $v\in U_{n,q}$ from Corollary
\ref{cor:generic}, and an attained restricted minimizer $\Lambda$ there.
Since $P(\Lambda)=J(v)=\Istd(v)=I(v)$, this cluster is a global
minimizer and must be standard. However, Lemma~\ref{lem:sphere-tube} gives
\[
 \dist_{\dcl}(\Lambda,\cS(v))
 \ge \dist_{\dcl}(\cK,\cS(v))-\dist_{\dcl}(\Lambda,\cK)
 >3r-r=2r,
\]
contradicting $\Lambda\in\cS(v)$. Thus no nonstandard minimizer exists.
\end{proof}

\section{Euclidean perimeter optimality}
\label{sec:flat}
We now use ordinary, unnormalized volume and area. Set
\begin{equation}\label{eq:space-form-metric}
 g_s=\rho_s^2g_E,\qquad
 \rho_s(x)=\left(1+\frac{s|x|^2}{4}\right)^{-1},\qquad s\ge0.
\end{equation}
For $s>0$ this is the round sphere of curvature $s$ with one point removed.
Let $h(s,v)$ denote the perimeter of the standard cluster whose first $k$
cells have raw volumes $v=(v_1,\ldots,v_k)$; the omitted point is in the
remaining cell. For $v$ in a fixed compact subset of $(0,\infty)^k$, these
volumes are admissible when $s>0$ is sufficiently small. The local smooth
extension of $h$ through $s=0$, proved in Lemma~\ref{lem:model}, is a fact
about the standard family and does not presuppose its optimality.

\begin{proposition}[Flat consequence]\label{prop:flat}
For $n\ge2$ and $1\le k\le n$, a standard Euclidean $k$-bubble minimizes
perimeter at every positive finite-volume vector.
\end{proposition}
\begin{proof}
First suppose that all finite cells of a competitor $\Omega$ lie in a
fixed ball $B_R$, up to null sets. Its entire interface measure is then
supported in $\overline B_R$. The conformal area and volume formulas give
\[
 V_{g_s,i}(\Omega)=\int_{\Omega_i}\rho_s^n\dd x,\qquad
 P_{g_s}(\Omega)=\sum_{i<j}\int_{\Sigma_{ij}(\Omega)}
                    \rho_s^m\dd\HH^m,\qquad m=n-1.
\]
On $\overline B_R$, the mean-value theorem applied to
$(1+s|x|^2/4)^{-a}$, for $a=n,m$, gives
\[
 \|\rho_s^a-1\|_{L^\infty(B_R)}\le C_{n,R}s.
\]
Consequently, for $1\le i\le k$,
\[
 \begin{aligned}
 |V_{g_s,i}(\Omega)-V_{E,i}(\Omega)|&\le C_{n,R}s\,|\Omega_i|,\\
 |P_{g_s}(\Omega)-P_E(\Omega)|&\le C_{n,R}s\,P_E(\Omega).
 \end{aligned}
\]
The volume integrals have finite Lebesgue measure and the interface
integrals have finite $\HH^m$ measure, so dominated convergence yields
\[
 V_{g_s}(\Omega)\longrightarrow V_E(\Omega)=v,
 \qquad P_{g_s}(\Omega)\longrightarrow P_E(\Omega).
\]
Complete the chart to the round sphere, assigning the omitted point to
the exterior cell. This adds neither volume nor perimeter. The spherical
profile equality, with raw measures rescaled to curvature $s$, implies
\[
 P_{g_s}(\Omega)\ge h(s,V_{g_s}(\Omega)).
\]
The joint continuity supplied by Lemma~\ref{lem:model} therefore gives
$P_E(\Omega)\ge h(0,v)$.

We next justify the reduction to bounded finite cells explicitly. For an
arbitrary competitor let $U=\bigcup_{i=1}^k\Omega_i$. This set has finite
volume. By the coarea formula,
\[
 \int_0^\infty \HH^m(U\cap\partial B_t)\dd t=|U|<\infty.
\]
Choose regular slicing radii $R_j\to\infty$ such that
$a_j:=\HH^m(U\cap\partial B_{R_j})\to0$. For example, one may choose a
radius in each interval $(j,j+1)$ where the slice area is at most twice
its integral on that interval, while avoiding the null set of exceptional
BV slicing radii. Define
\[
 \Omega_i^{(j)}=\Omega_i\cap B_{R_j}\quad(1\le i\le k),
 \qquad
 \Omega_q^{(j)}=\R^n\setminus\bigcup_{i=1}^k\Omega_i^{(j)}.
\]
The BV product and trace formulas show that
\[
 P_E(\Omega^{(j)})
 \le P_E(\Omega;B_{R_j})+a_j.
\]
Here each newly exposed finite-cell/exterior interface on the slicing
sphere is counted once; existing interfaces inside the ball are
unchanged. Moreover
\[
 \dcl(\Omega^{(j)},\Omega)=2|U\setminus B_{R_j}|\longrightarrow0,
 \qquad V_E(\Omega^{(j)})\longrightarrow v.
\]
The preceding upper bound and lower semicontinuity of cluster perimeter
give $P_E(\Omega^{(j)})\to P_E(\Omega)$. For large $j$ all finite volumes
are positive. Apply the bounded-competitor inequality at the actual
volume vector $V_E(\Omega^{(j)})$ and pass to the limit, using continuity
of $h(0,\cdot)$. Finally,
the standard competitor attains $h(0,v)$, proving equality of the true
and standard Euclidean profiles.
\end{proof}

For comparison with normalized measures, put $\omega_n=|\Sn|$. With raw
measures the resulting flat-limit identity is
\[
 I_{\R^n}(v)=\lim_{\eps\downarrow0}\eps^{-m/n}
 I^{\mathrm{raw}}_{\Sn}
 \left(\eps v_1,\ldots,\eps v_k,
       \omega_n-\eps\sum_{i=1}^k v_i\right).
\]
Indeed, set $s=\eps^{2/n}$, use spherical profile equality, and rescale
the sphere of curvature $s$ to the unit sphere. This is consistent with
\cite[Section 13.5]{MNP}. The argument transfers perimeter optimality;
an equality classification requires the separate argument below because
a strictly positive perimeter gap can disappear in a limit.

\section{The standard model through zero curvature}
\label{sec:model}
\begin{lemma}[Smoothness and strict concavity of the model]\label{lem:model}
Fix $n\ge2$, $1\le k\le n$, and $v_0\in(0,\infty)^k$. The standard
profile $h(s,v)$ in \eqref{eq:space-form-metric} extends smoothly to a
neighborhood of $(0,v_0)$ in $\R\times\R^k$, and
\[
 -D^2_{vv}h(s,v)\succ0
\]
there. Small negative $s$ is used only to extend the model smoothly on a
bounded chart; no assertion about negative-curvature minimization is
made here.
\end{lemma}
\begin{proof}
Put $q=k+1$. Choose vectors $u_1,\ldots,u_k\in\R^n$ satisfying
\[
 u_i\cdot u_j=\frac{1+\delta_{ij}}2,
 \qquad u_q=0.
\]
Indeed, for $a\in\R^k\setminus\{0\}$,
\[
 a^T\frac{I_k+\one\one^T}{2}a
 =\frac{|a|^2+(\one\cdot a)^2}{2}>0.
\]
A square root of this Gram matrix, followed by the inclusion
$\R^k\hookrightarrow\R^n$, gives the required vectors.
Let $U:\R^n\to\R^k$ have rows $u_i^T$; then
$UU^T=(I_k+\one\one^T)/2$ is invertible. For
$\kappa\in(0,\infty)^k$, set $\kappa_q=0$ and
\begin{equation}\label{eq:model-scores}
 Q_q=0,\qquad Q_i(x)=\frac{\kappa_i}{2}|x|^2-u_i\cdot x,
 \qquad \Omega_i(\kappa)=\{Q_i<Q_j\text{ for all }j\ne i\}.
\end{equation}
Since $|u_i|=1$, completing the square gives
\[
 Q_i(x)=\frac{\kappa_i}{2}
 \left(\left|x-\frac{u_i}{\kappa_i}\right|^2
                  -\frac1{\kappa_i^2}\right).
\]
Thus each finite cell lies in $B_{1/\kappa_i}(u_i/\kappa_i)$.
For distinct $i,j$, including $q$, put
$\alpha=\kappa_i-\kappa_j$ and $u=u_i-u_j$.
On $Q_i=Q_j$ we have $u\cdot x=\alpha|x|^2/2$ and $|u|=1$, so
\[
 |\nabla(Q_i-Q_j)|^2=|\alpha x-u|^2
 =\alpha^2|x|^2-2\alpha u\cdot x+|u|^2=1.
\]
This also covers the planar case $\alpha=0$.

To identify the family, set
\[
 b=\tfrac12U^T(UU^T)^{-1}\kappa,\qquad
 y=\frac{x}{|x|^2},\qquad y'=y-b.
\]
For $x\ne0$ and every $i$, including $q$, the score identity is
\[
 \frac{Q_i(x)}{|x|^2}=\frac{\kappa_i}{2}-u_i\cdot(y'+b)
 =-u_i\cdot y'.
\]
Write $\bar u=q^{-1}\sum_i u_i$ and $a_i=u_i-\bar u$.
Subtracting the common score $-\bar u\cdot y'$ leaves $-a_i\cdot y'$,
and the Gram identities give
\[
 a_i\cdot a_j=\frac12\left(\delta_{ij}-\frac1q\right),
 \qquad |a_i-a_j|=1\quad(i\ne j).
\]
These are the scores of a centered regular simplex. Each cell is nonempty:
at $y'=a_i$ its score is strictly smallest. Each pair face is nonempty as
well. For $q\ge3$, at $y'=a_i+a_j$ the scores of $i,j$ agree and are
strictly smaller than all the others; their difference has nonzero
gradient $a_i-a_j$. For $q=2$ the pair face is a hyperplane.
The inverse of $y'=x/|x|^2-b$ is a conformal transformation of the
one-point compactification. The point $x=\infty$ corresponds to $y'=-b$,
where the scores are $\kappa_i/2>0$ for $i\le k$ and $0$ for $i=q$.
Thus this pole lies strictly inside the exterior cell. It follows that
\eqref{eq:model-scores} defines a standard Euclidean cluster with positive
finite-cell volumes and positive area on every pair interface.

Conversely, write a centered regular-simplex spherical partition with
score vectors $a_i\in\R^{n+1}$. Its common meeting set is
\[
 \Sn\cap\operatorname{span}\{a_i-a_q:1\le i\le k\}^{\perp}
 \simeq\mathbb S^{n-k}\ne\varnothing,
\]
because those $k$ vectors are independent and $k\le n$.
Conformal transformations preserve this common meeting set. A
stereographic pole lying inside the exterior cell is disjoint from it;
hence the Euclidean cluster has a finite common meeting point.
Translate that point to $0$ and subtract the exterior score. The constant
terms vanish, and the normalized scores take the form
$Q_i(x)=\kappa_i|x|^2/2-u_i\cdot x$ with $u_q=0$ and
$|u_i-u_j|=1$. Consequently
\[
 |u_i|=1,\qquad
 u_i\cdot u_j=\tfrac12(|u_i|^2+|u_j|^2-|u_i-u_j|^2)
 =\frac{1+\delta_{ij}}2\quad(1\le i,j\le k).
\]
A Euclidean orthogonal transformation identifies these vectors with our
fixed choice. Choose $R$ so that every finite cell is contained in $B_R$. The exterior score
is minimal almost everywhere on $\{|x|>R\}$, so continuity gives
$Q_i(x)\ge0$ throughout that open set. In particular, for every $t>R$,
\[
 0\le Q_i(tu_i)=\frac{\kappa_i}{2}t^2-t,
\]
which forces $\kappa_i>0$. Existence and uniqueness of the standard family
at prescribed positive volumes \cite{Montesinos,MNS,MNP} therefore give
a parameter $\kappa^0$ representing $v_0$.

The same fixed subsets are standard clusters in the completed spherical
metric for $s>0$, because stereographic projection and the maps used
above are conformal. Their physical pressures can be checked directly.
For a normalized Euclidean interface
\[
  \begin{gathered}
  Q(x)=\frac{\kappa_E}{2}|x|^2+c\cdot x+d=0,
  \qquad |c|^2-2\kappa_Ed=1,\\
  \nu_E=\kappa_E x+c,
  \qquad |\nu_E|^2
   =\kappa_E^2|x|^2+2\kappa_Ec\cdot x+|c|^2=1.
  \end{gathered}
\]
the convention $A=D\nu$ and the conformal curvature rule give
\[
 \kappa_{g_s}
 =\rho_s^{-1}\bigl(\kappa_E+\partial_{\nu_E}\log\rho_s\bigr)
 =\kappa_E+\frac s2d.
\]
Indeed $x\cdot\nu_E=\kappa_E|x|^2/2-d$ and
$\partial_{\nu_E}\log\rho_s=-(s\rho_s/2)x\cdot\nu_E$.
In \eqref{eq:model-scores} every constant term vanishes, so the pressure
of the $i$-th finite cell relative to the exterior is $p_i=m\kappa_i$.
Triple junctions retain their $120$-degree angles under the conformal
change. The first-variation identity therefore holds for these model
clusters for every small $s$, including the locally defined negative
values used for extension.

Let $V_s(\kappa)$ be the vector of finite volumes. If
$\dot\kappa=a$ and $a_q=0$, differentiating the interface equation
shows that the outward Euclidean normal speed of cell $i$ at its
interface with $j$ is $-(a_i-a_j)|x|^2/2$. Multiplication of normal
speed by $\rho_s$ and of interface area by $\rho_s^m$ yields
\begin{equation}\label{eq:model-laplacian}
 D_\kappa V_s=-W_s,\qquad
 a^TW_sa=\sum_{i<j}(a_i-a_j)^2
 \int_{\Sigma_{ij}(\kappa)}\rho_s^n\frac{|x|^2}{2}\dd\HH^m.
\end{equation}
The interfaces have positive area and $|x|^2$ vanishes only at one
point, a set of zero $\HH^m$ measure since $m\ge1$. All weights are
strictly positive. If the quadratic form vanishes, all differences
$a_i-a_j$ vanish, and $a_q=0$ forces $a=0$. Hence $W_s\succ0$.

We justify all parameter differentiations on fixed integration domains.
For $b$ small, the special conformal map
\[
 T_b(x)=\frac{x-|x|^2b}{1-2b\cdot x+|b|^2|x|^2}
\]
maps the cluster with parameter $\kappa^0$ to that with parameter
$\kappa^0-2Ub$. In fact, substituting $T_b^{-1}=T_{-b}$ into each score
gives the latter score divided by the common positive denominator
$1+2b\cdot x+|b|^2|x|^2$. Set
\[
 b(\kappa)=-\tfrac12 U^T(UU^T)^{-1}(\kappa-\kappa^0).
\]
Choose a compact neighborhood of the closures of the finite model cells.
After reducing the parameter neighborhood, both $T_{b(\kappa)}$ and its
inverse are smooth there, their denominators are bounded below by a
positive constant, and all transformed finite cells stay in one compact
set. On this set $1+s|x|^2/4$ is also bounded below for $|s|$ small.
Put $T_\kappa=T_{b(\kappa)}$, $\Omega_i^0=\Omega_i(\kappa^0)$ and
$\Sigma_{ij}^0=\Sigma_{ij}(\kappa^0)$. The fixed-domain formulas are
\[
 \begin{aligned}
 V_{s,i}(\kappa)&=\int_{\Omega_i^0}
  \rho_s(T_\kappa x)^n|\det DT_\kappa(x)|\dd x,\\
 \widetilde h(s,\kappa)&=\sum_{i<j}\int_{\Sigma_{ij}^0}
 \rho_s(T_\kappa x)^m
 \operatorname{Jac}_m(DT_\kappa|T_x\Sigma_{ij}^0)\dd\HH^m(x).
 \end{aligned}
\]
The volume formula is used for $1\le i\le k$.
For each multi-index $\alpha$ in $(s,\kappa)$, the absolute value of
the corresponding derivative of either integrand is bounded by a
constant $C_\alpha$ on a smaller parameter neighborhood and the common
compact spatial set. Each $\Omega_i^0$ has finite volume and each
$\Sigma_{ij}^0$ has finite area. Dominated differentiation therefore
allows every $\partial_{(s,\kappa)}^\alpha$ to pass under these integrals.
This proves joint $C^\infty$ dependence, including at $s=0$.

Apply the inverse function theorem to
$(s,\kappa)\mapsto(s,V_s(\kappa))$ at $(0,\kappa^0)$. Its derivative
is block triangular with invertible lower-right block $-W_0$.
It supplies a smooth inverse $\kappa=\kappa(s,v)$. If
$\widetilde h(s,\kappa)=P_{g_s}(\Omega(\kappa))$, define
$h(s,v)=\widetilde h(s,\kappa(s,v))$. This agrees with the standard
profile for $s\ge0$. At fixed $s$, stationarity gives
$\dd_\kappa\widetilde h=m\kappa\cdot\dd_\kappa V_s$, whence
\[
 D_vh=m\kappa(s,v),\qquad
 D^2_{vv}h=-mW_s^{-1}\prec0.
\]
All inverses here are finite-coordinate $k\times k$ inverses. The strict
inequality persists throughout a smaller full neighborhood, completing
the proof.
\end{proof}

\section{A centered local problem near Euclidean equality cases}
\label{sec:euclidean-local}
Throughout the geometric construction in this section assume $n\ge3$;
the planar main theorem is handled by Proposition~\ref{ext:planar-bubbles}.
The rest of the proof uses only spherical profile equality and generic
spherical uniqueness, not the all-volume spherical uniqueness theorem.

\subsection{Uniform compactness of the initial family}
For a bounded-film cluster define
\[
 S_s(\Omega)=\sum_{i=1}^k\int_{\Omega_i}x\,\rho_s^n\dd x,
 \qquad
 \dcl(\Omega,\Lambda)=\sum_{i=1}^{k+1}|\Omega_i\symdiff\Lambda_i|.
\]
The distance uses Euclidean measure in the fixed chart. Differences of
exterior cells have finite measure for the clusters considered here.

\begin{lemma}[Centered equality cases are compact]\label{lem:euclidean-compact}
At a fixed $v_0\in(0,\infty)^k$, all Euclidean global minimizers with
$S_0=0$ form a compact $L^1$ family $\cM_0$, and their finite cells lie in
one common ball. Their standard subset $\cS_0$ is compact and relatively
open in $\cM_0$.
\end{lemma}
\begin{proof}
By Proposition~\ref{prop:flat}, the true profile equals the smooth model
$h(0,\cdot)$ near $v_0$. If $L$ is a local Lipschitz constant and a local
competitor $\Lambda$ has volumes $v'$ close to $v_0$, then
\[
 P_E(\Omega)=h(0,v_0)\le h(0,v')+L|v'-v_0|
 \le P_E(\Lambda)+L\dcl(\Lambda,\Omega).
\]
Thus global minimizers at $v_0$ have uniform quasiminimality and density
constants. Their total boundaries are connected and bounded by
Proposition~\ref{ext:connected-global}; boundedness is also given explicitly by
\cite[Theorem 2.3]{MNS}. These statements are used only for global
minimizers at this point.

Fix a radius smaller than the uniform density radius. Every ball of that
radius centered on the total boundary contains a uniform positive amount
of at least one finite phase: at a boundary point some finite label is
incident, and infiltration gives the lower bound. The total finite volume
therefore bounds the number of disjoint such balls. Choose a maximal family of pairwise disjoint balls $B_a(x_j)$ with
centers on the boundary, for a fixed sufficiently small $a$. Each contains
at least $c a^n$ of the union of finite cells, so their number is at most
$V_{\mathrm{tot}}/(c a^n)$. Maximality implies that the balls $B_{3a}(x_j)$
cover the boundary. The intersection graph of their nonempty intersections
with the boundary is connected: otherwise their unions would give a
separation of the connected boundary into disjoint relatively open sets.
An intersecting pair has centers at distance at most $6a$. Chaining at
most the uniformly bounded number of balls bounds the boundary diameter. The exterior of a sufficiently large ball containing the film is
connected and has no boundary of the cluster; it belongs to the infinite
cell. Every finite cell is thus inside a ball of uniformly bounded radius.
Its center is also bounded after imposing $S_0=0$, because the barycenter of
the union of the finite cells is $0$ and lies in the convex hull of that
ball.

BV compactness now gives $L^1$ compactness with no loss of volume. The common
support bound preserves the centroid constraint in the limit. Lower
semicontinuity and optimality show that every limit is again a minimizer,
and all convergent sequences have strict perimeter convergence.
Lemma~\ref{lem:pair-areas} gives convergence of every pair-interface area.
If nonstandard minimizers approached a standard one, all their pair
interfaces would eventually be present. The Euclidean structural theorem
and the full-interface characterization (Proposition~\ref{ext:full-interfaces})
would make them standard. This proves relative openness. The standard
family is unique at its positive volumes up to rigid motion
\cite{Montesinos}; after centering only a compact orthogonal orbit remains.
\end{proof}

Assume for contradiction that
\[
 \cK=\cM_0\setminus\cS_0\ne\varnothing.
\]
It is compact and has positive distance from $\cS_0$. Choose $R$ so large
that all finite cells of every member of $\cM_0$ lie in $B_{R/10}$, and
choose a small tube radius $r>0$ separated from $\cS_0$. For $s\ge0$ small
and $v$ near $v_0$ set
\begin{equation}\label{eq:curved-local-profile}
 J(s,v)=\inf\left\{
 P_{g_s}(\Lambda):
 \begin{array}{l}
 V_{g_s}(\Lambda)=v,\quad S_s(\Lambda)=0,\\
 |\Lambda_i\setminus B_R|=0\ (1\le i\le k),\\
 \dist_{\dcl}(\Lambda,\cK)\le r
 \end{array}\right\}.
\end{equation}
The admissible class has prescribed volumes, zero total first moment,
support in $B_R$, and an $L^1$ tube restriction. The next three lemmas
establish attainment, a positive distance from the support wall, and
isometric recentering. These properties allow ordinary small
volume-preserving variations of its minimizers.

\subsection{Uniform volume--moment correction}
\begin{lemma}[Attainment and interior margin]\label{lem:curved-attain}
On a sufficiently small parameter neighborhood, the infimum in
\eqref{eq:curved-local-profile} is attained. All restricted minimizers have
a uniform perimeter bound and approach $\cK$ uniformly in $L^1$ as
$(s,v)\to(0,v_0)$. In particular they can be required to have distance
less than $r/8$ from $\cK$. There are uniform smooth correction maps for
both volumes and moments, supported away from $\partial B_R$.
\end{lemma}
\begin{proof}
At each $\Omega\in\cK$ choose $k$ smooth regular-patch fields whose finite
volume derivative matrix $A$ is invertible. Add $n$ cutoff translations
$Y_a=\chi a$, where $\chi=1$ on $B_{R/5}$ and $\chi=0$ outside $B_{R/3}$.
At $(s,\Omega)=(0,\Omega)$ their volume derivatives vanish and
\[
 \delta_{Y_a}S_0=V_{\mathrm{tot}}a,
 \qquad V_{\mathrm{tot}}=\sum_i(v_0)_i>0.
\]
The combined volume--moment derivative is
\[
 \begin{pmatrix}A&0\\D&V_{\mathrm{tot}}I_n\end{pmatrix},
\]
which is invertible. For any finite-perimeter partition the derivatives
of these maps are integrals of smooth bounded functions:
\[
 (\delta_XV_{g_s})_i=\int_{\Lambda_i}\operatorname{div}_{g_s}X\dd\mu_s,
 \qquad
 \delta_XS_s=\sum_i\int_{\Lambda_i}
 (X+x\operatorname{div}_{g_s}X)\dd\mu_s.
\]
Thus they are continuous under $L^1$ convergence on $B_R$. Finite covering
of $\cK$ gives uniform inverse and $C^2$ bounds, exactly as in Lemma
\ref{lem:repair}. These maps are the identity near the support wall.
More precisely, there are an $L^1$ neighborhood $\mathcal N$ of
$\cK$, numbers $\eta,\delta,C>0$, and finitely many flow families.
Every $\Lambda\in\mathcal N$ whose finite cells satisfy the support
condition in \eqref{eq:curved-local-profile} admits one of these families for which
the volume--moment map has first and second derivatives bounded by $C$
and an $a$-derivative with inverse norm at most $C$ whenever
$|s|+|a|<\eta$. Its parameter inverse is defined for every target error
of norm less than $\delta$, with the same constants for all such
$\Lambda$. Here $\eta$ is also chosen so that
$1+s|x|^2/4>0$ throughout the common chart, including the small
negative $s$ used only to extend individual smooth support functions.

Correcting a member of $\cK$ for the small metric and volume changes gives
an admissible competitor with
\[
 J(s,v)\le h(0,v_0)+C(s+|v-v_0|).
\]
At fixed $(s,v)$, uniform metric equivalence and the preceding upper
bound give a Euclidean BV bound for every finite cell of a minimizing
sequence. These cells are supported in $\overline B_R$, so a subsequence
converges globally in $L^1$; the exterior is determined by their
complement. The bounded continuous weights $\rho_s^n$ and $x\rho_s^n$
preserve the volume and moment constraints. The support constraint and
the closed tube are also $L^1$-closed. Lower semicontinuity of weighted
perimeter therefore proves attainment.

For $(s_j,v_j)\to(0,v_0)$, write $\Lambda_j$ for any corresponding
restricted minimizers. On $B_R$,
\[
 \|\rho_{s_j}^n-1\|_\infty+
 \|\rho_{s_j}^{m}-1\|_\infty\le C_Rs_j.
\]
The support and perimeter bounds therefore imply
\[
 V_E(\Lambda_j)\to v_0,\qquad S_0(\Lambda_j)\to0,
 \qquad |P_E(\Lambda_j)-P_{g_{s_j}}(\Lambda_j)|\to0.
\]
Every subsequential $L^1$ limit $\Omega$ has those volumes and moment,
and the competitor bound gives
\[
 h(0,v_0)\le P_E(\Omega)\le\liminf_jP_E(\Lambda_j)
 \le\limsup_jP_E(\Lambda_j)\le h(0,v_0).
\]
Thus $\Omega\in\cM_0$ and perimeter convergence is strict. The limit
belongs to the closed tube; separation excludes $\cS_0$, so
$\Omega\in\cK$. If the approach to $\cK$ were not uniform, there
would be $\varepsilon>0$ and such minimizers with distance at least
$\varepsilon$ from $\cK$; a convergent subsequence would contradict
that bound. This proves the asserted approach and the $r/8$ tube margin.
\end{proof}

\subsection{Removing the artificial support wall}
\begin{lemma}[Support clearance]\label{lem:wall}
After shrinking the parameter neighborhood, every restricted minimizer has
all finite cells in $B_{R/5}$, up to null sets.
\end{lemma}
\begin{proof}
Let $U=\bigcup_{i=1}^k\Lambda_i$. Erase all finite labels in a small
Euclidean ball $B_t(x)$, assigning the erased region to the exterior.
This competitor respects the support constraint even if the ball crosses
$\partial B_R$. Its volume and moment errors are bounded by
$C_R\mu_s(U\cap B_t(x))$. Repair them using Lemma~\ref{lem:curved-attain}.
For $t\le t_0$ uniformly small, both the deletion and its correction stay
inside the tube. If the uncorrected perimeter is uniformly bounded, the
repair costs at most $C\mu_s(U\cap B_t)$. If its new slice area already
exceeds the original uniform perimeter bound, the following inequality is
immediate instead. Thus for almost every $t\le t_0$,
\begin{equation}\label{eq:wall-deletion}
 P_{g_s}(\Lambda;B_t(x))
 \le\HH_{g_s}^{n-1}(U\cap\partial B_t(x))
 +C\mu_s(U\cap B_t(x)).
\end{equation}
The repairs do not have to avoid this ball; only their global cost enters
this comparison. They do vanish near the artificial wall.

Put $f(t)=\mu_s(U\cap B_t(x))$. The local isoperimetric inequality for
$U\cap B_t$, the bound $\Per_{g_s}(U;B_t)\le P_{g_s}(\Lambda;B_t)$,
and coarea imply
\[
 c f(t)^{(n-1)/n}\le C_1f'(t)+C_2f(t).
\]
Uniform metric equivalence accounts for the Euclidean radii. Decrease
$t_0$ to absorb the last term, using $f(t)^{1/n}\le Ct$. Integrating as in
Lemma~\ref{lem:infiltration} gives a uniform $\eta>0$ with
\[
 f(t_0)<\eta t_0^n\quad\Longrightarrow\quad f(t_0/2)=0.
\]
Choose $t_0<R/20$. Every such ball centered outside $B_{R/5}$ is disjoint
from $B_{R/10}$. Uniform $L^1$ closeness to $\cK$ makes its finite mass
smaller than the threshold. A countable cover of the complement of
$B_{R/5}$ proves the claim. This proof uses neither boundary regularity at
the wall nor an infiltration estimate for the exterior at the wall.
\end{proof}

\subsection{Isometric recentering}
\begin{lemma}[Exact isometric recentering]\label{lem:gauge}
Every restricted minimizer is stationary and stable for the ordinary
finite-volume constraints in $g_s$, and has ordinary uniform local
quasiminimality, density and regularity. The moment constraint contributes
no multiplier.
\end{lemma}
\begin{proof}
For $a\in\R^n$ the fields
\begin{equation}\label{eq:killing-gauge}
 X_{s,a}(x)=\left(1-\frac{s|x|^2}{4}\right)a
            +\frac s2(a\cdot x)x
\end{equation}
are Killing fields of $g_s$. Their Euclidean conformal factor is
$s(a\cdot x)/2$, canceled by $X_{s,a}(\log\rho_s)$. In particular
$\operatorname{div}_{g_s}X_{s,a}=0$. Write $U=\bigcup_{i=1}^k\Lambda_i$ and
$V_{\mathrm{tot}}(v)=\sum_{i=1}^k v_i$. At a centered cluster,
\[
 \delta_{X_{s,a}}S_s=\mathsf M_s(\Lambda)a,\qquad
 \mathsf M_s(\Lambda)=V_{\mathrm{tot}}(v)I_n+
 \frac{s}{4}\int_U(2x\otimes x-|x|^2I_n)\dd\mu_s.
\]
The matrix $2x\otimes x-|x|^2I_n$ has operator norm $|x|^2$.
Since $U\subset B_R$ up to null sets,
\[
 \|\mathsf M_s(\Lambda)-V_{\mathrm{tot}}(v)I_n\|
 \le\frac{sR^2}{4}V_{\mathrm{tot}}(v).
\]
For $sR^2\le2$ and $v$ sufficiently near $v_0$,
\[
 \mathsf M_s(\Lambda)\succeq\tfrac12V_{\mathrm{tot}}(v)I_n,
 \qquad
 \|\mathsf M_s(\Lambda)^{-1}\|\le\frac{2}{V_{\mathrm{tot}}(v)}
 \le\frac{4}{V_{\mathrm{tot}}(v_0)}.
\]
Compose the short flows in the $n$ coordinate directions to obtain an
isometry $T_{s,b}$. At each centered restricted minimizer the map
$b\mapsto S_s(T_{s,b}\Lambda)$ has the displayed invertible derivative.
Its derivatives depend continuously on the finite-cell indicators in
$L^1(B_R)$, so the inverse function theorem gives a uniform recentering
parameter for every sufficiently small physical perturbation. For a
smooth perturbation this parameter is smooth in its time variable.
The isometry is considered on the completed round sphere; its restriction
to the support stays in the chart by the support clearance. Thus it
preserves both perimeter and all volumes exactly and returns the
competitor to the moment gauge, while the small displacement preserves
the tube margin.

In particular every sufficiently small ordinary volume-preserving
variation is admissible after an area-preserving recentering. Its
perimeter as a function of time is identical before and after
recentering, not merely to first order. Constrained stationarity and
stability therefore hold for the ordinary volume constraints by the
argument of Lemma~\ref{lem:quasimin}; no extra moment multiplier enters
the index form.

For local perimeter comparisons meeting the actual film, Lemma
\ref{lem:wall} leaves a uniform distance from the artificial wall.
Volume--moment corrections then prove ordinary quasiminimality by the
argument of Lemma~\ref{lem:quasimin}. A small comparison ball disjoint
from the film has zero original local perimeter and needs no argument.
The regularity and density conclusions follow from Lemmas
\ref{lem:infiltration} and \ref{lem:regularity}.
\end{proof}

\subsection{Uniform reflection in the bounded chart}
\begin{lemma}[Great-sphere reflections through the flat limit]
\label{lem:chart-reflection}
Put $D_s(x)=1+s|x|^2/4$. A great sphere in the $g_s$ chart has equation
\[
 f_{s,\beta,a}(x)=\beta\cdot x+a(1-s|x|^2/4)=0,
 \qquad |\beta|^2+sa^2=1.
\]
Its reflection is
\begin{equation}\label{eq:chart-reflection}
 R_{s,\beta,a}(x)=
 \frac{x-2f_{s,\beta,a}(x)\beta}{1-sa f_{s,\beta,a}(x)}.
\end{equation}
If the great sphere meets $B_A$ and $sA^2<1$, then
$|a|\le A/(1-sA^2/4)$. Thus, along $s\downarrow0$, these reflections have
subsequences converging smoothly on every fixed bounded chart to affine
Euclidean reflections.
\end{lemma}
\begin{proof}
The isometric embedding into the radius-$s^{-1/2}$ sphere is
\[
 x\longmapsto\left(\frac{x}{D_s(x)},
 \frac{1-s|x|^2/4}{\sqrt s\,D_s(x)}\right).
\]
Reflect in the unit normal $(\beta,\sqrt s\,a)$ and apply inverse
stereographic projection to obtain \eqref{eq:chart-reflection}. At an
intersection point with $B_A$, the equation $f=0$ gives the bound on $a$.
It bounds $a$, forces $|\beta|\to1$, and keeps the denominator uniformly
away from zero on fixed bounded sets for small $s$. Differentiation of
the explicit formula proves the asserted smooth convergence. At $s=0$
the formula is $x\mapsto x-2(\beta\cdot x+a)\beta$.
\end{proof}

\subsection{Structure at positive curvature}
\begin{lemma}[Structured local minimizers and their PDI]\label{lem:curved-structure}
For every sufficiently small $s>0$ and $v$ near $v_0$, there is a restricted
minimizer which, on the completed round sphere, is a regular stable
perpendicular $3$-Plateau spherical Voronoi cluster with equatorial cells.
For almost every $v$ at each such $s$,
\begin{equation}\label{eq:curved-pdi}
 B=-J_{vv}\succ0,\qquad
 \frac{ms}{2}\tr(MB^{-1})+\frac1m p^TB^{-1}p\le J,
 \qquad p=J_v,\quad M=I_k+\one\one^T.
\end{equation}
The same inequality holds at almost every joint Alexandrov point $(s,v)$
with $s>0$.
\end{lemma}
\begin{proof}
We verify the two geometric issues caused by the bounded chart.

\emph{Reflection.} A common great-sphere bisector of the $k$ finite cells
exists since $k\le n$; it also bisects the exterior. Such a bisector meets
the common support ball. Lemma~\ref{lem:chart-reflection} makes the
following convergence uniform. If $s_j\downarrow0$ and restricted minimizers
approach $\Omega\in\cK$, the bisectors subconverge in the bounded chart
to an affine Euclidean hyperplane $x\cdot N=d$ bisecting every finite
cell of $\Omega$.
Write its normalized global quadratic scores as
\[
 Q_i(x)=\kappa_i|x|^2/2+c_i\cdot x-d_i,\qquad Q_q=0.
\]
Set $a_i=c_i\cdot N+\kappa_i d$, including $a_q=0$. Reflection $R$ in
that hyperplane gives
\[
 Q_i(Rx)-Q_i(x)=-2(x\cdot N-d)a_i.
\]
At almost every point of the positive half-space, let $i,j$ be the unique
winners at $x,Rx$. The winner inequalities imply $a_j\ge a_i$; if $i\ne j$,
their strict form gives
\[
 0<Q_j(x)-Q_i(x)<2(x\cdot N-d)(a_j-a_i),
\]
so $a_j>a_i$. Write $H^\pm=\{\pm(x\cdot N-d)>0\}$ and let
$\iota(x)$ denote the winning label. Since $a_q=0$ and the finite cells
are bounded, the following integrand is integrable, and bisection gives
\[
 \int_{H^+}\bigl(a_{\iota(Rx)}-a_{\iota(x)}\bigr)\dd x
 =\sum_{i=1}^k a_i\bigl(|\Omega_i\cap H^-|-|\Omega_i\cap H^+|\bigr)=0.
\]
The integrand is nonnegative and strictly positive wherever the two
winning labels differ. Thus $R\Omega_i=\Omega_i$ up to null sets for
every $i$. Set $U=\bigcup_{i=1}^k\Omega_i$ and
$\mathsf R_N=I_n-2N\otimes N$. The moment of the reflected union is
\[
 0=S_0(\Omega)=\int_{RU}x\dd x
 =\mathsf R_N\int_Ux\dd x+2dN|U|=2dN|U|.
\]
Since $|U|=\sum_i(v_0)_i>0$, this implies $d=0$.

We also detail the moving-reflection limit in the fixed Euclidean
measure used by the tube. Put $\chi_{j,i}=\mathbf1_{\Lambda_{j,i}}$
and $\chi_i=\mathbf1_{\Omega_i}$ for $i\le k$.
Lemma~\ref{lem:chart-reflection} gives $R_j\to R$ smoothly on bounded
sets. Their images of a fixed ball $B_L$ lie in another fixed ball
$B_{L'}$, and their inverse Jacobians are uniformly bounded there.
Changing variables yields
\[
 \|\chi_{j,i}\circ R_j-\chi_i\circ R\|_{L^1(B_L)}
 \le C\|\chi_{j,i}-\chi_i\|_{L^1(B_{L'})}
       +\|\chi_i\circ R_j-\chi_i\circ R\|_{L^1(B_L)}\longrightarrow0.
\]
For the second term, approximate $\chi_i$ in $L^1(B_{L'})$ by
continuous functions: the uniform Jacobian bound controls both
approximation errors, and uniform convergence of $R_j$ controls the
continuous-function term. The two half-space indicators converge
almost everywhere, since their limiting hyperplane has zero Lebesgue
measure. The reflected-cell indicator formula from Lemma~\ref{lem:structured}
therefore proves local $L^1$ convergence of both hemisphere copies to
$\Omega$. All their finite cells are supported in one fixed ball,
by the reflection bound below, so this convergence is global in the
cluster distance, including the exterior label. On the same ball
$\rho_{s_j}^n\to1$ uniformly; consequently both reflected moments
tend to $S_0(\Omega)=0$. The uniform inverse bound in Lemma~\ref{lem:gauge}
now provides recentering isometries $T$ tending to the identity.
If the reflected cluster is invariant under $R$, its recentered image is
invariant under the great-sphere reflection $TRT^{-1}$.
The reflected clusters remain in the tube. They also respect support:
a Euclidean reflection in a plane meeting $B_a$ maps $B_a$ into $B_{3a}$;
the nearby spherical reflections obey the same bounded-chart estimate
with an $o(1)$ error by Lemma~\ref{lem:chart-reflection}. Taking $a=R/5$ leaves a margin inside $B_R$ before and
after the small recentering. Their two perimeters have sum at most twice
the original perimeter. Restricted minimality therefore allows a symmetric
restricted minimizer to be chosen. By contradiction and compactness these
statements hold uniformly for all parameters in a small neighborhood of
$(0,v_0)$.

\emph{Connectedness.} By Lemma~\ref{lem:finite-components}, each restricted
minimizer has finitely many open cell components. Only the exterior
component containing the complement of the support ball can leave that
ball; the uniform component-infiltration argument applies to all others. Uniform density
and the support clearance give Hausdorff convergence of films to the
connected film of a member of $\cK$. If nearby films were disconnected,
apply Lemma~\ref{lem:collision} on the completed sphere and move the
closed group not containing the omitted point. It is contained in the
bounded chart, and its first-contact displacement tends to zero. Recenter
by the isometries just constructed. All these paths retain their support
and tube margins. The collision contradiction applies unchanged. Thus the
film is connected for every sufficiently nearby restricted minimizer.

The symmetric restricted minimizer is now bounded, regular and stable,
with finitely many components and connected total boundary. First rescale the completed sphere to unit curvature. A constant
metric rescaling preserves regularity, stationarity, volume-constrained
stability, reflection symmetry, connectedness and the number of cell
components. Apply Proposition~\ref{ext:stable-structure} and then return
to the original metric, which preserves the resulting Voronoi structure.
Its tangent partitions are minimizing by Lemma~\ref{lem:regularity},
and Lemma~\ref{lem:local-plateau} gives the $3$-Plateau property.
The finitely many spherical faces now give bounded curvature, and
\cite[Section 5.1]{MNP} supplies the face and triple-stratum Lipschitz
conditions. Thus all hypotheses of Proposition~\ref{ext:physical-approximation} hold before its use below.

For fixed $s>0$ the local upper-support argument gives semiconcavity in $v$.
The exact index theorem then yields a negative Alexandrov Hessian, and
hence local concavity, exactly as in Proposition~\ref{prop:local-pdi}.
The exceptional-pressure argument and the bounded conformal
representatives give the sharp profile inequality on the completed sphere.
More explicitly, for fixed $s>0$ the centered intrinsic pressure and the
unit-sphere normalized pressure are respectively
\[
 \lambda_s=LM^{-1}p,\qquad
 \kappa_{\rm unit}=\frac{1}{m\sqrt s}LM^{-1}p.
\]
Here $L:\R^k\to E$ is given by $L\zeta=(\zeta,-\one^T\zeta)$,
and $M=L^TL$.
The map $p\mapsto\kappa_{\rm unit}$ is a linear isomorphism from
$\R^k$ onto $E$. Thus the exceptional set in finite pressure coordinates is
\[
 N_s=\left\{p\in\R^k:\frac{1}{m\sqrt s}LM^{-1}p
                  \in\mathcal B_{n,q}\right\},\qquad |N_s|=0.
\]
Apply Lemma~\ref{lem:null-pullback} to the convex function $u=-J(s,\cdot)$
and the target null set $-N_s$. Since $Du=-p$, the resulting conclusion
is exactly $p\notin N_s$ almost everywhere at this fixed $s$.
All physical variations can be cut off outside $B_{R/2}$ while remaining
unchanged near the entire film, and can then be recentered. For each fixed
$s>0$ and each physical-approximation error, take the small-flow limit
first, and only then let that error tend to zero. No uniform bound on
these individual approximating fields as $s\downarrow0$ is used.

Finally convert to finite-volume coordinates. The tangent-coordinate
inclusion is $L\zeta=(\zeta,-\one^T\zeta)$, with $L^TL=M$. The intrinsic
trace of an inverse Hessian becomes $\tr(MB^{-1})$, whereas the pressure
contraction becomes $p^TB^{-1}p$. Rescaling the unit-sphere inequality to
curvature $s$ gives exactly \eqref{eq:curved-pdi}; Appendix
\ref{app:normalization} checks the factors explicitly. Joint semiconcavity is proved next. On compact parameter cylinders its
Alexandrov Hessian is measurable (it is the density of the absolutely
continuous part of the distributional Hessian). At a joint Alexandrov
point, restricting the quadratic expansion to increments $(0,h)$
identifies the slice's Hessian with the joint $vv$ block. For every fixed
$s>0$, the exceptional volume set for the preceding slice argument is
null. Fubini, applied to this measurable exceptional set and then to a
countable exhaustion by compact cylinders, makes it jointly null.
The unit-sphere statements of Propositions~\ref{ext:physical-approximation}
and \ref{ext:exact-index} are applied after rescaling the completed
sphere to unit curvature; the rescaling of the resulting profile
inequality is computed in Appendix~\ref{app:normalization}.
The claimed inequality therefore holds almost everywhere jointly.
\end{proof}

\section{Joint semiconcavity and curvature lifting}
\label{sec:lift}
The pressure term in the flat profile equation has a rank-one
linearization in the volume Hessian; Appendix~\ref{app:normalization}
records the calculation. To retain ellipticity, we keep the curvature
variable until after linearization. The scaling identity converts the
pressure term into derivatives in that variable. A radial substitution
then turns these derivatives into a Laplacian, with the number of new
variables fixed by the coefficient of $J_s$. Finally we linearize at
the smooth model, so the ellipticity constants involve only the model
and remain controlled at zero curvature.

\subsection{Upper supports in all parameters}
\begin{lemma}[Joint supports and scaling]\label{lem:joint}
The restricted profile $J$ is locally Lipschitz and jointly semiconcave
in $(s,v)$ up to $s=0$. The radial lift
\[
 W(z,v)=J(|z|^2,v),\qquad z\in\R^{n+2},
\]
is jointly semiconcave on a full neighborhood of $(0,v_0)$. At almost
every point with $s>0$,
\begin{equation}\label{eq:euler}
 2sJ_s=nv\cdot J_v-mJ.
\end{equation}
The same identity holds for the smooth standard profile $h$.
\end{lemma}
\begin{proof}
Fix a restricted minimizer $\Lambda$ at $(s_0,v')$. Apply the uniform
volume--moment correction family from Lemma~\ref{lem:curved-attain}
while varying the metric and the target volume. More explicitly, if
$F_a$ is its finite-dimensional flow family, the map
\[
 (s,a)\longmapsto
 \bigl(V_{g_s}(F_a\Lambda),S_s(F_a\Lambda)\bigr)
\]
has uniformly bounded derivatives through order two and an invertible
$a$-derivative with a uniform inverse bound near $(s_0,0)$. Its parameter
inverse gives $a=a(s,v)$ imposing the value $(v,0)$. Consequently
\[
 H_\Lambda(s,v)=P_{g_s}(F_{a(s,v)}\Lambda)
\]
is a $C^2$ function on a parameter neighborhood of uniform radius.
The uniform perimeter bound is used here: after pulling back to the
fixed interfaces of $\Lambda$, each differentiated area Jacobian is
bounded by a constant times the original interface measure. Thus
\[
 J\le H_\Lambda,\qquad H_\Lambda(s_0,v')=J(s_0,v'),
 \qquad |DH_\Lambda|+|D^2H_\Lambda|\le C.
\]
The constants and neighborhood radius are independent of the chosen
restricted minimizer. Support clearance and the positive tube margin
make all corrected clusters admissible after a single reduction of the
parameter neighborhood.

Denote the common parameter radius by $\delta>0$ and, if necessary,
decrease it once so that every corrected competitor also retains the
support clearance and tube margin. Choose a convex half-neighborhood
$D\subset\{s\ge0\}$ of diameter less than $\delta/4$.
The supports based at all points of
$D$ are defined on $D$, majorize $J$ there, and touch it at their base
points. Hence $J=\inf_\Lambda H_\Lambda$ on $D$.
At two points $x,y\in D$, use a support based at $x$ to obtain
$J(y)-J(x)\le C|y-x|$, and then interchange the points. This proves
the Lipschitz estimate. For $x,y\in D$ and $0\le t\le1$, the common Hessian bound gives,
for every support in this same family,
\[
 \begin{split}
 H_\Lambda(tx+(1-t)y)
 &\ge tH_\Lambda(x)+(1-t)H_\Lambda(y)
            -\tfrac C2t(1-t)|x-y|^2\\
 &\ge tJ(x)+(1-t)J(y)-\tfrac C2t(1-t)|x-y|^2.
 \end{split}
\]
Taking the infimum over $\Lambda$ on the left proves the same
inequality for $J$, hence joint semiconcavity on $D$.

This construction also gives semiconcavity across the center of the
radial lift, without introducing a minimization problem at $s<0$.
Let $r=n+2$ and shrink a product neighborhood in $(z,v)$ so that
$(|z|^2,v)\in D$. For $\widehat H_\Lambda(z,v)=H_\Lambda(|z|^2,v)$,
\[
 D^2_{zz}\widehat H_\Lambda
 =2(H_\Lambda)_s I_r+4z\otimes z\,(H_\Lambda)_{ss},\qquad
 D^2_{zv}\widehat H_\Lambda=2z\otimes(H_\Lambda)_{sv}.
\]
Together with the unchanged $vv$ block these formulas give a uniform
upper bound for its full Hessian. The formulas hold at $z=0$ as well;
the smooth correction construction extends each support to small
negative $s$ on the bounded chart, although its majorization property
is needed only at $s\ge0$. Since
$W=\inf_\Lambda\widehat H_\Lambda$, the lift is semiconcave on a full
open neighborhood of $(0,v_0)$.

For the scaling identity fix a differentiability point $(s,v)$ with
$s>0$ and one restricted minimizer $\Lambda$ there. Put
$s_t=e^{-2t}s$ and $\Lambda_t=e^t\Lambda$. Direct substitution in
\eqref{eq:space-form-metric} gives
\[
 (x\mapsto e^t x)^*g_{s_t}=e^{2t}g_s,
 \quad V_{g_{s_t}}(\Lambda_t)=e^{nt}v,
 \quad P_{g_{s_t}}(\Lambda_t)=e^{mt}J(s,v),
\]
and
$S_{s_t}(\Lambda_t)=e^{(n+1)t}S_s(\Lambda)=0$.
Dilation is continuous in $L^1$ on the bounded finite cells; for BV
cells its displacement is bounded by $C|t|P_E(\Lambda)$ on the common
support ball. Thus the support and tube margins make $\Lambda_t$
admissible for both signs of sufficiently small $t$. We have
\[
 J(e^{-2t}s,e^{nt}v)\le e^{mt}J(s,v),
\]
with equality at $t=0$. Differentiating this two-sided touching
inequality gives $-2sJ_s+nv\cdot J_v=mJ$.
Rademacher's theorem proves the identity almost everywhere. Only a
single minimizer was varied; no selection of differentiable minimizers
or global invariance of the restricted tube under scaling is required.
Scaling the standard family gives the same identity for $h$.
\end{proof}

\subsection{Differentiating the scaling identity}
On the open region $s>0$, joint semiconcavity implies that $DJ$ is
locally BV and its distributional derivative is a matrix-valued Radon
measure \cite[Proposition 1.1.3 and Theorem 2.3.1]{Cannarsa}. Its absolutely continuous density is the
Alexandrov Hessian at almost every point. Therefore
\eqref{eq:euler}, an equality between locally integrable functions,
may be differentiated distributionally, using the usual product rule
for multiplication of a measure by the smooth functions $s$ and $v$.
Taking the absolutely continuous parts is legitimate. Slicing the
distributional $vv$ derivatives and applying Fubini shows that their
densities agree with the slice Alexandrov Hessians for almost every
$(s,v)$. Thus Lemma~\ref{lem:curved-structure} and Alexandrov's theorem
provide a measurable set $G\subset D=(0,\eps_0)\times U_0$ with
$|D\setminus G|=0$, on which all identities below and
\eqref{eq:curved-pdi} hold simultaneously. Here $U_0$ is a sufficiently
small open volume ball and $\eps_0>0$.

At such a point write $p=J_v$ and $B=-J_{vv}\succ0$. Differentiation
in $v$ and $s$, respectively, gives
\begin{equation}\label{eq:euler-differentiated}
 p=nBv+2sp_s,
 \qquad nv\cdot p_s=2sJ_{ss}+(m+2)J_s.
\end{equation}
The first formula yields $B^{-1}p=nv+2sB^{-1}p_s$. Using it twice and
then using the second formula and \eqref{eq:euler}, we obtain
\begin{align}
 p^TB^{-1}p
 &=np\cdot v+2ns\,v\cdot p_s+4s^2p_s^TB^{-1}p_s\nonumber\\
 &=mJ+2s(n+2)J_s+4s^2J_{ss}+4s^2p_s^TB^{-1}p_s.
 \label{eq:pressure-elimination}
\end{align}
Substitute this equality into \eqref{eq:curved-pdi}, multiply by $m$,
cancel $mJ$, and divide by $2s$. This division is performed only on
$s>0$. The resulting inequality is
\begin{equation}\label{eq:curvature-time}
 2sJ_{ss}+(n+2)J_s+2sp_s^TB^{-1}p_s
 +\frac{m^2}{4}\tr(MB^{-1})\le0.
\end{equation}

\subsection{The dimension of the radial lift}
Set $r=n+2$, $s=|z|^2$ and $K=W_{zv}$, an $r\times k$ matrix.
On the set
\[
 \widehat G=\{(z,v)\in B_{\sqrt{\eps_0}}^{n+2}\times U_0:
           z\ne0,\ (|z|^2,v)\in G\},
\]
the second-order chain rule gives
\[
 \Delta_zW=2(n+2)J_s+4sJ_{ss},\qquad
 K=2z\otimes p_s,\qquad
 \tr(KB^{-1}K^T)=4s p_s^TB^{-1}p_s,
 \quad B=-W_{vv}.
\]
Thus twice \eqref{eq:curvature-time} is
\begin{equation}\label{eq:lifted-pdi}
 \Delta_zW+\Phi(B,K)\le0,
 \qquad
 \Phi(B,K)=\tr(KB^{-1}K^T)+\frac{m^2}{2}\tr(MB^{-1}).
\end{equation}
The coefficient of $J_s$ determines the number $n+2$ of auxiliary
variables. To justify the almost-everywhere assertion under this change
of variables, polar coordinates give, for every nonnegative measurable
$f$,
\[
 \int_{|z|^2<\eps} f(|z|^2,v)\dd z
 =\frac{|\mathbb S^{n+1}|}{2}
   \int_0^\eps f(s,v)s^{n/2}\dd s.
\]
Fubini therefore shows that the inverse image of any null subset of
$(0,\eps_0)\times U_0$ is null in $(z,v)$. In particular,
$|(B_{\sqrt{\eps_0}}^{n+2}\times U_0)\setminus\widehat G|=0$.
Composition of the second-order
expansion of $J$ with the smooth map $(z,v)\mapsto(|z|^2,v)$ justifies
the displayed chain rule on the full-measure good set. The flat locus
$\{z=0\}$ is itself null, so \eqref{eq:lifted-pdi} holds almost
everywhere on the full lifted neighborhood. The smooth model
$W_*(z,v)=h(|z|^2,v)$ satisfies equality there.

\begin{theorem}[Analytic comparison through zero curvature]\label{thm:lift}
Let $n\ge2$, $m=n-1$, $k\ge1$, and let $J,h$ be defined on
$[0,\eps)\times U$, where $U\subset\R^k$ is open and $v_0\in U$.
Assume the following locally near $(0,v_0)$:
\begin{enumerate}[label=(\roman*)]
\item $J$ is locally Lipschitz and jointly semiconcave, and the lift
$W(z,v)=J(|z|^2,v)$, $z\in\R^{n+2}$, is jointly semiconcave on a full
open neighborhood of $(0,v_0)$;
\item $h$ extends $C^4$ across $s=0$, with $-h_{vv}\succ0$;
\item $J\ge h$ and $J(0,v_0)=h(0,v_0)$;
\item $B=-J_{vv}\succ0$ and \eqref{eq:curved-pdi} hold almost
everywhere for $s>0$, for a fixed symmetric matrix $M\succ0$, and
$h$ satisfies equality in that inequality;
\item $J$ and $h$ satisfy \eqref{eq:euler} almost everywhere for $s>0$.
\end{enumerate}
Then $J=h$ on a smaller product $[0,\eps')\times U'$ containing
$(0,v_0)$, where $U'$ is connected and open.
\end{theorem}
\begin{proof}
The preceding distributional and chain-rule calculations establish
\eqref{eq:lifted-pdi} and its model equality. We now obtain coefficients
whose ellipticity is controlled entirely by the smooth model.

The function $\Phi$ is jointly convex in $(B,K)$ for $B\succ0$.
Indeed, choose a constant matrix $C_0$ with $C_0^TC_0=m^2M/2$ and
stack $K$ and $C_0$ into $T$. Completing the square gives
\[
 \Phi(B,K)=\tr(TB^{-1}T^T)
 =\sup_Y\{2\tr(Y^TT)-\tr(YBY^T)\}.
\]
The supremum is attained at $Y=TB^{-1}$; for fixed $Y$ the expression
is affine in $(B,K)$, proving convexity. Write
\[
 B_*=-W_{*,vv},\qquad K_*=W_{*,zv},\qquad w=W-W_*.
\]
The derivatives of $\Phi$ at the model are
\begin{align*}
 D_B\Phi[H]&=-\tr\left[
 B_*^{-1}(K_*^TK_*+m^2M/2)B_*^{-1}H\right],\\
 D_K\Phi[Z]&=2\tr((K_*B_*^{-1})^TZ).
\end{align*}
Subtract the model equality from \eqref{eq:lifted-pdi} and apply the
supporting-plane inequality for $\Phi$. Since $B-B_*=-w_{vv}$ and
$K-K_*=w_{zv}$, it follows that
\begin{equation}\label{eq:elliptic-lift}
 \mathcal A:D^2_{(z,v)}w\le0\quad\text{almost everywhere},
\end{equation}
where
\begin{equation}\label{eq:lift-matrix}
 \mathcal A=
 \begin{pmatrix}
 I_r&K_*B_*^{-1}\\
 B_*^{-1}K_*^T&B_*^{-1}(K_*^TK_*+m^2M/2)B_*^{-1}
 \end{pmatrix}.
\end{equation}
In particular the two mixed blocks together give the factor $2$ in
the derivative with respect to $K$. For $\eta\in\R^r$ and
$\zeta\in\R^k$,
\begin{equation}\label{eq:schur-square}
 (\eta,\zeta)^T\mathcal A(\eta,\zeta)
 =|\eta+K_*B_*^{-1}\zeta|^2
 +\frac{m^2}{2}\zeta^TB_*^{-1}MB_*^{-1}\zeta.
\end{equation}
If $\zeta\ne0$, then $B_*^{-1}\zeta\ne0$ and the second term equals
$(m^2/2)(B_*^{-1}\zeta)^TM(B_*^{-1}\zeta)>0$.
If $\zeta=0$ and $(\eta,\zeta)\ne0$, the first term is
$|\eta|^2>0$. Thus $\mathcal A\succ0$, using $m>0$ and $M\succ0$.
Choose a compact neighborhood and constants $K_0,b_1<\infty$ such that
\[
 \|K_*B_*^{-1}\|\le K_0,\qquad \|B_*\|\le b_1,
 \qquad c_0=\frac{m^2\lambda_{\min}(M)}{2b_1^2}>0.
\]
For $D=K_*B_*^{-1}$,
\[
 |\eta|^2+|\zeta|^2
 \le2|\eta+D\zeta|^2+(1+2K_0^2)|\zeta|^2.
\]
Combining this with \eqref{eq:schur-square} gives
\[
 \mathcal A\succeq\lambda I,\qquad
 \lambda=\min\{\tfrac12,c_0/(1+2K_0^2)\}>0.
\]
The smooth model also bounds $\|\mathcal A\|$ and its first derivatives
on a smaller compact neighborhood. This gives
$\lambda I\preceq\mathcal A\preceq\Lambda I$ with finite $\Lambda$,
including at $z=0$, where $K_*=0$. All these bounds depend only on the
smooth model and the fixed positive matrix $M$.

It remains to justify passage from the almost-everywhere inequality to
the weak inequality across possible singularities. Semiconcavity gives
the measure decomposition
\[
 D^2w=(D_A^2w)\,\mathcal L^{r+k}+D_{\mathrm{sing}}^2w,
 \qquad D_{\mathrm{sing}}^2w\preceq0.
\]
The sign of the singular part is unchanged when the smooth $W_*$ is
subtracted. Contracting a negative-semidefinite matrix-valued measure
with the continuous positive-definite matrix $\mathcal A$ yields a
nonpositive scalar measure. Combining this observation with
\eqref{eq:elliptic-lift} gives $\mathcal A:D^2w\le0$ in distributions
on the full neighborhood, including $\{z=0\}$. Since $w$ is locally
Lipschitz, the product rule rewrites this as
\[
 \operatorname{div}(\mathcal A Dw)
 -(\operatorname{div}\mathcal A)\cdot Dw\le0.
\]
All hypotheses of Lemma~\ref{lem:weak-minimum} now hold with
$A=\mathcal A$, $b=0$ and $c=0$ on a small connected ball about
$(0,v_0)$. The nonnegative locally semiconcave function $w$ has a zero
at its center, so $w$ vanishes throughout that ball.
Let $\delta>0$ be the radius of this ball. Choose
$0<\eps'<\delta^2/4$ and $U'=B_{\delta/2}(v_0)$, decreasing $\delta$
first if necessary to stay in the original parameter domain. Then
$\{|z|^2<\eps'\}\times U'$ lies in the zero set of $w$. The
definition of the radial lift gives $J=h$ on
$[0,\eps')\times U'$.
\end{proof}

\begin{remark}[Application to the localized profile]
For the profile in \eqref{eq:curved-local-profile},
Lemma~\ref{lem:model} gives the smooth model, Lemmas
\ref{lem:curved-attain}--\ref{lem:curved-structure} give attained
interior minimizers and their positive-curvature inequality, and
Lemma~\ref{lem:joint} gives the semiconcavity and scaling assumptions.
The remaining comparison $J\ge h$, with equality at $(0,v_0)$,
is proved below when Theorem~\ref{thm:lift} is applied.
\end{remark}

\section{Euclidean uniqueness}
\label{sec:euclidean-unique}
\begin{lemma}[The tube contains no standard curved clusters]\label{lem:separation}
For sufficiently small $s\ge0$ and $v$ near $v_0$, no standard $g_s$
cluster satisfying the support and moment conditions in
\eqref{eq:curved-local-profile} lies in its closed tube.
\end{lemma}
\begin{proof}
Suppose otherwise and choose such clusters $\Lambda_j$ with
$(s_j,v_j)\to(0,v_0)$. Their finite cells lie in $B_R$. On this ball
$\rho_{s_j}^n=1+O(s_j)$ uniformly, with a constant depending only on
$n,R$. The weighted volume and zero-moment constraints therefore imply
\[
 w_j:=V_E(\Lambda_j)=v_j+O(s_j)\longrightarrow v_0,
 \qquad
 S_0(\Lambda_j)=O(s_j)\longrightarrow0.
\]
In particular their Euclidean finite volumes remain bounded below by
positive constants. If $s_j>0$, the omitted spherical point is in the
interior of the exterior cell, since that cell contains the complement
of $\overline B_R$. Stereographic projection makes $\Lambda_j$ a
standard Euclidean cluster. This is immediate when $s_j=0$ as well.

The local model in Lemma~\ref{lem:model}, followed by translation to
zero total finite-cell moment, supplies a continuous family $\Theta(w)$
of centered standard Euclidean clusters for $w$ near $v_0$. By
uniqueness within the standard family \cite{Montesinos,MNS,MNP},
\[
 \Lambda_j=R_j\Theta(w_j)+a_j
 \quad\text{up to null sets},\qquad R_j\in O(n).
\]
Taking moments gives
$a_j=S_0(\Lambda_j)/(\sum_i(w_j)_i)\to0$.
Compactness of $O(n)$ supplies a subsequence with $R_j\to R$.
The fixed-domain parametrization in Lemma~\ref{lem:model} gives
$\Theta(w_j)\to\Theta(v_0)$ in $L^1$; continuity of translations and
rotations on finite-measure BV cells then gives
$\Lambda_j\to R\Theta(v_0)\in\cS_0$ in cluster distance.
This contradicts the fixed positive separation between $\cS_0$ and
the closed tube about $\cK$.
\end{proof}

\begin{proof}[Proof of Theorem~\ref{thm:euclidean-main}]
For $n=2$ the assertion is Proposition~\ref{ext:planar-bubbles}, with
$2\le q=k+1\le3$. We may therefore assume $n\ge3$.
Perimeter optimality is Proposition~\ref{prop:flat}. Suppose that a
nonstandard minimizer exists at some positive volume vector $v_0$.
Lemma~\ref{lem:euclidean-compact} supplies a nonempty compact family
$\cK$ separated from the centered standard orbit, and hence the
restricted problem \eqref{eq:curved-local-profile}. Every competitor
for that problem is also an unrestricted competitor. Spherical profile
equality for $s>0$ and Euclidean profile equality at $s=0$ therefore give
\[
 J(s,v)\ge h(s,v),\qquad J(0,v_0)=h(0,v_0).
\]
The equality at the base point follows by choosing a member of $\cK$.
Lemma~\ref{lem:model} supplies the smooth model and its strictly
negative volume Hessian. Lemma~\ref{lem:curved-structure} supplies the
sharp inequality for $s>0$, while Lemma~\ref{lem:joint} supplies joint
semiconcavity, semiconcavity of the lift, and the scaling identity.
Theorem~\ref{thm:lift} thus gives $\eps'>0$ and an open volume
neighborhood $U_0$ of $v_0$ such that $J=h$ on $[0,\eps')\times U_0$.

Choose a nonempty open volume ball $U'$ with
$\overline{U'}\subset U_0\cap(0,\infty)^k$. Put
$M_V=\sup_{v\in U'}\sum_{i=1}^k v_i<\infty$ and fix
\[
 0<s<\min\{\eps',(\omega_n/M_V)^{2/n}\}.
\]
With $a_s=s^{n/2}/\omega_n$, rescaling raw volumes to normalized
unit-sphere volumes sends the vectors $v\in U'$ to
\[
 v\longmapsto\left(
 \frac{s^{n/2}}{\omega_n}v_1,\ldots,
 \frac{s^{n/2}}{\omega_n}v_k,
 1-\frac{s^{n/2}}{\omega_n}\sum_{i=1}^k v_i\right).
\]
Denote this map by $A_s$. Its derivative is $DA_s=a_sL$, of rank $k$,
and its last coordinate is positive by the choice of $s$.
Thus $A_s$ is an affine diffeomorphism onto a relatively open subset of the
interior of the $q$-cell simplex. The dense generic-uniqueness set in
Corollary~\ref{cor:generic} intersects that subset. Choose such a volume
and an attained restricted minimizer $\Lambda$ at its inverse image
$v\in U'$. Writing $w=A_s(v)$, its perimeter satisfies
\[
 P_{g_s}(\Lambda)=J(s,v)=h(s,v)
 =\omega_n s^{-m/2}\Istd(w)=\omega_n s^{-m/2}I(w).
\]
Thus after rescaling it is a global spherical minimizer. Generic
uniqueness makes it standard, contrary to Lemma~\ref{lem:separation}.
Hence $\cK$ is empty, proving uniqueness at every positive Euclidean
volume vector.
\end{proof}

\section{Hyperbolic double bubbles: a force-potential argument}
\label{sec:hyperbolic}
This section is independent of the preceding spherical and Euclidean
arguments. We use unnormalized hyperbolic volume and area and fix sectional
curvature $-1$. Write $\Omega=(\Omega_0,\Omega_1,\Omega_2)$, with $\Omega_0$ the exterior cell and $\Vol(\Omega_i)=v_i\in(0,\infty)$ for $i=1,2$. Put $m=n-1$, and assume throughout this section that $n\ge3$.
The second fundamental form convention is $A=\nabla\nu$; thus an outward
geodesic sphere has positive principal curvatures. Normalized mean curvature
is $\kappa=\tr A/m$, and the stationary pressures satisfy
\[
 p_0=0,\qquad m\kappa_{ij}=p_i-p_j.
\]
All interface sums count each unoriented interface once, after choosing a
normal and retaining that choice in signed formulas.

\subsection{The imported minimizing-cluster structure}
\begin{proposition}[Cotton--Freeman: hyperbolic double-bubble inputs]
\label{ext:cf-structure}
Let $n\ge3$, let $g$ be the metric on $\Hn$ of sectional curvature
$-1$, and let $v_1,v_2>0$. Define the true profile by
\[
 I_{\mathbb H}(v_1,v_2)
 =\inf\{P_g(\Omega):\Omega=(\Omega_0,\Omega_1,\Omega_2),\
              \Vol_g(\Omega_i)=v_i\ (i=1,2)\},
\]
where the infimum is over finite-perimeter clusters and $\Omega_0$
is the exterior cell. In the geometric conclusions below, minimizers are
taken in their regular representatives. The following statements hold.
\begin{enumerate}[label=(\roman*),leftmargin=*]
\item The infimum is attained. Every minimizing double bubble is a
surface of revolution about a geodesic axis
\cite[Propositions 2.3--2.4]{CF}.
\item The function $I_{\mathbb H}$ is strictly concave on
$(0,\infty)^2$. For each minimizer, the exterior cell is connected and
each finite cell has finitely many connected components
\cite[Proposition 4.3, Corollary 4.5, and Proposition 4.11]{CF}.
\item The interfaces of a minimizer consist of finitely many smooth
constant-mean-curvature hypersurfaces of revolution. Their meridian
generators are smooth arcs, meet in threes with $120$-degree sectors
at interior junctions, and meet the axis orthogonally
\cite[Proposition 6.1]{CF}.
\item If a minimizer is nonstandard, it consists topologically of a
sphere with a finite tree of toroidal bands attached. The root has one
branch
\cite[Theorem 6.10]{CF}. The non-root nodes represent physical chamber
components disjoint from the axis; see the Case 1 construction in
\cite[proof of Theorem 6.5, pp.~678--679]{CF}.
\item For the prescribed pair $(v_1,v_2)$, there is exactly one standard
double bubble up to ambient isometry. It consists of three totally
umbilical caps meeting at $120$ degrees, with its outer two caps lying
on geodesic spheres \cite[Proposition 2.6]{CF}.
\end{enumerate}
\end{proposition}

We use the geometric-measure-theoretic existence and regularity
assertion of \cite[Proposition 2.3]{CF} for finite-perimeter minimizers
in their regular representatives. For such a representative, write
\[
 \Sigma=\bigcup_{i=0}^2\partial\Omega_i
       =\Bigl(\bigcup_{0\le i<j\le2}\Sigma_{ij}\Bigr)\cup S_{\rm sing},
 \qquad \HH_g^{n-1}(S_{\rm sing})=0.
\]
The regular interfaces are pairwise disjoint, and hence
\[
 \HH_g^{n-1}(\Sigma)
 =\sum_{0\le i<j\le2}\HH_g^{n-1}(\Sigma_{ij})
 =\frac12\sum_{i=0}^2\Per_g(\Omega_i)=P_g(\Omega).
\]
Thus the area of the bubble in \cite{CF} is exactly the cluster perimeter
used here. The infimum remains over the full finite-perimeter class
specified in Proposition~\ref{ext:cf-structure}.

The compactness of the finite chamber closures is established in
Lemma~\ref{lem:h-pressure}; the closure of the non-root chamber used
below is shown to avoid the axis in Lemma~\ref{lem:h-off-axis}.
The structural statements of Section 6 of \cite{CF} are formulated for
$n\ge3$. No equal-volume assumption, and no bound of two on the number
of components of a finite phase, is used below. The new assertion of the
present section is the exclusion of the nonstandard alternatives for
all positive volume pairs.
The connected rotation group $SO(n-1)$ preserves each complementary
chamber of the rotational film: the orbit of a point is connected and
cannot cross that invariant film. Thus the rotations preserve the labels,
and each meridian edge carries the full angular fiber used below.

The proof separates three questions. The rotational CMC equation first
gives a signed force on every meridian edge. Conservation at a triple
point turns these forces into differences of a scalar potential on
the complementary faces. Geometry bounds that potential from above;
an ambient conformal variation then computes its weighted average
with the opposite variational sign. The pressure and compactness
lemmas are placed before the potential argument to justify all global
minima and integrations by parts used in that comparison.

\subsection{Meridian geometry and the conserved force}
In cylindrical coordinates about the symmetry geodesic,
\begin{equation}\label{eq:h-cylindrical}
 g=d\rho^2+\cosh^2\rho\,dz^2+\sinh^2\rho\,g_{\mathbb S^{m-1}},
 \qquad r=\sinh\rho,\quad c=\cosh\rho,\quad e_z=c^{-1}\partial_z.
\end{equation}
The meridian half-plane has metric $g_B=d\rho^2+c^2dz^2$ and positive
orientation $d\rho\wedge dz$. Orient a chamber boundary with its chamber
on the right and its outward normal on the left. Write
\[
 T=\cos\theta\,e_\rho+\sin\theta\,e_z,
 \qquad \nu=-\sin\theta\,e_\rho+\cos\theta\,e_z.
\]
Here $T$ has unit length in $g_B$ and its parameter is meridian arc length
$s$. Define
\begin{equation}\label{eq:h-psi}
 \psi=c\langle\nu,e_\rho\rangle=-c\sin\theta,
 \qquad \beta=\frac\psi r,
 \qquad \psi\,\dd s=-c^2\,\dd z.
\end{equation}
The meridional principal curvature is denoted by $\alpha$; the other
$m-1$ principal curvatures all equal $\beta$. Consequently
$m\kappa=\alpha+(m-1)\beta$.

\begin{lemma}[Conserved meridian force]\label{lem:h-force}
Along an oriented CMC meridian edge, the quantity
\begin{equation}\label{eq:h-force}
 C=r^{m-1}(\psi-\kappa r)=\frac{r^m}{m}(\beta-\alpha)
\end{equation}
is constant. Reversing the normal and tangent reverses $C$. An edge
whose closure meets the axis has $C=0$.
\end{lemma}
\begin{proof}
Since $r=\sinh\rho$, the meridian connection satisfies
\[
 \nabla_Br=ce_\rho,\qquad \nabla_B^2r=rg_B,
 \qquad \nabla_T\nu=\alpha T,\qquad \nabla_TT=-\alpha\nu.
\]
Here the Hessian identity is an identity on the two-dimensional base:
$\nabla_B^2r(e_\rho,e_\rho)=r$ and
$\nabla_B^2r(e_z,e_z)=c\tanh\rho=r$.
Consequently
\[
 \psi'=\frac{\dd}{\dd s}\langle\nabla_Br,\nu\rangle
 =\nabla_B^2r(T,\nu)+\langle\nabla_Br,\nabla_T\nu\rangle
 =\alpha r'.
\]
The pressure difference, and hence $\kappa$, is constant along the edge.
Differentiating now gives
\[
 \begin{split}
 C'&=r^{m-2}r'\bigl((m-1)(\psi-\kappa r)+r(\alpha-\kappa)\bigr)\\
 &=r^{m-1}r'\bigl((m-1)\beta+\alpha-m\kappa\bigr)=0.
 \end{split}
\]
At an axis contact, $|\psi|\le c$ and $\kappa$ is the finite constant
curvature of that edge. Since $m-1\ge1$,
\[
 |C|\le r^{m-1}(c+|\kappa|r)\longrightarrow0.
\]
This argument does not require choosing a limiting angular coordinate
or limiting normal at the axis.
\end{proof}

\subsection{Strict pressure bounds and compactness}
\begin{lemma}[The actual exterior pressures]\label{lem:h-pressure}
At a minimizing hyperbolic double bubble,
\begin{equation}\label{eq:h-pressure}
 p_1,p_2>m,\qquad \kappa_{10},\kappa_{20}>1.
\end{equation}
Its two finite cells have compact closures.
\end{lemma}
\begin{proof}
Let $I(v_1,v_2)$ be the true minimum-perimeter function.
Lemma~\ref{lem:volume-onto} supplies two compactly supported ambient
fields with independent first finite-volume vectors. Its distributional
proof applies on every connected smooth manifold, without a curvature
or cluster-compactness assumption; only compactly supported volume
increments are used. The implicit
function theorem applied to their composed flows gives exact local
volume correction at the chosen minimizer. Hence, for either sign of small
$t$,
\[
 I(v+te_i)\le I(v)+p_it+O(t^2).
\]
Set $f(t)=I(v+te_i)$ on $(-v_i,\infty)$. Concavity implies
$f'_+(0)\le p_i\le f'_-(0)$: divide the displayed inequality by $t$
and let $t$ tend to zero from the two sides. Every slope between these
one-sided derivatives is a supergradient of a concave function. Thus
\[
 f(t)\le f(0)+p_it\qquad(t>-v_i).
\]
This establishes a global supporting inequality without presupposing
that the two-volume profile is differentiable.

To fix the sign of the field used in the perimeter estimate, use the
upper-half-space model $g=y^{-2}(|\dd x|^2+\dd y^2)$, with
$x\in\R^{n-1}$ and $y>0$, and set $b=-\log y$. Then
\[
 \nabla b=-y\partial_y,\qquad |\nabla b|_g=1,\qquad
 \operatorname{div}_g\nabla b
 =y^n\partial_y(-y^{1-n})=n-1=m.
\]
Let $\chi_R$ be smooth,
equal to one on $B_R$, vanish outside $B_{R+1}$, satisfy
$0\le\chi_R\le1$, and have uniformly bounded gradient. The divergence
formula for a finite-volume finite-perimeter set $U$ gives
\[
 m\int_U\chi_R+\int_U\langle\nabla\chi_R,\nabla b\rangle
 =\int_{\partial^*U}\chi_R\langle\nabla b,\nu_U\rangle
 \le\Per(U).
\]
The second integral tends to zero, since its absolute value is bounded
by a fixed constant times $\Vol(U\setminus B_R)$. Dominated convergence
in the first integral therefore proves $\Per(U)\ge m\Vol(U)$. Apply it to the union of the two
finite cells. The cluster perimeter is at least the perimeter of that
union, so
\[
 I(v+te_i)\ge m(v_1+v_2+t).
\]
After division by $t$ and passage to $t\to\infty$, the supporting
inequality implies $p_i\ge m$. Suppose $p_i=m$. Strict concavity gives
$f(t_0)<f(0)+mt_0$ for each fixed $t_0>0$; equality would force $f$ to
coincide with its supporting line on $[0,t_0]$. Choose $t_0$ and write
$a=(f(t_0)-f(0))/t_0<m$. The secant-slope monotonicity of a concave
function yields $f(t)\le f(0)+at$ for $t\ge t_0$, whereas the preceding
lower bound has asymptotic slope $m$. This contradiction proves
\eqref{eq:h-pressure}.

We include a compactness check before using global meridian minima.
From \eqref{eq:h-force}, $\psi=\kappa r+Cr^{-(m-1)}$. Since
$r'=c\cos\theta$ and $\psi=-c\sin\theta$, one has
$(r')^2=c^2-\psi^2$. Expanding gives
\begin{equation}\label{eq:h-radial-integral}
 (r')^2=1-(\kappa^2-1)r^2-2\kappa Cr^{2-m}
                 -C^2r^{-2(m-1)}.
\end{equation}
On an exterior edge, $\kappa>1$. If $C\ne0$, the right-hand side is
negative both near zero and for large $r$. Hence $r$ remains in a compact
subinterval of $(0,\infty)$. Finite interface area then gives finite
meridian length, since the area element is
$|\mathbb S^{m-1}|r^{m-1}\dd s$. Both $\rho$ and $z$ are bounded on
such an edge: $\rho=\operatorname{arsinh}r$ is bounded and
$|z'|\le c^{-1}\le1$. If $C=0$, then
$\alpha=\beta=\kappa>1$. One can verify compactness directly in the
hyperboloid model: the ambient vector $a=\kappa x-\nu$ is constant on
the hypersurface, has Lorentz norm $\langle a,a\rangle=1-\kappa^2<0$,
and satisfies $\langle x,a\rangle=-\kappa$. The resulting level set is
a geodesic sphere. There are finitely many exterior edges by \cite[Proposition 6.1]{CF}.
Their closures are compact, so choose a ball $B_R$ containing all of them
and put $E=\Omega_1\cup\Omega_2$. The boundary of this union consists
only of exterior interfaces; consequently
\[
 |D\mathbf1_E|(\Hn\setminus\overline B_R)=0.
\]
The complement of $\overline B_R$ is connected, so
$\mathbf1_E=a$ almost everywhere there for some $a\in\{0,1\}$.
Since
\[
 \Vol(E)=v_1+v_2<\infty,\qquad
 \Vol(\Hn\setminus\overline B_R)=\infty,
\]
we have $a=0$. Choose the chambers as connected open components of the
complement of the regular film, with their prescribed labels.
For these open representatives, a nonempty set
$\Omega_i\setminus\overline B_R$ would contain a ball of positive
volume, contradicting $a=0$. Hence
\[
 \Omega_i\subset\overline B_R,\qquad
 \overline{\Omega_i}\subset\overline B_R\quad(i=1,2).
\]
The finite cells therefore have compact closures. Every meridian edge
and finite meridian face also has compact closure, so the minima used
below are attained on these representatives.
\end{proof}

\begin{lemma}[An off-axis chamber in the nonstandard case]
\label{lem:h-off-axis}
If a minimizing hyperbolic double bubble is nonstandard, it has a
finite meridian chamber component $D$ with
$\overline D\Subset\{\rho>0\}$.
\end{lemma}
\begin{proof}
The nonstandard tree in Proposition~\ref{ext:cf-structure}\textup{(iv)} has a non-root physical toroidal
chamber component $U_D$ with $U_D\cap L=\varnothing$, where $L$ is
the rotation axis. Let $D$ be its meridian face.
Lemma~\ref{lem:h-pressure} makes its closure compact. Suppose that
$\overline{U_D}$ met the axis at $x$. If $x$ were interior to a
chamber component, a small connected ball about $x$ meeting $U_D$
would put $x$ in $U_D$, a contradiction. Otherwise Proposition~\ref{ext:cf-structure}\textup{(iii)}
gives a regular orthogonal axis contact. In smooth local coordinates
about $x$, its rotational sheet has the form
\[
 z=z_x+O(\rho^2).
\]
Each of its two local chamber sides therefore contains an open interval
of the axis. A sequence in $U_D$ tending to $x$ belongs to one of these
connected local sides, forcing that side and its axis interval into
$U_D$. This contradicts $U_D\cap L=\varnothing$.
Hence $\overline D$ misses the axis, and
compactness gives $\min_{\overline D}\rho>0$.
\end{proof}

\subsection{A scalar potential on complementary faces}
The complementary faces below are connected meridian chamber components,
not merely the two original labels. Distinct components carrying the same
label are allowed to have distinct potential values.

\begin{lemma}[Planar force potential and axis grounding]\label{lem:h-potential}
There is a number $\eta_D$ on each connected meridian complementary face,
with exterior value $\eta_0=0$, such that
\begin{equation}\label{eq:h-potential}
 C_{D\to E}=\eta_D-\eta_E
\end{equation}
when the normal points from $D$ to $E$. Every face whose closure meets
the axis has potential zero.
\end{lemma}
\begin{proof}
At a triple vertex use the cyclic normal orientations $(01),(12),(20)$.
The normal vectors sum to zero, the pressure differences sum to zero,
and $r,c$ have common values. Formula \eqref{eq:h-force} gives
\[
 C_{01}+C_{12}+C_{20}=0.
\]
Rotating all three normals by the same right angle gives tangents which
are all incoming or all outgoing at that vertex. Thus this equality is
Kirchhoff conservation for the oriented scalar edge flow $C$.

Regard the meridian half-plane as a closed half-plane in $\R^2$.
Give every physical edge its chosen tangent orientation and coefficient
$C_e$. At an axis endpoint the coefficient vanishes by
Lemma~\ref{lem:h-force}; at an interior vertex the preceding identity
makes the sum of incoming coefficients equal the sum of outgoing ones.
Thus the finite weighted network is a compactly supported closed
one-current in the full plane. Closed-loop edges are allowed and do not
change this conclusion.

Here is an elementary construction of its face potential. Declare the
value to be zero in the unbounded complementary region. For a point
off the network choose a polygonal path from that region, transverse to
the edges and avoiding all vertices. The path can be chosen with
finitely many crossings: after a small perturbation it is transverse
to each smooth edge and misses the finitely many vertices; the
intersection set is then compact and discrete. Add $C_e$ when crossing an edge
from its left face to its right face, and subtract $C_e$ in the opposite
direction. This is path-independent. Indeed a closed difference path
can first be resolved into simple closed polygonal curves; summing the
vertex conservation equations inside each such curve shows that its
signed transverse crossing sum is zero. Closed graph loops entirely
inside that curve contribute no crossing term. The resulting potential
is constant on each face and satisfies
$C_e=\eta_{\mathrm{right}(e)}-\eta_{\mathrm{left}(e)}$.
Adding axis segments with coefficient zero subdivides the faces without
changing their potential. This is exactly \eqref{eq:h-potential}.
Equivalently, the current is the boundary of the facewise constant
two-current with coefficient $-\eta_D$ in orientation
$\dd\rho\wedge\dd z$.

Axis grounding can be seen without assigning an angular coordinate at
an endpoint. Every nonzero-force edge has positive distance from the
axis, since $|C|\le r^{m-1}(c+|\kappa|r)$ near it. Finiteness gives a
whole strip along the axis containing no nonzero flow. Its potential is
constant and equals the unbounded exterior value. Faces meeting an axis
segment or endpoint have that same value because all incident forces
there vanish.

The assertion that the exterior has a single face value also deserves
verification. We verify path connectedness after removing the axis.
The Klein projection $(x_0,\mathbf x)\mapsto\mathbf x/x_0$ identifies
$\Hn$ with an open Euclidean ball and sends the axis into a straight
line $\ell$. Join two off-axis points by a finite polygonal path in the
open exterior, subdividing it so that every segment has an internal
vertex as at least one endpoint. Its compact image has positive distance from the
complement of that open set. Perturb its internal vertices by a smaller
amount, keeping the endpoints fixed. For a fixed off-axis endpoint $a$,
\[
 \{b:[a,b]\cap\ell\ne\varnothing\}
 \subset\operatorname{aff}(\{a\}\cup\ell),
 \qquad \dim\operatorname{aff}(\{a\}\cup\ell)\le2<n.
\]
In the finite product of small balls for the internal vertices, first
exclude vertices on $\ell$. Each forbidden edge condition is then
a null set by Fubini and the displayed proper-affine-subspace bound.
Choose vertices outside this finite union. All segments remain in the open
exterior and miss the axis. This proves the required path connectedness.
The same construction applies to every connected open chamber, with
its two chosen off-axis endpoints fixed. Thus removing the axis does
not disconnect any connected open chamber.
The quotient map to $(\rho,z)$ is continuous on the complement of the
axis. Hence the image of the connected physical exterior is connected.
Conversely, the inverse image of a connected open meridian face is
connected because it is the product of that face with the connected
sphere $\mathbb S^{n-2}$. Consequently the meridian faces are precisely
the off-axis images of the connected chamber components. The unique
physical exterior therefore has the single grounded value $0$.
\end{proof}

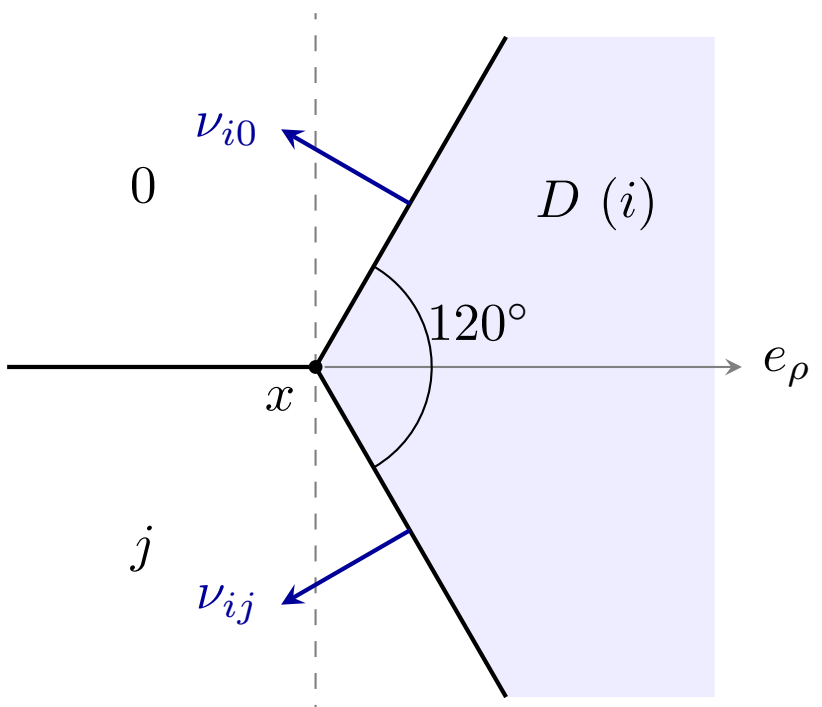
\begin{figure}[htbp]
\centering
\begin{tikzpicture}[x=1.2cm,y=1.2cm,>=stealth,font=\small]
  \fill[blue!7] (0,0)--(1.05,1.819)--(2.2,1.819)--(2.2,-1.819)--(1.05,-1.819)--cycle;
  \draw[thick] (-1.7,0)--(0,0)--(1.05,1.819);
  \draw[thick] (0,0)--(1.05,-1.819);
  \draw[dashed,gray] (0,-1.95)--(0,1.95);
  \draw[->,gray] (.05,0)--(2.35,0) node[right,black] {$e_\rho$};
  \draw[->,thick,blue!60!black] (.52,.90)--(-.19,1.31) node[left] {$\nu_{i0}$};
  \draw[->,thick,blue!60!black] (.52,-.90)--(-.19,-1.31) node[left] {$\nu_{ij}$};
  \draw (.32,-.554) arc[start angle=-60,end angle=60,radius=.64];
  \node at (.9,.25) {$120^\circ$};
  \node at (1.55,.9) {$D\ (i)$};
  \node at (-.95,1.0) {$0$};
  \node at (-.95,-1.0) {$j$};
  \fill (0,0) circle (1.3pt) node[below left] {$x$};
  \node[anchor=north,align=center] at (.65,-2.03)
  {$\mathcal C_x(D)\subset\{v:\langle v,e_\rho\rangle\ge0\}$\qquad
   $\langle\nu_{i0},e_\rho\rangle,\ \langle\nu_{ij},e_\rho\rangle\le-\tfrac12$};
\end{tikzpicture}
\caption{The tangent-plane geometry at a minimum-radius triple point
in Lemma~\ref{lem:h-maximum}. Here $\mathcal C_x(D)$ is the tangent
sector of $D$ at $x$, shown in the symmetric position
$\gamma=0$. The two rays bounding $D$ may tilt together by
$|\gamma|\le\pi/6$. Both outward normals still have negative radial
component. With exactly two finite labels, one neighboring sector is
the exterior, of label $0$; this identifies the edge on which the
positive exterior pressure gives the contradiction.}
\label{fig:minimum-radius}
\end{figure}

\begin{lemma}[Maximum principle for the force potential]\label{lem:h-maximum}
Every finite meridian face satisfies $\eta_D\le0$.
\end{lemma}
\begin{proof}
If some potential were positive, choose a face $D$ with maximal positive
potential. Every outward force from it satisfies $C_{D\to E}\ge0$.
Its closure is disjoint from the axis by
Lemma~\ref{lem:h-potential}. Compactness of that closure gives a
strictly positive minimum of $r$. This minimum lies on the boundary,
since $\nabla_Br=ce_\rho$ never vanishes. The finite regular network
leaves only a smooth-edge point and a triple vertex as possibilities.

Suppose first the minimum is on a smooth edge. The outward normal is
$-e_\rho$, so $\psi=-c$ and $\beta=-c/r$. On the \emph{two-dimensional
meridian metric}, one has $\nabla_B^2r=rg_B$. Hence
\[
 0\le r''=r-\alpha\psi=r+c\alpha,\qquad \alpha\ge-r/c.
\]
Using $c^2-r^2=1$ gives
\begin{equation}\label{eq:h-minimum-sign}
 C=\frac{r^m}{m}(\beta-\alpha)
 \le \frac{r^m}{m}\left(-\frac cr+\frac rc\right)
 =-\frac{r^{m-1}}{mc}<0,
\end{equation}
contrary to the outward-force inequality.

Suppose instead that the minimum is a triple corner; see Figure~\ref{fig:minimum-radius}. In the tangent
plane the $120$-degree sector occupied by $D$ lies in
$\{\dd\rho\ge0\}$. To make the sign explicit, measure tangent angles
from $e_\rho$ and write the two outgoing rays as
$\gamma-\pi/3$ and $\gamma+\pi/3$. Containment in this half-plane means
$-\pi/6\le\gamma\le\pi/6$. The outward normals along these two rays
have radial components
$-\cos(\gamma+\pi/6)$ and $-\cos(\gamma-\pi/6)$, respectively.
Both are at most $-1/2$. Since there are only two finite labels,
one of those two edges separates $D$ from the exterior. Its outward
curvature is positive by Lemma~\ref{lem:h-pressure}, and its force is
nonnegative. But
\[
 \psi=\kappa_{D0}r+C_{D0}r^{-(m-1)}>0,
\]
contradicting the sign of that radial component. The imported regular
network has no other boundary-vertex type away from the axis. The
maximal positive face cannot exist.
\end{proof}

\subsection{One transverse conformal variation}
We compute the constrained second variation by differentiating an
ambient flow. The flow acts on the entire cluster with its junctions.
After pullback to the fixed chambers and interfaces, the volume and
area integrands are smooth Jacobian factors. On their common compact
support these factors and their first two parameter derivatives are
uniformly bounded, so differentiation under the integral is valid.

Use the hyperboloid model with Lorentz metric
$\langle\cdot,\cdot\rangle_L$ of signature $(-,+,\ldots,+)$, and let
$a$ be a unit spacelike vector orthogonal to the two-plane containing the
rotation axis. For $x\in\Hn$, put $\phi(x)=\langle a,x\rangle_L$.
Choosing the first transverse coordinate in the direction $a$ gives
$\phi=rY(\omega)$, where $Y(\omega)=\omega_1$. Since
$\langle x,x\rangle_L=-1$, tangent projection gives
\[
 \nabla\phi=a+\phi x,\qquad
 |\nabla\phi|_g^2=1+\phi^2.
\]
If $D$ is the flat ambient connection and $v\in T_x\Hn$, then
$D_v(a+\phi x)=(v\phi)x+\phi v$. Its tangential part is $\phi v$;
hence
\[
 \nabla_v\nabla\phi=\phi v,\qquad \nabla^2\phi=\phi g.
\]
Its gradient $X=\nabla\phi$ has normal speed
\[
 u=\langle X,\nu\rangle=\psi Y.
\]
The area element on a rotational interface is
$r^{m-1}\dd s\dd\omega$, while $\psi$ depends only on $s$.
Since $\int_{\mathbb S^{m-1}}Y\dd\omega=0$, the flux through every
interface vanishes separately. In particular $\delta_XV_1=\delta_XV_2=0$.
The field is smooth in ambient hyperboloid coordinates even on the axis;
the apparent angular coordinate singularity does not create a boundary
term. Choose a smooth cutoff equal to one on a neighborhood of the
closures of both finite chambers and all their interfaces. Replacing
$X$ by the resulting compactly supported field leaves all derivatives
below unchanged, including the volume integrals over the chamber
interiors.

\begin{lemma}[Direct conformal-flow index]\label{lem:h-conformal-index}
For the preceding field, the actual constrained second variation is
\begin{equation}\label{eq:h-index}
 Q(u)=m\sum_{i<j}\int_{\Sigma_{ij}}u(u-\kappa_{ij}\phi)\,\dd A.
\end{equation}
It is nonnegative at the minimizing cluster.
\end{lemma}
\begin{proof}
Write
\[
 S=\sum_{i<j}\int_{\Sigma_{ij}}\phi^2,\qquad
 U=\sum_{i<j}\int_{\Sigma_{ij}}u^2,\qquad
 H_* =\sum_{i<j}\kappa_{ij}\int_{\Sigma_{ij}}\phi u,
\]
with $P$ the total perimeter. The infinitesimal conformal factor of
$X$ is $\phi$. Differentiating the area Jacobian along its flow gives
\[
 P''=\int_\Sigma\bigl(mX\phi+m^2\phi^2\bigr)
     =mP+mnS.
\]
For a finite chamber the volume Jacobian has first logarithmic
derivative $\operatorname{div}X=n\phi$. Therefore
\[
 V_i''=n\int_{\Omega_i}\bigl(X\phi+n\phi^2\bigr)
       =n\int_{\partial^*\Omega_i}\phi\langle X,\nu_i\rangle.
\]
The last equality is the finite-perimeter divergence theorem for
$\phi X$; every integrand is bounded on the common compact support.
Triple loci have zero hypersurface measure and introduce no extra term
in this volume divergence formula.
Consequently $\sum_i p_iV_i''=mnH_*$.

Apply first-order stationarity to the different ambient field $\phi X$.
Since
\[
 \operatorname{div}_\Sigma(\phi X)
   =|\nabla_\Sigma\phi|^2+m\phi^2=1+n\phi^2-u^2,
\]
we obtain $mH_*=P+nS-U$. Elimination now gives
\[
 \begin{split}
 Q(u)&=P''-\sum_i p_iV_i''
      =mP+mnS-mnH_*\\
     &=nU-nS-P=m(U-H_*),
 \end{split}
\]
which is \eqref{eq:h-index}. The computation uses stationary ambient
variations of the entire cluster, so junction contributions have already
been included. For clarity, the constrained interpretation requires an exact correction,
not merely first-order flux cancellation. Choose compactly supported
smooth fields $Z_1,Z_2$ whose first finite-volume vectors form a basis
of $\R^2$, as supplied by Lemma~\ref{lem:volume-onto}, and let $G_a$
be their composed short flows. The implicit
function theorem produces a $C^2$ parameter $a(t)$ such that
$G_{a(t)}\Phi_t\Omega$ has exactly the original two volumes. Since
$\delta_XV=0$, one has $a(0)=a'(0)=0$ and $a(t)=O(t^2)$.
Differentiating the exact volume identity twice gives
$D_aV\,a''(0)=-V''(0)$. First-order stationarity gives
$D_aP=p\cdot D_aV$. Hence the corrected perimeter has second derivative
$P''(0)+D_aP\,a''(0)=P''(0)-p\cdot V''(0)=Q(u)$.
Its nonnegativity follows from minimality. This argument also confirms
that the ambient-flow calculation has included all junction terms.
\end{proof}

\subsection{Stokes' theorem and exclusion of the nonstandard tree}
Let $\omega_{m-1}=|\mathbb S^{m-1}|$. By invariance under coordinate
permutations, the $m$ integrals of $\omega_j^2$ are equal. Their sum is
$\omega_{m-1}$ because $\sum_{j=1}^m\omega_j^2=1$. Thus
\[
 \int_{\mathbb S^{m-1}}Y^2\,\dd\omega=\frac{\omega_{m-1}}m,
 \qquad
 r^{m-1}u(u-\kappa\phi)
 =r^{m-1}\psi(\psi-\kappa r)Y^2=C\psi Y^2.
\]
The factor $1/m$ cancels the factor $m$ in \eqref{eq:h-index};
integrating over the angular variables therefore gives
\begin{equation}\label{eq:h-edge-index}
 Q(u)=\omega_{m-1}\sum_e C_e\int_e\psi\,\dd s
     =-\omega_{m-1}\sum_e C_e\int_e c^2\,\dd z.
\end{equation}
This is the actual index form, without any additional angular
normalization. Using \eqref{eq:h-potential} and the induced boundary
orientations, the edge sum becomes
\[
 \sum_e C_e\int_e c^2\,\dd z
 =\sum_D\eta_D\oint_{\partial D}c^2\,\dd z.
\]
The exterior coefficient is zero. Any face adjacent to the axis also
has zero potential by Lemma~\ref{lem:h-potential}. Consequently, if an
axis segment is added to complete such a face boundary, its integral is
multiplied by zero; no axis term is omitted from the weighted face sum.
Every face with a nonzero coefficient is a bounded off-axis domain
with finitely many smooth boundary arcs and ordinary corners.
The planar Gauss--Green formula applies to it, including all inner
boundary components; it is also a case of
\cite[Theorem 9.6]{Maggi}. Relative to the ambient orientation
$\dd\rho\wedge\dd z$, our right-face convention is the negative of
the usual induced boundary orientation. Since
$\dd(c^2\dd z)=2rc\,\dd\rho\wedge\dd z$,
\[
 \oint_{\partial D}c^2\,\dd z
 =-2\int_Drc\,\dd\rho\,\dd z
 =-2\int_Dr\,\dd r\,\dd z.
\]
Inner boundary loops have the opposite orientation, so the same formula
holds for a face with holes. For a rectangle $a\le\rho\le b$,
$z_0\le z\le z_1$, the clockwise integral is
$(\cosh^2a-\cosh^2b)(z_1-z_0)<0$, which fixes the orientation explicitly.
Thus
\begin{equation}\label{eq:h-stokes-index}
 \boxed{\quad
 Q(u)=2\omega_{m-1}\sum_D\eta_D\int_Dr\,\dd r\,\dd z\le0.
 \quad}
\end{equation}
The sums here are finite. Each finite face has compact closure and is
a nonempty open subset of $\{r>0\}$, so its weighted area
$\int_D r\dd r\dd z$ is finite and strictly positive. Every coefficient is nonpositive by Lemma~\ref{lem:h-maximum}, whereas
stability gives $Q(u)\ge0$. Thus every individual summand must vanish. Hence
\[
 \eta_D=0\quad\hbox{for every finite face},\qquad
 C_e=0\quad\hbox{on every edge},\qquad \alpha=\beta.
\]
All regular interfaces are therefore umbilical.

\begin{proof}[Proof of Theorem~\ref{thm:hyperbolic-main}]
The previous argument applies to every global minimizer supplied by
\cite{CF}. It remains to exclude nonstandard ways of joining umbilical
pieces; umbilicity alone is not a classification theorem.
If a minimizer were nonstandard, Lemma~\ref{lem:h-off-axis} would supply
an off-axis finite meridian component with compact closure in
$\{\rho>0\}$. All its outward forces now vanish.
At its minimum radius, a smooth boundary point is impossible by
\eqref{eq:h-minimum-sign}. A triple corner is impossible by the
exterior-edge argument in Lemma~\ref{lem:h-maximum}, using $C=0$ and
$\kappa_{D0}>0$. This exhausts the boundary alternatives and is a
contradiction. The minimizer is standard. Existence and uniqueness within
the standard family, \cite[Proposition 2.6]{CF}, give the stated
uniqueness for the prescribed two volumes.
\end{proof}

\begin{remark}[Limits of this argument]
The restriction to two finite phases is geometric, not merely notational.
At an internal triple junction of three finite phases, the face at
minimum radius need not touch the exterior and may have lower pressure
than both neighbors. Then $C\ge0$ does not force
$\kappa r+Cr^{-(m-1)}>0$. Also, the available symmetry of a general
multi-bubble does not supply this two-dimensional meridian graph and its
ordered scalar face potential. Neither Theorem~\ref{thm:hyperbolic-main}
nor the positive-curvature lift proves a hyperbolic three-bubble theorem.
\end{remark}

\appendix
\section{Coordinates, scaling and profile identities}
\label{app:normalization}
\subsection{Intrinsic and finite-volume coordinates}
Let $E=\one^\perp\subset\R^{k+1}$, and define
\[
 L:\R^k\to E,\qquad L\zeta=(\zeta,-\one^T\zeta),
 \qquad M=L^TL=I_k+\one\one^T.
\]
The columns of $O=LM^{-1/2}$ form an orthonormal basis of $E$.
Let $A_E$ be an intrinsic Hessian with $-A_E\succ0$, and put
$\widetilde A=O^TA_EO$ and $H=L^TA_EL$. Then
\[
 \begin{gathered}
 H=M^{1/2}\widetilde A M^{1/2},\qquad
 (-H)^{-1}=M^{-1/2}(-\widetilde A)^{-1}M^{-1/2},\\
 D_vf=L^T\nabla_E f=M^{1/2}O^T\nabla_E f.
 \end{gathered}
\]
Substituting these identities and using cyclic invariance of the trace
gives
\[
 \tr_E((-A_E)^{-1})=\tr(M(-H)^{-1}),
 \qquad
 (\nabla_Ef)^T(-A_E)^{-1}\nabla_Ef
 =(D_vf)^T(-H)^{-1}D_vf.
\]
The factor $M$ in the first trace accounts for the nonorthonormal
finite-volume coordinates.

Dividing both unit-sphere volume and area by $\omega_n$ leaves the pressure
multiplier unchanged. The response operator and perimeter are both divided
by $\omega_n$, which preserves the trace identities. Fix $s>0$ and an
admissible family on the sphere of curvature $s$. Transport this family to
the unit sphere, and denote their raw profiles by $\mathcal J^{(s)}$ and
$\mathcal J^{(1)}$, respectively. Raw volume and area scale by $s^{-n/2}$
and $s^{-m/2}$, so
\[
 \mathcal J^{(s)}(v)=s^{-m/2}\mathcal J^{(1)}(s^{n/2}v).
\]
Here the superscripts label the metric; the parameter $s$ is fixed during
differentiation in $v$. Writing $w=s^{n/2}v$, the chain rule gives
\[
 D_v\mathcal J^{(s)}=s^{1/2}D_w\mathcal J^{(1)},\qquad
 D_v^2\mathcal J^{(s)}=s^{(n+1)/2}D_w^2\mathcal J^{(1)}.
\]
Substitution into the unit-sphere PDI, with
$p=D_v\mathcal J^{(s)}$ and $B=-D_v^2\mathcal J^{(s)}$, yields
\[
 \frac{ms}{2}\tr(MB^{-1})+\frac1m p^TB^{-1}p\le\mathcal J^{(s)}(v),
\]
which is \eqref{eq:curved-pdi}. This rescaling calculation is applied to
the local variational inequality at a fixed $s$; it does not require the
entire admissible tube to scale exactly.

\subsection{The single-bubble model on a two-sphere}
Let $s>0$ and $0<v<4\pi/s$, so that both cells on the two-sphere of
curvature $s$ have positive area. The standard cap boundary has length
\[
 h(s,v)=\sqrt{4\pi v-sv^2}>0.
\]
At $s=0$ the same formula extends to $h(0,v)=\sqrt{4\pi v}$ for
$v>0$, the planar circle profile. Here $m=1$, $M=2$, and
\[
 p=h_v=\frac{2\pi-sv}{h},\qquad B=-h_{vv}=\frac{4\pi^2}{h^3}.
\]
Direct substitution gives
\[
 \frac s2 M B^{-1}+p^2B^{-1}=h,
 \qquad 2sh_s=2vp-h.
\]
The numerator in the first equality reduces to
$sh^2+(2\pi-sv)^2=4\pi^2$. This checks simultaneously the curvature
factor, the finite-coordinate factor and the Euler identity.

On the normalized unit two-sphere the profile is
$i(t)=\sqrt{t(1-t)}$, with $t$ one cell's normalized volume. The intrinsic
orthonormal coordinate is $r=\sqrt2(t-1/2)$, so
$i_r=i'(t)/\sqrt2$ and $i_{rr}=i''(t)/2$. Hence
\[
 (-i_{rr})^{-1}(1+2i_r^2)=2i,
\]
checking \eqref{eq:model-pde} in a case where every quantity is explicit.

\subsection{The flat homogeneous expression is not elliptic}
Let $h_0>0$ be $C^2$ on $(0,\infty)^k$, with
$h_0(tv)=t^{m/n}h_0(v)$ for $t>0$, and set
$p=Dh_0(v)$ and $B=-D^2h_0(v)\succ0$. Euler's identity and its derivative
are
\[
 v\cdot p=\frac mn h_0,\qquad
 D^2h_0\,v+p=\frac mn p.
\]
Since $m=n-1$, these imply
\[
 Bv=p/n,\qquad B^{-1}p=nv,\qquad
 p^TB^{-1}p=n\,p\cdot v=mh_0.
\]
To linearize in the Hessian, hold $p$ fixed and define
$F(p,H)=p^T(-H)^{-1}p/m$ for $-H\succ0$. For each symmetric matrix $Z$,
at $H=D^2h_0(v)$ one has
\[
 D_HF(p,H)[Z]
 =\frac1m(B^{-1}p)^TZ(B^{-1}p)
 =\frac{n^2}{m}v^TZv.
\]
Thus the coefficient matrix is exactly $(n^2/m)v\otimes v$, of rank one.
When $k\ge2$, it annihilates every direction orthogonal to $v$ and is not
positive definite. The uncorrected flat-limit profile equation therefore
does not supply a uniformly elliptic uniqueness argument in that range.
The additional curvature variables in Theorem~\ref{thm:lift} address this
specific loss.

\section*{Acknowledgment of AI assistance}
	The author used ChatGPT in preparing and revising this manuscript and takes responsibility for its mathematical content and references.

\end{document}